\documentclass{article}
\usepackage[dvipsnames]{xcolor}
\usepackage{hyperref}
\usepackage{dcolumn}
\usepackage{bbm,dsfont}
\usepackage{authblk}
\usepackage{scalerel,stackengine}
\usepackage{amssymb,amsmath,amsfonts}
\usepackage{tabularray,tabularx,multirow}
\usepackage[framemethod=tikz]{mdframed} 
\usepackage[sorting=nyvt]{biblatex}
\bibliography{references} 
\usepackage{empheq}
\usepackage{tikz}

\usepackage[normalem]{ulem}

\newcommand\reallywidehat[1]{
   \savestack{\tmpbox}{\stretchto{
       \scaleto{\scalerel*[\widthof{\ensuremath{#1}}]{\kern.1pt\mathchar"0362\kern.1pt}
                  {\rule{0ex}{\textheight}}
               }{\textheight}}{2.4ex}}
   \ensurestackMath{\stackon[-6.9pt]{#1}{\tmpbox}}
}\parskip 1ex

\newcolumntype{Y}{>{\centering\arraybackslash}X}

\begin{document}

\title{Novel spectral methods for shock capturing and the removal of tygers in computational fluid dynamics}
\author[1]{Sai Swetha Venkata Kolluru}
\author[2]{Nicolas Besse}
\author[1]{Rahul Pandit}
\affil[1]{Department of Physics, Indian Institute of Science, 560012 Bangalore, India}
\affil[2]{Laboratoire J.-L. LAGRANGE, Observatoire de la Côte d'Azur, 06304 Nice, France}
\date{\today}
\maketitle

\begin{abstract} 
Spectral methods yield numerical solutions of the Galerkin-truncated versions of nonlinear partial differential equations (PDEs) involved especially in fluid dynamics. In the presence of discontinuities, such as shocks, spectral approximations develop Gibbs oscillations near the discontinuity. This causes the numerical solution to deviate quickly from the true solution. For spectral approximations of the inviscid Burgers equation in one-dimension (1D), nonlinear wave resonances lead to the formation of localised oscillatory structures called \textit{tygers} in well-resolved areas of the flow, far from the shock. 
Recently, the author of~\cite{bessemain} has proposed novel spectral relaxation (SR) and spectral purging (SP) schemes for the removal of tygers and Gibbs oscillations in spectral approximations of nonlinear conservation laws. For the 1D inviscid Burgers equation, it is shown in~\cite{bessemain} that the novel SR and SP approximations of the solution converge strongly in $L^2$ norm to the entropic weak solution, under an appropriate choice of kernels and related parameters.
In this work, we carry out a detailed numerical investigation of SR and SP schemes when applied to the 1D inviscid Burgers equation and report the efficiency of shock capture and the removal of tygers. 
We then extend our study to systems of nonlinear hyperbolic conservation laws --- such as the 2x2 system of the shallow water equations and the standard 3x3 system of 1D compressible Euler equations. For the latter, we generalise the implementation of SR methods to non-periodic problems using Chebyshev polynomials. We then turn to singular flow in the 1D wall approximation of the 3D-axisymmetric wall-bounded incompressible Euler equation~\cite{luo2014potentially}. Here, in order to determine the blowup time of the solution,  we compare the decay of the width of the analyticity strip, obtained from the pure pseudospectral method (PPS), with the improved estimate obtained using the novel spectral relaxation scheme. 
\end{abstract}

\section{Introduction}
\label{sec:intro}
Spectral methods are routinely used for numerical investigations of fluid dynamics PDEs, particularly for studying turbulence in the 3D incompressible Navier--Stokes equations~\cite{canuto2007spectral}. When solutions of these PDEs are smooth, e.g., in elliptic or parabolic problems, it is well-established that spectral and pseudospectral methods exhibit superior accuracy and efficiency of implementation when compared to local (in real space) schemes, such as finite-difference and finite-element methods. This is because the complexity of computing linear gradients is significantly reduced with a spectral basis. Since spectral approximations use global and smooth functions as the approximation basis, they converge to the exact smooth (with analytic or Gevrey regularity) solution with an exponential rate of convergence, also called \textit{spectral convergence}~\cite{canuto1982approximation,kreiss1973methods,tadmor1986exponential}. The spectral solution of the Galerkin-truncated PDE is then, roughly speaking, the classical solution of the full PDE~\cite{gottlieb2001spectral}. 

However, traditional spectral methods are not well-suited for solving nonlinear hyperbolic PDEs, where we frequently encounter discontinuities (shocks) in finite time. When we resolve discontinuous solutions using high-order spectral schemes, the approximations develop Gibbs oscillations near the discontinuity and suffer from non-uniform convergence. The nonlinear terms then cause Gibbs oscillations to mix with and corrupt the solution far away from the discontinuity, resulting in aliasing errors. The situation worsens for nonlinear hyperbolic conservation laws such as the 1D inviscid Burgers and the 1D compressible Euler equations. The absence of physical dissipation in these equations allows for uninhibited growth of Gibbs errors and round-off errors; this makes spectral schemes highly prone to destabilisation~\cite{gottlieb1997gibbs,kreiss1979stability,majda1978fourier}. Even if we mitigate aliasing instabilities, Gibbs oscillations destroy pointwise convergence and render the scheme only first-order accurate. For example, $2/3$-dealiasing prevents the spectral approximation of the 1D inviscid Burgers equation from converging to the exact entropic weak solution with energy-dissipating shocks~\cite{canuto2007spectral}. Stabilisation of spectral schemes is particularly challenging for standard Riemann problems in gas dynamics, modelled by the 1D compressible Euler equations. Discontinuous initial data and the interaction of shocks with fine structures in these systems render spectral schemes virtually unusable~\cite{tadmor1989convergence}. Thus, spectral methods have received limited attention for modelling shock capture in nonlinear hyperbolic conservation laws.  

On the other hand, low-order schemes based on finite-difference and finite-volume methods have had more success. The inherent numerical dissipation of these schemes is carefully tuned to allow for accurate resolution of the shocks~\cite{brehm2015comparison,zhao2019general,moura2016eigensolution}. However, when used to resolve specific flows that contain discontinuities \textit{and} smoothly varying fine structure, these schemes can overdamp the latter. Spectral methods contribute very little or  no inherent numerical dissipation, when used to approximate derivatives of the flow field; this means that, although we have to suppress Gibbs oscillations, we can achieve better control of the dissipation mechanism. Moreover, far from the shock, low-order schemes exhibit reduced accuracy and slow convergence to the smooth solution, whereas spectral methods can achieve much higher-order convergence~\cite{tadmor2007filters}. High-order and high-resolution finite-element schemes, such as the weighted essentially non-oscillatory (WENO) and discontinuous Galerkin (DG) methods are state-of-the-art schemes which perform well for some of the most challenging problems in computational fluid dynamics, such as the interaction of blast waves~\cite{woodward1984numerical,hesthaven2017numerical}. However, a significant limitation of these methods is the complexity of implementation and computational cost~\cite{cai1989essentially,cai1993uniform,kumar2019efficient,klockner2011viscous,hartmann2006adaptive}. Thus, we would like to extend the use of spectral methods to shock-capture problems, as they have the potential to provide high-order approximations with a comparatively straightforward and low-cost implementation~\cite{gottlieb1981spectral,majda1978fourier,hussaini1983spectral}.\par

To model nonlinear hyperbolic conservation laws using spectral methods, we should control Gibbs oscillations near the discontinuity and thereby recover spectral or high-order convergence, at least in the piecewise smooth areas of the flow. Early attempts to control the high-frequency growth used post-processing methods based on filtering, wherein a high-order spectral filter was applied to improve the accuracy of the solution away from the discontinuity~\cite{tadmor2005adaptive,vandeven1991family,sun2002practical,hussaini1985spectral,tadmor2007filters,gottlieb1985recovering}. Filtering techniques have since been adapted to applications in the physical domain and for systems of nonlinear conservation laws~\cite{tadmor2007filters,sun2006windowed}. Since solutions to different problems have different spectral distributions, the filter chosen for optimal shock capture in one problem might perform poorly in another. Then, it is desirable to have only a few parameters for the filtering scheme that can be tuned effectively across implementations. Time evolution of the solution presents another challenge as the filtering needs to be updated dynamically to apply just enough dissipation for accurate shock capture. For non-periodic problems, applying filtering at the boundaries can be challenging, especially for the 1D compressible Euler equations, where we require characteristic-compatibility (CCM) methods~\cite{canuto2006spectral}. Here we note that in hyperbolic problems, any error at the boundary travels into the domain with incoming waves, unlike elliptic PDEs, where the boundary errors remain localised at the boundary. The lack of physical dissipation and the conservative nature of the spectral scheme can quickly amplify these errors and cause temporal stability to be lost~\cite{canuto2007spectral,thompson1990time}. 

Since spectral schemes for nonlinear hyperbolic conservation laws conserve energy, introducing a small amount of numerical (vanishing) viscosity in the schemes is a natural mechanism to obtain entropic weak solutions that dissipate the energy. Tadmor~\cite{tadmor1989convergence} developed the spectral vanishing viscosity (SVV) method by adding dissipation to the PDE through an elliptic operator of second order, on the right-hand side of the PDE, and a viscous regularisation kernel for building the numerical diffusion coefficient. The kernel was chosen to provide optimal damping of higher wavenumber modes of the solution, leaving low wavenumber modes unchanged. The stability and convergence analysis is more tractable for the SVV method than for filtering schemes. For the 1D inviscid Burgers equation,~\cite{tadmor1989convergence} established convergence of the SVV approximation to the entropic weak solution that dissipates energy. The kernel choice for the viscous term was improved by~\cite{maday1989analysis,tadmor1990shock}, and the convergence analysis (under the assumption that the approximate solution is uniformly bounded in $L^\infty$) was extended to hyperviscous kernels by~\cite{tadmor1993super}. The convergence analysis of the SVV method is yet to be established mathematically for systems of nonlinear hyperbolic conservation laws, but numerical explorations  have been carried out by~\cite{tadmor1989convergence}. The modification for non-periodic problems, using Legendre polynomials, was initiated in~\cite{maday1993legendre,guo2001spectral}. Despite viscous regularisation, the high-order SVV approximations still exhibit small oscillations with a bounded amplitude near the shock~\cite{karamanos2000spectral}. We present a brief note on the implementation of the SVV method for the 1D inviscid Burgers equation in Appendix~\ref{app:filters}. 

Recently, ~\cite{ray2011resonance} investigated the 2/3-dealiased spectral approximation of the 1D-periodic inviscid Burgers equation, starting from harmonic initial data. The authors reported spontaneous growth of localised oscillatory structures called \textit{tygers} away from the shock in smooth regions of the flow. In the Fourier spectrum, tygers manifest as a `boundary layer' close to the time of singularity when spectral convergence is lost. These structures collapse and spread throughout the domain and serve as precursors to thermalisation in the energy-conserving Galerkin-truncated Burgers equation~\cite{venkataraman2017onset}. Tygers have also been reported in the case of the 2D incompressible  Euler equation, where there is no finite-time singularity, but the flow develops ever-decreasing length scales~\cite{ray2011resonance,sulem1983tracing}. Further discussion on tygers and thermalisation in Galerkin-truncated PDEs can be found in~\cite{di2018dynamics,ray2015thermalized,frisch2008hyperviscosity,krstulovic2009cascades,cartes2021galerkin,feng2017thermalized}. Subsequently, ~\cite{murugan2020suppressing} developed the first \textit{tyger-purging} (TP) scheme for the 1D inviscid Burgers equation. Here, dissipation is provided by selectively purging Fourier modes of the velocity field that lie above the cut-off purging wavenumber $K_P$. The purging step is performed at discrete times separated by a purging timescale $\tau_P$. The parameters controlling the dissipation rate, $K_P$ and $\tau_P$, are tuned to match the energy decay of the purged spectral approximation with that of the actual entropic solution. The \textit{a priori} knowledge of the existence and energy decay of an entropic solution is crucial for the successful implementation of this scheme. Thus, extensions to systems of nonlinear conservation laws and to multidimensional nonlinear PDEs are non-trivial. Furthermore, the stability and convergence of the tyger-purging scheme are yet to be investigated mathematically. The removal of tygers, by using wavelet analysis, has been investigated numerically in~\cite{pereira2013wavelet,farge2017wavelet,pereira2023adaptive}.\par

There has been a recent renewal of interest in controlling the high-frequency growth in spectral approximations of fluid dynamics equations; this is partly because of the limiting role of tygers and Gibbs oscillations in the spectral search for finite-time singularities. Historically, numerical investigations of potentially singular flows have proved to be delicate and challenging tasks~\cite{gibbon2008euler}. The numerical schemes have to resolve accurately the ever-decreasing flow scales, while operating within the computational constraints of fixed precision and finite resolution. Spectral methods have been used to study the finite-time-blowup problem in many systems: the 1D-periodic inviscid Burgers equation~\cite{sulem1983tracing,caflisch2015complex,frisch2003singularities}, the 3D incompressible Euler equations in a periodic cube~\cite{cichowlas2005evolution,kida1986study}, the 3D incompressible, ideal magneto-hydrodynamic equations in a periodic cube~\cite{brachet2013ideal}, and more recently, to investigate a potentially singular solution of the 3D-axisymmetric wall-bounded incompressible Euler equation~\cite{kolluru2022insights,barkley2020fluid}. Spectral methods achieve reliable singularity tracking through the analyticity-strip technique~\cite{sulem1983tracing,caflisch2015complex}, in which the rate of exponential decay of Fourier spectra can be used to deduce the distance $\delta$ of the nearest complex singularity from the real-axis~\cite{carrier2005functions} (for extensions to Chebyshev spectral methods see~\cite{trefethen2019approximation,kolluru2022insights}). This decay rate yields the width $\delta(t)$ of the analyticity strip; if $\delta(t)$ vanishes at $t=t_*$, the complex singularity manifests itself as a finite-time singularity. For the 1D inviscid Burgers equation, spectral estimates of $\delta(t)$ exhibit power-law behaviour with an intercept on the $t$ axis at $t_*$, the time of singularity (or blowup time). However, we can only use the analyticity strip method as long as we have spectral convergence of the solution. After this time,  we cannot extract $\delta(t)$ reliably as the emergence of tygers and subsequent thermalisation contaminate the Fourier spectrum with the back-flow from higher wavenumbers. Recent investigations~\cite{kolluru2022insights} of the potential singularity in the 3D-axisymmetric wall-bounded incompressible Euler equations report tygers in the flow variables. In particular, the authors of
Ref.~\cite{kolluru2022insights} use a $2/3$-dealiased pseudospectral scheme to control the aliasing instability. Dealiasing leads to a significant reduction in the number of usable Fourier modes, which limits the time for which the potential singularity can be tracked. The addition of controlled dissipation to the spectral scheme can help to mitigate the aliasing instability, while avoiding dealiasing; this can extend the time until which $\delta$ can be extracted reliably. This dissipative mechanism, which can resemble the spectral shock-capturing schemes discussed earlier, might be able to control Gibbs oscillations and tygers, and to delay thermalisation. \par

Recently, Besse~\cite{bessemain} has proposed novel spectral schemes for approximating the solution of nonlinear conservation laws --- namely, the spectral relaxation (SR) and spectral purging (SP) methods. For the 1D inviscid Burgers equation, ~\cite{bessemain} presents a mathematical proof for convergence in the $L^2$ norm, of the SR and SP approximations to the unique entropic weak solution under a suitable choice of kernels and related parameters. In the SR scheme, dissipation is built into the PDE through a  Bhatnagar-Gross-Krook (BGK)-type relaxation operator; the kernel used for the relaxation mechanism damps the high wavenumber modes of the solution. The SR scheme is similar to the SVV method, in that the dissipation is built into the PDE. The dissipative mechanism, however, is very different; SVV methods use a viscous term, while SR methods use a relaxation-type mechanism. Indeed, in the SVV method, the monotonicity property and the maximum principle of the associated solution operator (on which $L^\infty$ bounds and Krushkov BV theory are built) are obtained thanks to a dissipative operator given by the Laplacian (or viscosity) regularisation. Unlike the viscosity regularisation, the solution operator associated with SR and SP methods (using a BGK-type relaxation operator for the dissipative operator on the right-hand of the PDE) lack the monotonicity property and the maximum principle. Let us note in passing that the loss of these properties is also present in hyperviscosity regularisations, for which a proof of an $L^\infty$ bound for the approximate solution, and thus a complete convergence analysis, are still open problems. In SR and SP methods, the monotonicity property and the maximum principle are recovered, not by the special structure of the dissipative operator (like the Laplacian in the SVV method) but by the property of the convolution kernels involved in the definition of the relaxation operator (especially the positivity of kernels). Moreover, the criteria for the selection of convolution kernels are also very different for SR and SVV schemes. In the SP method, the spectral approximation is convolved with a regularising kernel at discrete times; this purging step damps the high-wavenumber modes and removes Gibbs oscillations and tygers in the approximate solution. The SP scheme is the mathematically rigorous counterpart of the tyger-purging (TP) scheme~\cite{murugan2020suppressing}; the former can also be extended to problems where we do not know \textit{a priori} the entropic solution. 

In this work, we conduct rigorous numerical studies of the SR and SP approximations applied to nonlinear hyperbolic PDEs. For the 1D inviscid Burgers equation, we perform an extensive comparative study of SR and SP schemes and the pure pseudospectral scheme (PPS) at various times related to the formation and development of shocks. We obtain the SR and SP approximations by using different kernels; in particular, we see good visual convergence with the high-order de La Vall\'ee Poussin kernel which is not positive. For the 1D inviscid Burgers equation, we perform numerical convergence analysis of the SR approximation using $L^1$ and $L^2$ errors. We then compare the shock capture achieved with SR and SP schemes and the standard SVV scheme. We numerically extend the application of SR schemes to $2\times2$ and $3\times3$ systems of nonlinear hyperbolic conservation laws --- namely, the 1D shallow water equations and the 1D compressible Euler equations, respectively. Here, we conduct exhaustive tests with multiple standard problems per model. Finally, we use the SR method to improve the tracking of the finite-time singularity for the 1D wall approximation~\cite{1d_1,choi2017finite} of the 3D-axisymmetric wall-bounded incompressible Euler equation~\cite{luo2014potentially}.\par

The remainder of this article is organised as follows. In Section~\ref{sec:methods}, we describe in detail the implementation of pure pseudospectral (PPS), spectral relaxation (SR), and spectral purging (SP) methods. We also outline the extension of these methods to systems of nonlinear hyperbolic conservation laws. 
In Section~\ref{subsec:Results_burg}, we present results for the 1D inviscid Burgers equation. Here, we assess the shock capture and the tygers removal of SR and SP schemes and, furthermore, we examine quantitatively the numerical convergence of the novel spectral approximations to the exact solution. 
In Section~\ref{subsec:Results_burg_ics}, we present an exhaustive report on how the choice of initial conditions can affect the choice of optimal parameters for non-dealiased spectral schemes.
In Section~\ref{subsec:Results_shalwat}, we use SR schemes to model the solution of the $2\times2$ system of 1D shallow water equations. We perform tests using two initial conditions for this model --- namely, the hump of water and the dam break. In Section~\ref{subsec:Results_Euler}, we present analogous results using Chebyshev spectral relaxation schemes for the $3\times3$ system of 1D compressible Euler equations on a non-periodic domain. Here, we use an array of challenging Riemann problems to test the convergence and stability of SR approximations. We pay particular attention to the implementation of boundary conditions using the characteristic compatibility methods (CCM). In Section~\ref{subsec:Results_lhmodel}, we discuss SR methods for singular flows in the 1D wall approximation of the 3D-axisymmetric wall-bounded incompressible Euler equation. We compare the time-decay of the analyticity-strip width $\delta(t)$ obtained using PPS and SR schemes. Finally, we summarise our findings and draw attention to the scope and future directions of our study in Section~\ref{sec:conclusions}.\par

There are 2 appendices to this article. In Appendix~\ref{app:filters}, we profile the kernels employed in SR and SP schemes. In Appendix~\ref{app:cbc}, we give an accessible discussion of the characteristic boundary conditions for the 1D compressible Euler equations. For the 1D shallow water and 1D compressible Euler equations, the reference solutions are derived from well-tested and standard open-source libraries such as  \texttt{CLAWPACK}~\cite{clawpack} (for high-resolution finite-volume approximations) and \texttt{HyPar}~\cite{hypar} (for high-order finite-difference  methods).\par

\section{Numerical methods}
\label{sec:methods}

In this section, we discuss the pure pseudospectral method and its modifications, namely, the spectral relaxation and spectral purging methods, in that order. First, we describe the implementation of these methods for the 1D inviscid Burgers equation with periodic boundary conditions, as formulated by~\cite{bessemain}. We then describe our notation for the extension of these schemes to systems of nonlinear conservation laws and non-periodic domains.\par

\subsection{The pure pseudospectral method (PPS)}
\label{subsec:Methods_pure_spec}
Consider the 1D inviscid Burgers equation written in the following conservative form,
\begin{subequations}    
\begin{align}
    \partial_t u(x,t) + \partial_x \left( \tfrac{1}{2}{u^2(x,t)} \right) = 0, \quad t \in \mathbb{R}^{+}, \quad x \in \mathbb{R}, \label{eq:burg_x_pps}
\end{align}
where $u(x,t)$ is the velocity field for which the initial data is given by $u(x,t=0)=u_0(x)$. In an $L$-periodic domain, we use Fourier pseudospectral methods to solve Eq.~\eqref{eq:burg_x_pps} over a uniform collocation grid of $N_x=2N+1$ points given by $\mathbf{X}_{L}=\{x_j = j\Delta x$ : $j=0, 1 ,..., 2N\}$, where $\Delta x=\tfrac{L}{2N+1}$. The Fourier-pseudospectral projection operator $P_N$ acts on $u(x,t)$ to give the pseudospectral approximation $u_N(x,t)$ as follows,
\begin{align}
    u_N(x,t) \coloneqq {P}_N u(x,t)= \sum_{|k|\leq N} \widehat{u}(k,t) e^{i\frac{2\pi}{L} kx}. \label{eq:pn_coeffs} 
\end{align}
Thus, the exact solution $u(x,t)$ is approximated by its projection on a finite basis of trigonometric polynomials with wavenumber $k \in [-N,N] \subset \mathbb{Z}$. The discrete Fourier coefficient $\widehat{u}(k,t)$ is computed as a quadrature approximation on the collocation grid $\mathbf{X}_{L}$. In our calculations, we use the \texttt{FFTW3}~\cite{frigo1999fftw} library for this step:
\begin{align}
    \widehat{u}(k,t) = \frac{1}{2N+1} \sum_{j=0}^{2N} u(x_j,t) e^{-i\tfrac{2\pi}{L}kx_j}.     \label{eq:fourcoeffs}
\end{align}
Equation~\eqref{eq:burg_x_pps} can be rewritten for the pseudospectral interpolant $u_N(x,t)$ as
\begin{align}
    \partial_t u_N +  \partial_x (P_N \tfrac{1}{2} u_N^2) =0, \label{eq:burg_k_pps}
\end{align}
where the initial data is given by $u_N(x,0)=P_N u_0(x)$.
Equation~\eqref{eq:burg_k_pps} is the Galerkin-truncated version of the full 1D Burgers equation in Eq.~\eqref{eq:burg_x_pps}.
\label{eq:burg_pps}
\end{subequations}

In the pure pseudospectral scheme described above, we use $2/3$-dealiasing to mitigate the growth of aliasing errors which can otherwise cause the scheme to destabilise quickly. The Fourier series in Eq.~\eqref{eq:pn_coeffs} is truncated up to the dealiasing cut-off wavenumber $K_G = 2N/3$ before computing the nonlinear term $P_N\tfrac{1}{2}u_N^2$ in real space, at every time step. 

We now describe the PPS scheme for the $d \times d$ system of 1D nonlinear hyperbolic conservation laws governing the $d$-component vector of conserved quantities $\mathbf{q}(x,t) \in \mathbb{R}^d$. The equations are given in the following conservative form,
\begin{subequations}
\begin{align}
    \partial_t \mathbf{q} + \partial_x \mathbf{F}(\mathbf{q}(x,t)) =0;\label{eq:sys_x_cons}
\end{align}
here $\mathbf{F}(\mathbf{q})\in \mathbb{R}^{d}$ is the flux vector for the quantity $\mathbf{q}$. In the advective formulation, Eq.~\eqref{eq:sys_x_cons} can be rewritten as follows,
\begin{align}
    \partial_t \mathbf{q} +  A(\mathbf{q}) \ \partial_x \mathbf{q} =0, \label{eq:sys_x}
\end{align}
where $A(\mathbf{q})=\nabla_{\mathbf{q}} \mathbf{F}$ is the $d\times d$ Jacobian matrix. Equations~\eqref{eq:sys_x_cons} and \eqref{eq:sys_x} are \textit{hyperbolic} at a point $(\textbf{q},x,t)$, if the matrix $A(\mathbf{q})$ is diagonalisable with real eigenvalues, at this point. 

In $L$-periodic domains, the pseudospectral Fourier projection interpolates the component variables of $\mathbf{q}(x,t)$, i.e., $\{ q^{(i)}; \ i=0, 1,..., (d-1)\}$ over a uniform grid of $2N+1$ points given by $\mathbf{X}_L=\{x_j= j\Delta x : j=0, 1,..., 2N\}$ where $\Delta x=\tfrac{L}{2N+1}$ as before. We have the analogues of Eqs.~\eqref{eq:pn_coeffs} and~\eqref{eq:fourcoeffs} for the components $q^{(i)}$ as follows:
\begin{align}
    q^{(i)}_N(x,t) \coloneqq {P}_N q^{(i)}(x,t)= \sum_{|k|\leq N} \widehat{q}^{(i)}(k,t) e^{i\frac{2\pi}{L}kx} , \\ \widehat{q}^{(i)}(k,t) = \frac{1}{2N+1} \sum_{j=0}^{2N} q^{(i)}(x_j,t) e^{-i\frac{2\pi}{L}kx_j}.
\end{align}
The system of conservation laws in Eq.~\eqref{eq:sys_x} is now rewritten for the pseudospectral projection $\mathbf{q}_N$ as
\begin{align}
\partial_t \mathbf{q}_N + P_N( A(\mathbf{q}_N) \ \partial_x \mathbf{q}_N )= 0.  \label{eq:sys_k} 
\end{align}
\label{eq:sys}
\end{subequations}
We use Chebyshev pseudospectral methods to solve the 1D compressible Euler flows described in Section~\ref{subsec:Results_Euler} for non-periodic domains. Boundary conditions (BCs) for these flows are problem-specific and differ for the left and right boundaries. Direct implementation of the physical boundary conditions can quickly destabilise the spectral schemes here. We then use characteristic boundary conditions (CBC) which take into account the eigen-structure of the Jacobian matrix. We detail their implementation by CCM methods in Appendix~\ref{app:cbc}. However, once we implement the BCs in a mathematically consistent way, we do not need local treatment at the boundaries, even for the novel spectral schemes discussed below. Further details are provided in Section~\ref{subsec:Results_Euler} and Appendix~\ref{app:cbc}.  

\subsection{The Spectral Relaxation method (SR)}
\label{subsec:Methods_spec_rel}
In the spectral relaxation scheme of~\cite{bessemain}, the pseudospectral approximation of the 1D inviscid Burgers equation in Eq.~\eqref{eq:burg_k_pps} is modified by adding a relaxation term to the right-hand side (RHS), as written below,
\begin{subequations}
\begin{align}
    \partial_t u_N + \partial_x (P_N \tfrac{1}{2} u_N^2) & = \frac{1}{\tau} (\mathcal{M}_m u_N - u_N).
    \label{eq:burg_x_sr}
\end{align}
The convolution operator defined as $\mathcal{M}_m u_N (x,t) \coloneqq K_m(x) * u_N(x,t) $ performs singular convolution of the approximate solution $u_N$ with the chosen kernel $K_m(x)$; this step produces a smoother approximation of $u_N$ by damping high-$k$ modes. 
Then, the modified RHS of Eq.~\eqref{eq:burg_x_sr} acts as a Bhatnagar--Gross--Krook (BGK) relaxation mechanism; it dissipates energy from the high-$k$ modes of the pseudospectral approximation $u_N$ at an exponential decay rate $1/\tau$ and forces $u_N$ to converge to the generalised ``Maxwellian" $\mathcal{M}_m u_N$. The parameter $m$ and the relaxation time scale $\tau$ are controlled by the parameters $(\alpha, \gamma)$ as defined by~\cite{bessemain} and given by
\begin{align}
        m  = N^\gamma ,  \quad 0 < \gamma <1 \qquad \text{and} \qquad \tau  = N^{-\alpha},  \quad {\alpha >0}.
        \label{eq:sr_params}
\end{align} 
Since the convolution operation on the RHS of Eq.~\eqref{eq:burg_x_sr} simplifies to a product in Fourier spectral space, we rewrite the equation in Fourier space as follows,
\begin{align}
    \partial_t \widehat{u}(k,t) + i \left( \frac{2\pi}{L}k \right) \reallywidehat{P_N (\tfrac{1}{2}u_N^2)}(k,t) = \frac{1}{\tau} \left[ \widehat{K}_m (k) -1 \right] \widehat{u}(k,t), \label{eq:burg_k_sr}
\end{align}
where $\reallywidehat{P_N (\tfrac{1}{2}u_N^2)}(k,t)$ corresponds to the discrete Fourier coefficient of the nonlinear term as defined in Eq.~\eqref{eq:fourcoeffs} and $\widehat{K}_m(k)$ is the discrete Fourier coefficient of the kernel $K_m(x)$. Initial data for Eq.~\eqref{eq:burg_x_sr} is given by $u_N(x,0)=P_N u_0(x)$, the same as in the PPS method. 
\label{eq:burg_sr}
\end{subequations}

In the present work, we use, among other kernels, the positive kernels of Fej\'er--Korovkin, Jackson and Jackson--de La Vall\'ee Poussin (see Appendix~\ref{app:filters} and also \cite{butzer1971fourier} for a detailed discussion of all filters used here) for which it has been proved in \cite{bessemain} that the approximate solution  $u_N$ converges strongly in $L^2$, as $N\rightarrow \infty$, to the entropic weak solution of the 1D inviscid Burgers equation. 

For the system of nonlinear hyperbolic conservation laws given in Eqs.~\eqref{eq:sys}, the relaxation operation is applied to the RHS of all the component equations. The chosen kernel $K_m(x)$ is used with the same parameter set $(\alpha,\gamma)$ for all the component equations. SR schemes can be extended in a straight-forward way to Chebyshev-pseudospectral approximations for non-periodic domains. We present the results of Chebyshev SR approximations for the 1D compressible Euler equations in Section~\ref{subsec:Results_Euler}; here, we note that application of the relaxation kernel at the sensitive boundary points goes through smoothly within the spectral relaxation scheme. 

Furthermore, we need not use $2/3$-dealiasing or other dealiasing strategies for the spectral relaxation and spectral purging schemes. For all the test problems discussed in this article, a suitable choice of the parameters $(\alpha,\gamma)$ suffices to control aliasing errors and prevent destabilisation of the new schemes. We explore the effect of dealiasing on SR and SP schemes in Section~\ref{subsec:Results_burg_ics} in detail.

\subsection{The Spectral Purging method (SP)}
\label{subsec:Methods_spec_purg}
In the spectral purging schemes of~\cite{bessemain}, the pseudospectral approximation $u_N$ is mollified with the kernel $K_m(x)$ at discrete times, separated by a time interval $\tau$. The definitions of the parameters $m$ and $\tau$ for the spectral purging method are the same as described for the spectral relaxation method in Eq.~\eqref{eq:sr_params}; we continue our usage of the parameters $(\alpha,\gamma)$ for SP methods.
Consider the discretisation of the time interval $t \in [0,T]$ into intervals of length $\tau$, using the index $n=0, 1 ,..., M$, where $\tau=T/(M+1)$. We use the notation $\mathbbm{1}_{]n\tau,(n+1)\tau]} (t)$ for the indicator function over the subset ${t \in ]n\tau,(n+1)\tau]}$. Then, the spectral purging approximation is given by 
\begin{subequations}    
\begin{align}
    u_N(x,t) = \sum_{n=0}^{M} u^n_N(x,t) \mathbbm{1}_{]n\tau,(n+1)\tau]} (t).  
\end{align}
Here, $u^n_N(x,t)$ is the exact solution of the following spectral purging equation:
\begin{align}
    &\partial_t u^n_N + \partial_x  (P_N \tfrac{1}{2}(u^n_N)^2) = 0, \qquad \qquad \forall t \in ]n \tau, (n+1)\tau] \label{eq:burg_x_sp}, \\
    &u^n_N(x,n\tau) = \mathcal{M}_m u^{n-1}_N (x,n\tau). \label{eq:burg_k_sp}
\end{align}
The convolution operator $\mathcal{M}_m ( \cdot )=K_m * ( \cdot )$ performs the singular convolution integral discussed in the previous section. Here, $u^{-1}_N(x,t=0) \coloneqq u_N(x,t=0)$ is the pseudospectral approximation of the initial condition $u_0(x)$.  The pseudospectral approximation of $u^n_N(x,t)$ is given by
\begin{align}
    u^n_N(x,t) = \sum_{|k|\le N} \widehat{u}^n(k,t) e^{i\frac{2\pi}{L}kx}, \qquad \forall t \in  ]n \tau, (n+1) \tau],
\end{align}
\end{subequations}
where the discrete Fourier coefficients $\widehat{u}^n(k,t)$ are computed as in Eq.~\eqref{eq:fourcoeffs}. Thus, the dissipative mechanism in the SP method is discontinuous-in-time unlike in the SR method, where the dissipation is enforced continuously-in-time via the penalty term on the right-hand side of Eq.~\eqref{eq:burg_x_sr}. We do not use $2/3$-dealiasing strategies for the SP method, as elaborated in the previous section for SR methods.

For systems of nonlinear conservation laws, we apply the purging operation to each component variable $q^{(i)} \in \mathbf{q}(x,t)$ of Eq.~\eqref{eq:sys}. The choice of kernel $K_m(x)$ and parameters $(\alpha,\gamma)$ remain the same for purging all the component variables of $\mathbf{q}$.  

\section{Results}
\label{sec:Results}
In this Section, we present the results of several numerical experiments carried out using the novel spectral schemes described in Sections \ref{subsec:Methods_spec_rel} and \ref{subsec:Methods_spec_purg}. We apply the SR and SP schemes to a collection of standard benchmark problems for nonlinear hyperbolic conservation laws, summarised here: $(a)$ 1D inviscid Burgers equation with a single-mode initial condition in a periodic domain in Section~\ref{subsec:Results_burg}, $(b)$ 1D shallow water equations with water-hump and dam-break initial configurations in Section~\ref{subsec:Results_shalwat}, $(c)$ 1D compressible Euler equations with Riemann initial configurations for the Sod-shock tube, the Lax-shock tube, the shock-entropy interaction (Shu--Osher problem), and the blast-waves interaction in Section~\ref{subsec:Results_Euler}.   

We examine the convergence of the SR and SP approximations to the exact solutions when they are available. Otherwise, we perform visual convergence tests using high-resolution approximations from readily available non-spectral schemes as reference solutions. We use the Godunov finite-volume schemes from \texttt{CLAWPACK} (Conservation law package~\cite{clawpack}) for the 1D shallow water problems as well as the blast-waves interaction problem for the 1D compressible Euler equations. For the remaining problems in Section~\ref{subsec:Results_Euler}, we use finite-difference schemes from \texttt{HyPar} (Hyperbolic-Parabolic Partial Differential Equations Solver~\cite{hypar}) In Section~\ref{subsec:Results_burg_ics}, we discuss the variation in the growth of aliasing instabilities depending on the choice of initial data for the 1D-periodic inviscid Burgers equation. The impact of this variation on the SR, SP, and SVV schemes is then discussed in detail. Finally, we apply the SR scheme to the solution of the 1D wall approximation of the 3D-axisymmetric wall-bounded incompressible Euler equation in Section~\ref{subsec:Results_lhmodel}. We compare the estimate of the width of the analyticity-strip $\delta(t)$ obtained using SR approximations at a lower resolution and PPS approximations at a higher resolution. We conclude that SR schemes allow for reliable extraction of $\delta(t)$ for much longer times than the $2/3$-dealiased PPS schemes. Therefore, SR schemes can remedy the constraints on spectral schemes employed for tracking finite-time singularities.  

We use the $4$th-order Runge--Kutta scheme for time integration of the finite-dimensional ODEs resulting from spatial discretisation using PPS schemes of Section~\ref{subsec:Methods_pure_spec}, SR schemes of Section~\ref{subsec:Methods_spec_rel}, and SP schemes of Section~\ref{subsec:Methods_spec_purg}. The choice of time step  $\Delta t$ is such that further decrease of $\Delta t$ does not increase the accuracy of the approximation; therefore, the limit in accuracy is attributed solely to the spatial discretisation scheme.

\subsection{1D inviscid Burgers equation}
\label{subsec:Results_burg}
We begin our studies with the 1D inviscid Burgers equation (also in Eq.~\eqref{eq:burg_pps}) with the initial data given below:
\begin{subequations}    
\begin{align}
    \partial_t u + \partial_x \left( \tfrac{1}{2}u^2 \right) &= 0, \qquad t\in \mathbb{R}^+, \qquad x \in [0,1],  \\
    u(x,t=0) &=\sin \left( 2\pi x \right).    \label{eq:burg_ini}
\end{align}
\label{eq:burg}
\end{subequations}
For the single-mode initial condition in Eq.~\eqref{eq:burg_ini}, this model develops a finite-time singularity at time $t_* = 1/2\pi$, which manifests as a shock wave. Prior to shock formation, we have a unique classical solution for the velocity field which remains smooth and conserves energy, i.e., $||u(\cdot,t)||^2_{L^2}$ does not change with time. Once shock waves are formed, uniqueness and regularity properties are lost and we have to consider solutions in a weak distributional sense. An entropy argument can then be used to pick a physically meaningful weak solution and to restore uniqueness. This \textit{entropic} weak solution is known to be dissipative i.e., $||u(\cdot,t)||^2_{L^2}$ decays in time~\cite{whitham2011linear}. \par

\begin{figure}[htbp]
    \centering
        \begin{tikzpicture}
 	\draw (0,0) node[inner sep=0]{\includegraphics[width=0.32\linewidth]{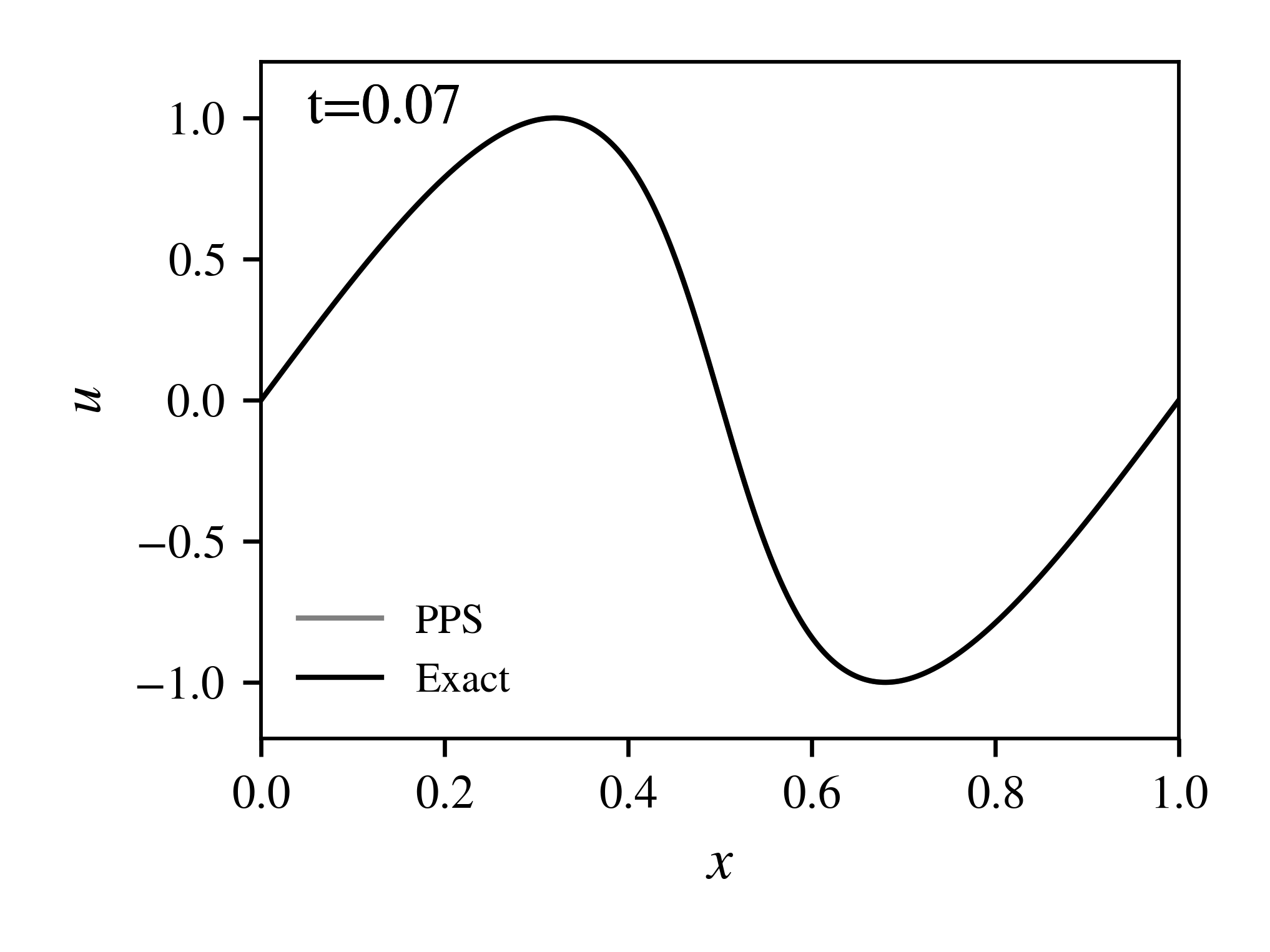}};
		\draw (1.77,1.3) node {$(a)$};
	\end{tikzpicture}
         \begin{tikzpicture}
 	\draw (0,0) node[inner sep=0]{\includegraphics[width=0.32\linewidth]{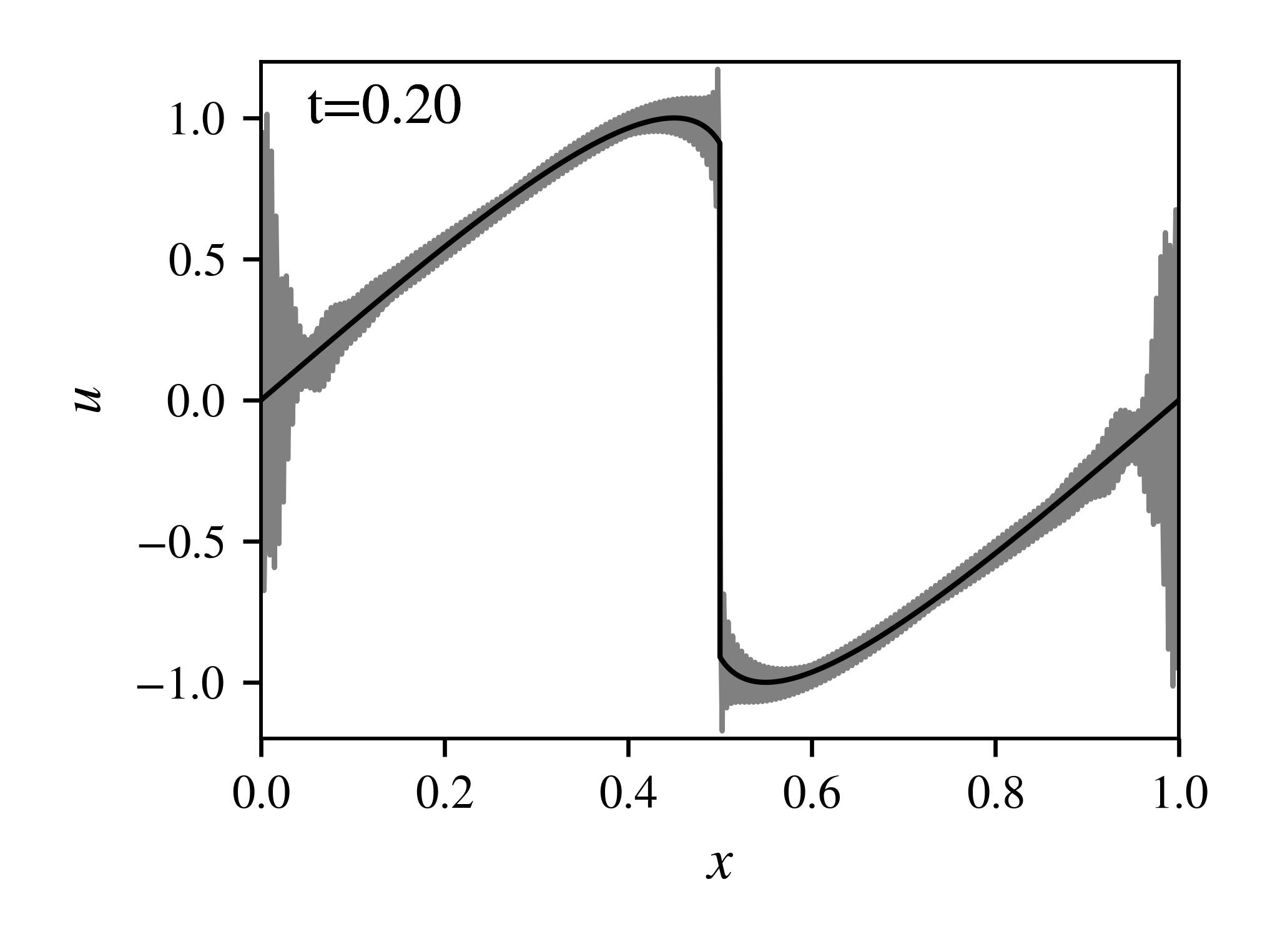}};
		\draw (1.77,1.3) node {$(b)$};
	\end{tikzpicture}
         \begin{tikzpicture}
 	\draw (0,0) node[inner sep=0]{\includegraphics[width=0.32\linewidth]{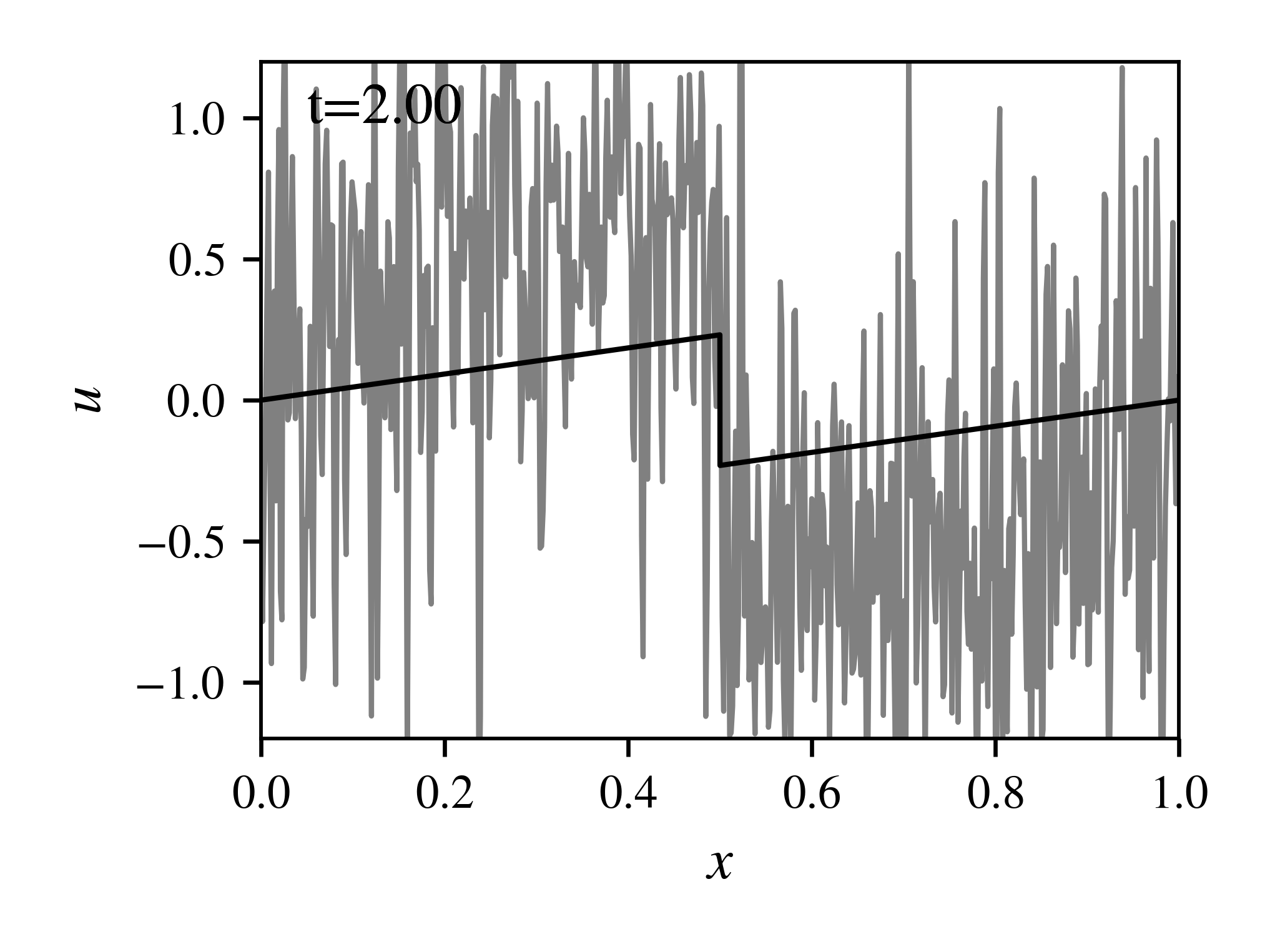}};
		\draw (1.77,1.3) node {$(c)$};
	\end{tikzpicture}
    \caption{Plots vs. $x$ of the velocity field $u(x,t)$ at increasing times from panels $(a)-(c)$. The exact solution (black line) and the $2/3$-dealiased PPS approximation on $N_x=615$ points (grey line) are superposed for comparison. We see the development of Gibbs oscillations near the shock and tygers in the smooth parts of the field which experience a positive strain. In panel $(c)$, we see that the tygers have spread throughout the domain.}
    \label{fig:1burg_tyg}
\end{figure}
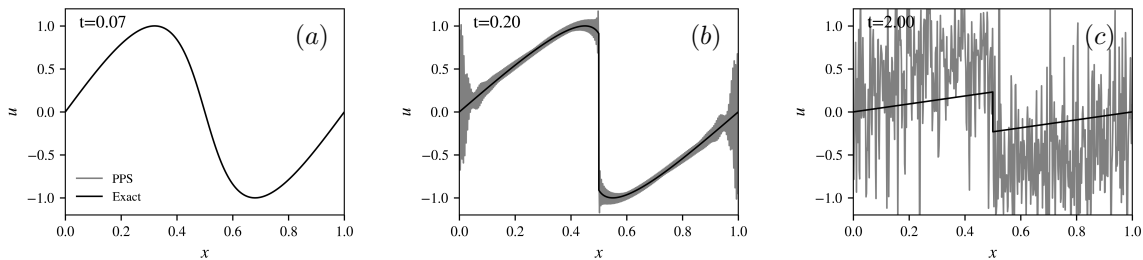

Since spectral approximations for inviscid PDEs are conservative, they cannot model the entropic solution which dissipates energy. In Fig.~\ref{fig:1burg_tyg}, we plot the $2/3$-dealiased PPS approximation of Eqs.~\eqref{eq:burg} for a spatial resolution of $N_x=2N+1=615$ at different times. We also plot the exact solution (black line) for comparison.
 In panel $(a)$, the PPS approximation converges to the exact entropic solution which is still smooth at early times. The finite-time singularity occurs at time $t_* = 1/2\pi \simeq 0.159$ and develops as a shock at $x_* = 0.5$. In panel $(b)$ of Fig.~\ref{fig:1burg_tyg}, we see the development of Gibbs oscillations around the shock as the system has just crossed the singularity time. We also observe tygers~\cite{ray2011resonance} as a rapidly-growing envelope of oscillations at the edges of the domain (smooth part of the solution).  
Note that the oscillations occur at the Galerkin truncation wavelength $\lambda_G = 2\pi/K_G$, where $K_G$ is the dealiasing cut-off. At long times, shown in panel $(c)$, the system is thermalising and the PPS approximation has almost settled down to a Gaussian noise~\cite{di2018dynamics}. 

\begin{figure}[h!]
 \centering
 \begin{tikzpicture}
\draw (0,0) node[inner sep=0]{\includegraphics[width=0.45\linewidth]{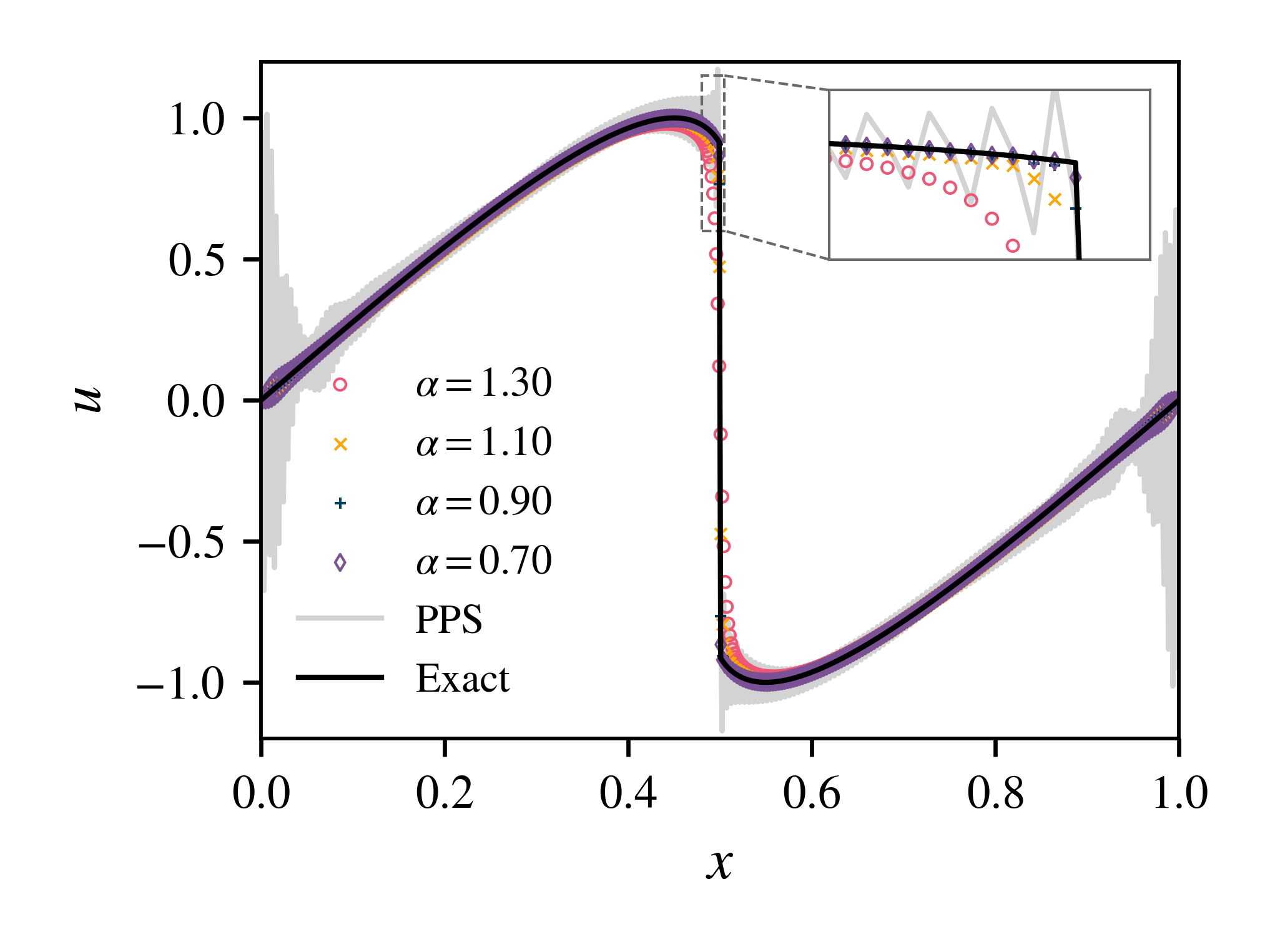}};
    \draw (2.1,-1.2) node {$(a)$};
\end{tikzpicture}
\begin{tikzpicture}
\draw (0,0) node[inner sep=0]{\includegraphics[width=0.45\linewidth]{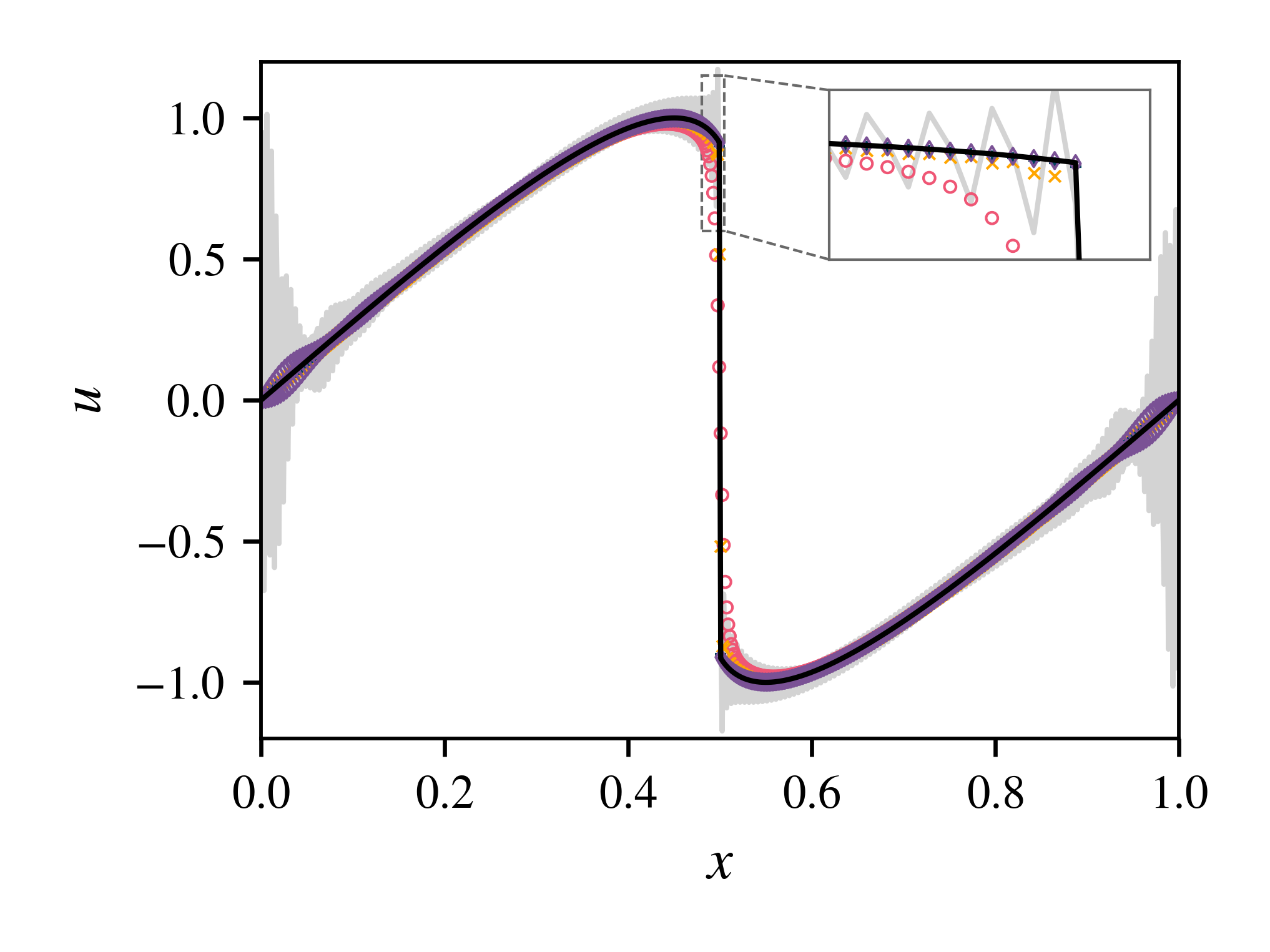}};
    \draw (2.1,-1.2) node {$(b)$};
\end{tikzpicture}
 \caption{Plots vs. $x$ of the velocity field $u(x,t)$ obtained by using $(a)$ spectral relaxation (SR) and $(b)$ spectral purging (SP) schemes with the Fej\'er--Korovkin kernel. The value of $\alpha$ used is given in the legend; we set $\gamma=0.99$. The plots are shown for time $t=0.2$ after the shock formation at $t_* \simeq 0.159$. We compare these approximations with the exact solution (black line) and the $2/3$-dealiased PPS solution (light grey line). The insets in each panel are used to zoom into the region of the shock.}
 \label{fig:2SPSR-diffalpha-feko}
\end{figure}

We now discuss the SR and SP approximations to solve Eqs.~\eqref{eq:burg}. A crucial step in implementing the new schemes is the choice of parameters $(\alpha,\gamma)$ defined in Eq.~\eqref{eq:sr_params}. In panels $(a)$ and $(b)$ of Fig.~\ref{fig:2SPSR-diffalpha-feko}, we show SR and SP approximations of the velocity field using the Fej\'er--Korovkin kernel, respectively, at $t=0.2$. The different keys in the legend represent different values of $\alpha$; we set $\gamma=0.99$. Approximations obtained with higher values of $\alpha$ suffer more dissipation and are overdamped near the shock. On the other hand, approximations obtained using lower values of $\alpha$ suffer from persistent tyger-like structures since the damping is insufficient. We note that the novel spectral approximations do not have oscillatory behaviour near the shock, as seen in the zoomed insets in panels $(a)-(b)$ of Fig.~\ref{fig:2SPSR-diffalpha-feko}. 

Hereafter, we use the abbreviation SR-FeKo to denote the spectral relaxation scheme using the Fej\'er--Korovkin kernel, and SP-FeKo for the spectral purging schemes using the same kernel. 

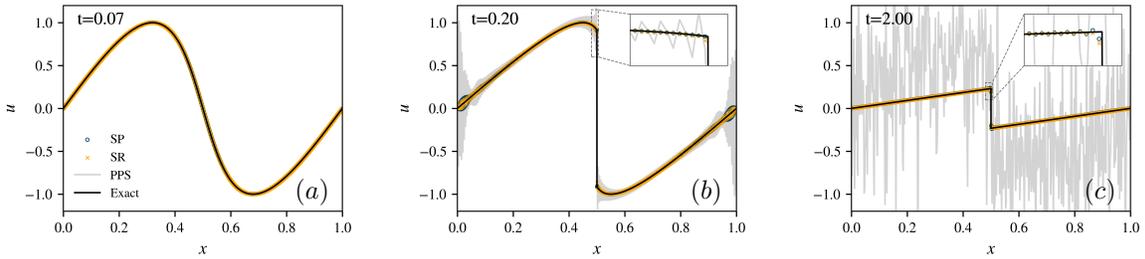
\begin{figure}[h!]
     \centering
        \begin{tikzpicture}
 	\draw (0,0) node[inner sep=0]{\includegraphics[width=0.32\linewidth]{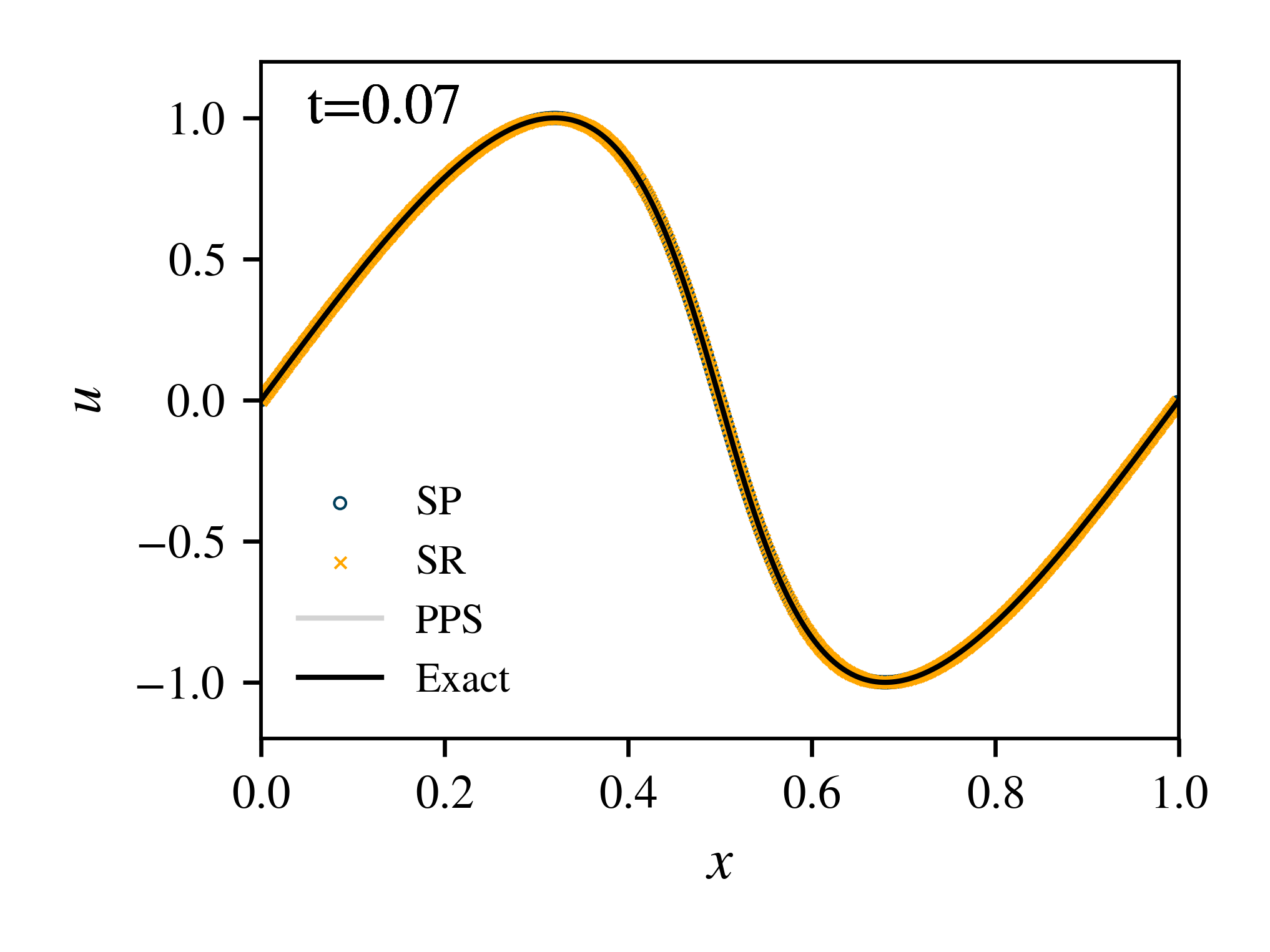}};
		\draw (1.8,-0.8) node {$(a)$};
	\end{tikzpicture}
         \begin{tikzpicture}
 	\draw (0,0) node[inner sep=0]{\includegraphics[width=0.32\linewidth]{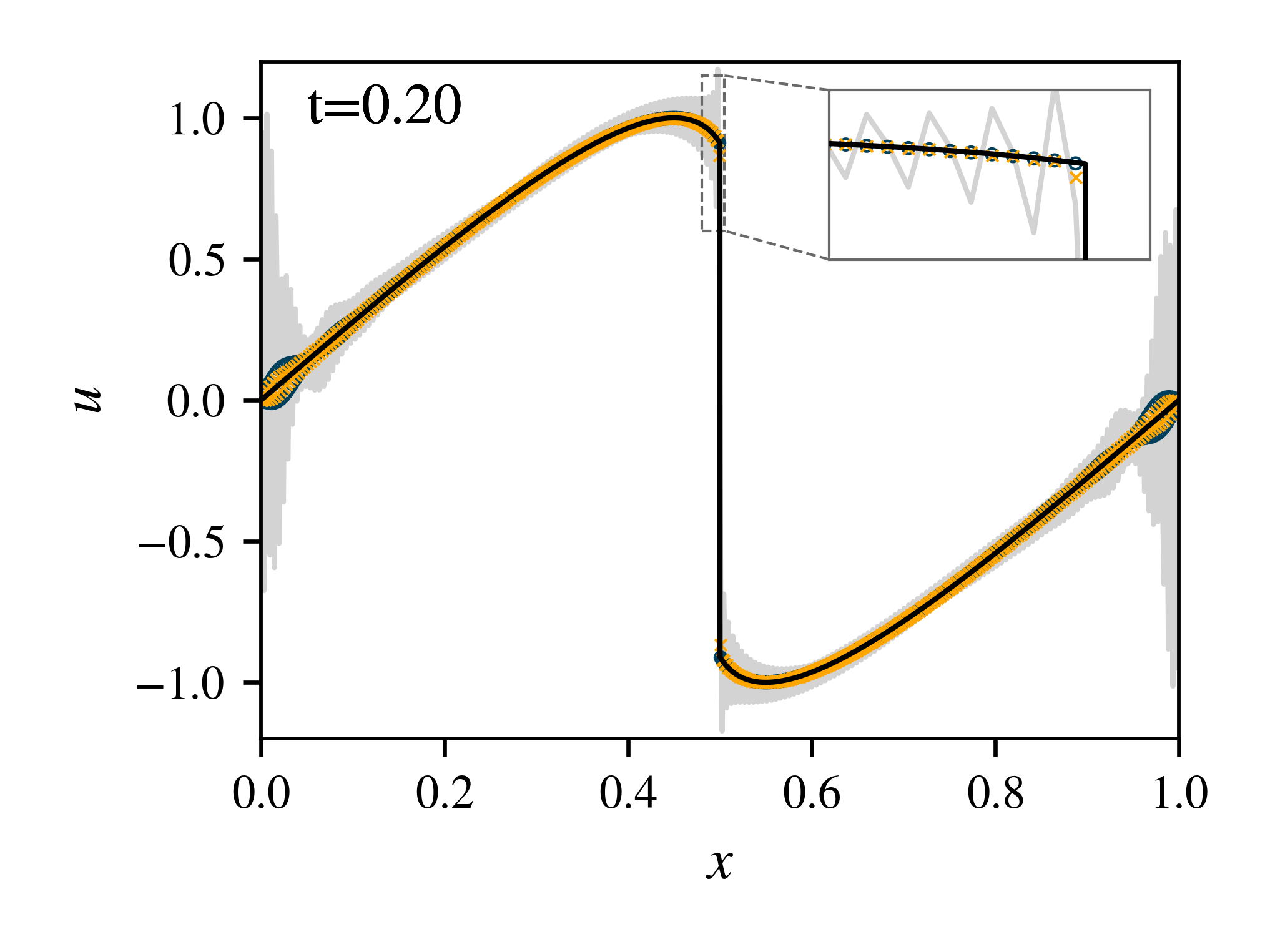}};
		\draw (1.8,-0.8) node {$(b)$};
	\end{tikzpicture}
         \begin{tikzpicture}
 	\draw (0,0) node[inner sep=0]{\includegraphics[width=0.32\linewidth]{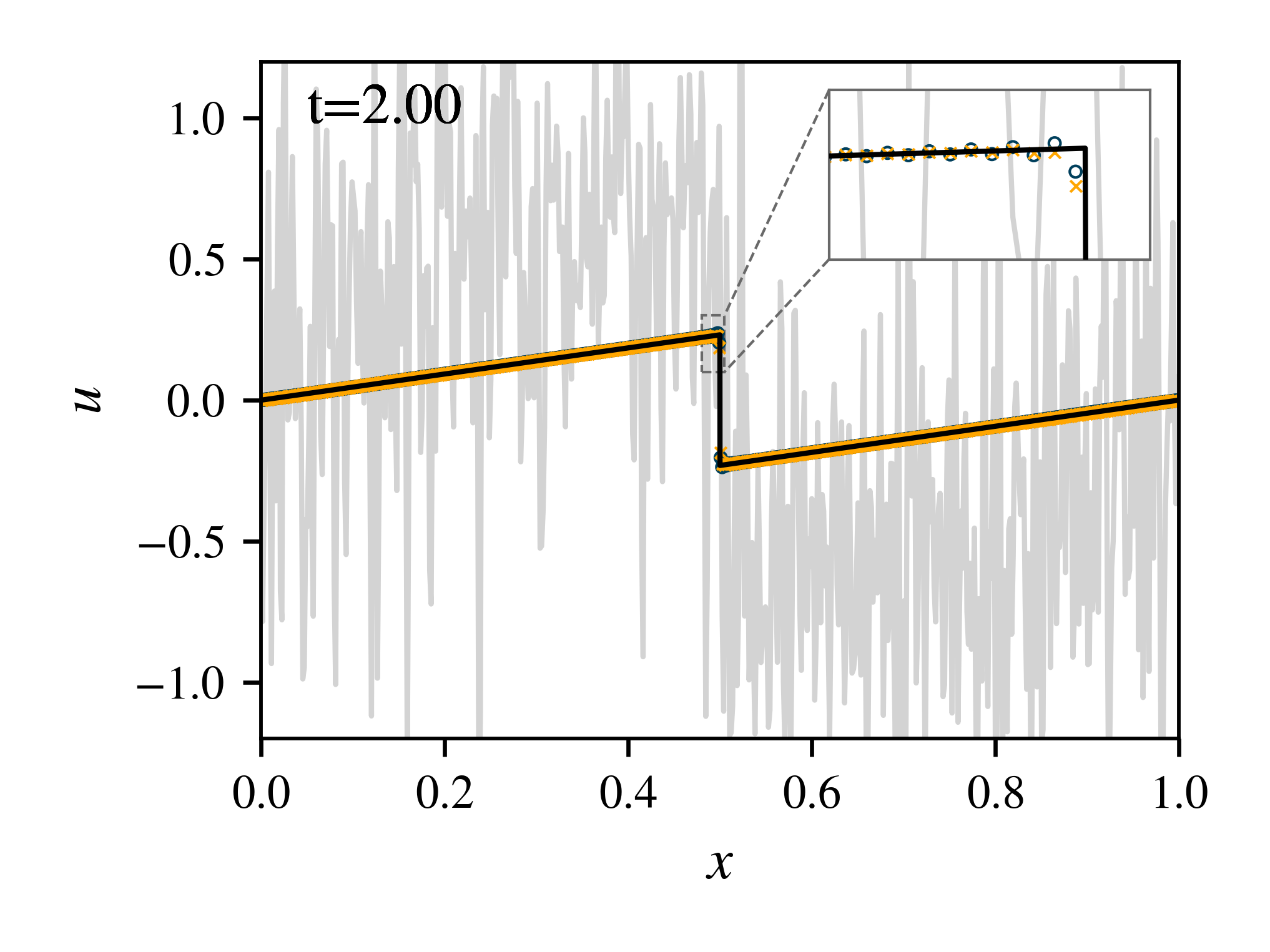}};
		\draw (1.8,-0.8) node {$(c)$};
	\end{tikzpicture}
     \caption{Plots vs. $x$ of the velocity field $u$ approximated using the best SR-FeKo (orange crosses) and SP-FeKo (blue circles) schemes. The parameter sets are, for SR-FeKo: $(\alpha=0.7, \gamma=0.99)$ and for SP-FeKo: $(\alpha=0.65, \gamma=0.99)$. The exact solution (black line) and the $2/3$-dealiased PPS solution (grey line) are plotted for comparison. We use a spatial resolution of $N_x = 615$ for all the runs shown here.}
     \label{fig:3SPSR_best_time}
 \end{figure}

From visual inspection of Fig.~\ref{fig:2SPSR-diffalpha-feko}, we can guess optimal parameter values for this problem as $(\alpha \simeq 0.7,\gamma=0.99)$ for both SR-FeKo and SP-FeKo schemes. To arrive at the best parameter set, we perform a more exhaustive parameter-search in $(\alpha,\gamma)$ space; we seek to minimise $(a)$ the time-integrated $L^2$ error from $t=0$ to $t=9.9$ of the SR and SP approximations and $(b)$ the maximum $L^2$ error over the same time interval. The best parameter sets obtained via this minimisation routine are given by: $(\alpha=0.7,\gamma=0.99)$ for the SR-FeKo scheme and $(\alpha=0.65,\gamma=0.99)$ for the SP-FeKo scheme, for a spatial resolution of $N_x=615$ points. 
In Fig.~\ref{fig:3SPSR_best_time}, we plot the SR-FeKo and SP-FeKo approximations obtained using the best parameter sets above, for three representative times. At early times $t=0.07$ shown in panel $(a)$, we have spectral convergence: PPS, SR-FeKo and SP-FeKo approximations represent the exact solution accurately. Note that the SR and SP approximations are not over-dissipated. In panel $(b)$, with the shock formation at $t=0.2$, we see the development of tygers in the PPS approximation. SR and SP approximations closely follow the shock (black line) and, furthermore, respect monotonicity (seen in the zoomed inset axis). We notice a difference far from the shock: the SR-FeKo scheme has successfully damped the tyger structures, thus giving an approximation that is very close to the exact solution. This is not the case for the SP-FeKo approximation, which still displays attenuated tyger-like structures. Since SP approximations are discontinuous-in-time, they can only approximate the exact solution well at discrete times. For long-time behaviour, in panel $(c)$, the SR and SP approximations are indistinguishable from the true solution. At this time, the PPS approximation is noisy and is moving towards thermalisation. 

\begin{figure}[h!]
\centering
    \includegraphics[width=0.5\linewidth]{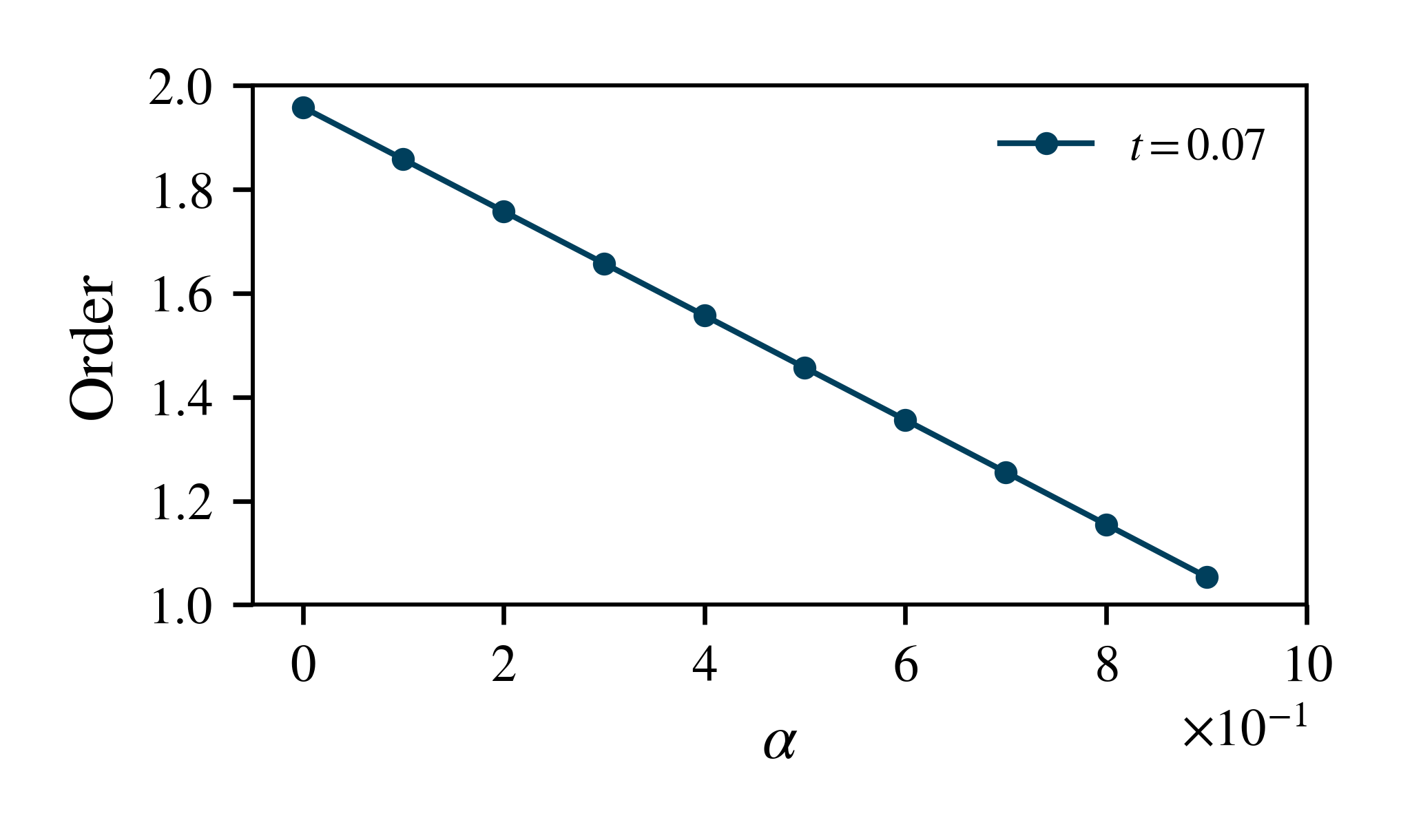}
    \caption{Plot vs. $\alpha$ of the order of convergence calculated from $L^\infty$-errors of the SR-FeKo approximation with $(\alpha,\gamma=0.99)$ at early time $t=0.07$ when the solution is still smooth.}
    \label{fig:t1_conv_sr_alpha}
\end{figure}

Reference~\cite{bessemain} shows that the SR and SP approximations, for a suitable choice of the parameters $(\alpha,\gamma)$, converge in the $L^2$ norm to the exact entropic solution, as $N \rightarrow \infty$. We perform numerical convergence analysis of the SR-FeKo approximation for the best parameter set $(\alpha=0.7,\gamma=0.99)$ and present the results in Table~\ref{tab:1_cs_ic0-conv-feko}. We evaluate the error w.r.t the exact entropic solution for 8 well-separated resolutions, at 3 different times; the order of the method is then calculated from the error of adjacent resolutions. As concluded by~\cite{bessemain}, SR-FeKo approximations indeed converge, in both norms, to the entropic solution. The asymptotic order of convergence is the same for $L^1$ and $L^2$ errors at early times (see column 2 of Table~\ref{tab:1_cs_ic0-conv-feko}), but differs at later times (see columns 3 and 4 of Table~\ref{tab:1_cs_ic0-conv-feko}). We attribute this difference to the increased sensitivity of the $L^2$ error to deviations.  
At early times $t=0.07$, we have polynomial convergence of the method with an asymptotic order of $\simeq 1.28$ using the $L^1$ error, and $\simeq 1.27$ for the $L^2$ error. For the times following the shock formation, we retain an order $\simeq 0.9$ for the  $L^1$ error, and $\simeq 0.78$ for the $L^2$ error (column 4, Table~\ref{tab:1_cs_ic0-conv-feko}). 
The order of convergence of the SR scheme is given by the scaling of the RHS of Eq.~\eqref{eq:burg_sr} with $N$ as $\tfrac{1}{\tau}(\tfrac{1}{m})^2 \simeq N^{-(2\gamma-\alpha)} $; for our choice of $\gamma \simeq 1$, we have $N^{-(2-\alpha)}$. The scaling is confirmed in Fig.~\ref{fig:t1_conv_sr_alpha}, where we plot the asymptotic order vs. the value of $\alpha$ used for SR-FeKo approximations. We set $\gamma=0.99$ for all the simulations performed here. 
Before the shock formation, since the solution is analytic in space, the order of convergence of the scheme is not limited by the spatial regularity of the exact solution but only by how we approximate it. But after the shock formation, since the solution $u(\cdot, t) \in {\rm BV}_x$ (the space of functions of bounded variation in space), the order of convergence is now limited by the spatial regularity of the solution $u$; this order is now bounded above by 1.





 \begin{table}[h!]
 \centering
  \begin{tabularx}{\linewidth}{@{}|YY|YY|YY|YY|@{}}
     \hline 
     \rule{0pt}{2.5ex}    
      SR & &  \multicolumn{2}{>{\hsize=\dimexpr1\hsize+4\tabcolsep+1\arrayrulewidth\relax}X|}{\centering $t=0.07$} & \multicolumn{2}{>{\hsize=\dimexpr2\hsize+4\tabcolsep+1\arrayrulewidth\relax}X|}{\centering $t=0.20$} &
      \multicolumn{2}{>{\hsize=\dimexpr2\hsize+4\tabcolsep+1\arrayrulewidth\relax}X|}{\centering $t=2.00$} \\
       FeKo & $N_x$ & Error & Order & Error & Order & Error & Order \\ 
 \hline
             \rule{0pt}{3ex}
\multirow{8}{*}{ $L^1$ error}&  39    & $4.5 \cdot 10^{-3}$ &      & $3.8 \cdot 10^{-2}$ &      & $9.0 \cdot 10^{-3}$ &      \\
&  65    & $2.6 \cdot 10^{-3}$ & 1.11 & $2.6 \cdot 10^{-2}$ & 0.78 & $6.3 \cdot 10^{-3}$ & 0.71 \\
&  123   & $1.2 \cdot 10^{-3}$ & 1.26 & $1.6 \cdot 10^{-2}$ & 0.72 & $3.3 \cdot 10^{-3}$ & 1.03 \\
&  205   & $6.1 \cdot 10^{-4}$ & 1.25 & $1.1 \cdot 10^{-2}$ & 0.74 & $1.9 \cdot 10^{-3}$ & 1.06 \\
&  615   & $1.5 \cdot 10^{-4}$ & 1.26 & $4.6 \cdot 10^{-3}$ & 0.81 & $6.5 \cdot 10^{-4}$ & 0.97 \\
&  1599  & $4.5 \cdot 10^{-5}$ & 1.28 & $2.0 \cdot 10^{-3}$ & 0.87 & $2.8 \cdot 10^{-4}$ & 0.90 \\
&  2665  & $2.4 \cdot 10^{-5}$ & 1.28 & $1.3 \cdot 10^{-3}$ & 0.90 & $1.8 \cdot 10^{-4}$ & 0.90 \\
&  7995  & $5.8 \cdot 10^{-6}$ & 1.28 & $4.6 \cdot 10^{-4}$ & 0.91 & $6.5 \cdot 10^{-5}$ & 0.90 \\
 \hline
             \rule{0pt}{3ex}
   \multirow{8}{*}{$L^2$ error}& 39   & $5.4 \cdot 10^{-3}$ &      &  $4.9 \cdot 10^{-2}$ &       & $1.9 \cdot 10^{-2}$&        \\
   & 65   & $3.1 \cdot 10^{-3}$ & 1.10 &  $3.3 \cdot 10^{-2}$ &  0.74 & $1.4 \cdot 10^{-2}$ & 0.59 \\
   & 123  & $1.4 \cdot 10^{-3}$ & 1.25 &  $2.2 \cdot 10^{-2}$ &  0.66 & $9.0 \cdot 10^{-3}$ & 0.73 \\
   & 205  & $7.4 \cdot 10^{-4}$ & 1.24 &  $1.6 \cdot 10^{-2}$ &  0.61 & $6.1 \cdot 10^{-3}$ & 0.75 \\
   & 615  & $1.8 \cdot 10^{-4}$ & 1.25 &  $7.8 \cdot 10^{-3}$ &  0.67 & $2.6 \cdot 10^{-3}$ & 0.76 \\
   & 1599 & $5.5 \cdot 10^{-5}$ & 1.27 &  $3.7 \cdot 10^{-3}$ &  0.76 & $1.3 \cdot 10^{-3}$ & 0.78 \\
   & 2665 & $2.8 \cdot 10^{-5}$ & 1.27 &  $2.5 \cdot 10^{-3}$ &  0.79 & $8.4 \cdot 10^{-4}$ & 0.78 \\
   & 7995 & $7.0 \cdot 10^{-6}$ & 1.27 &  $1.0 \cdot 10^{-3}$ &  0.81 & $3.6 \cdot 10^{-4}$ & 0.78 \\
 \hline
     \end{tabularx}
     \caption{Convergence analysis for the SR-FeKo approximations using $(\alpha=0.7,\gamma=0.99)$ at three representative times before and after the shock formation $t=0.07,0.2$ and $2.0$. Here, the second row gives convergence analysis based on the $L^1$ error; and the third row shows the $L^2$ error. All sets of runs do not employ $2/3$-dealiasing.} 
     \label{tab:1_cs_ic0-conv-feko}
 \end{table}%

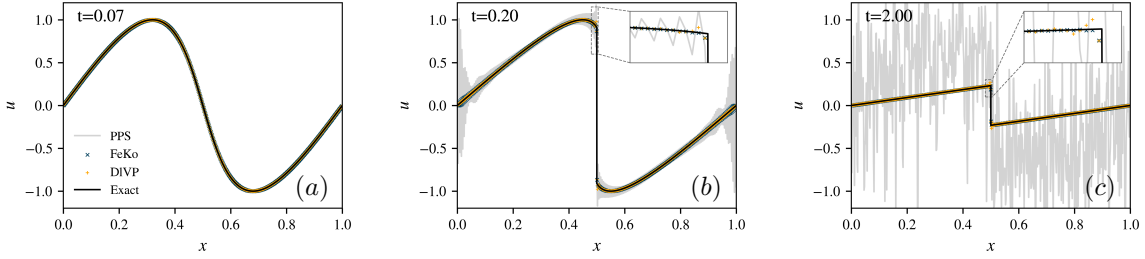
\begin{figure}[h!]
     \centering
        \begin{tikzpicture}
 	\draw (0,0) node[inner sep=0]{\includegraphics[width=0.32\linewidth]{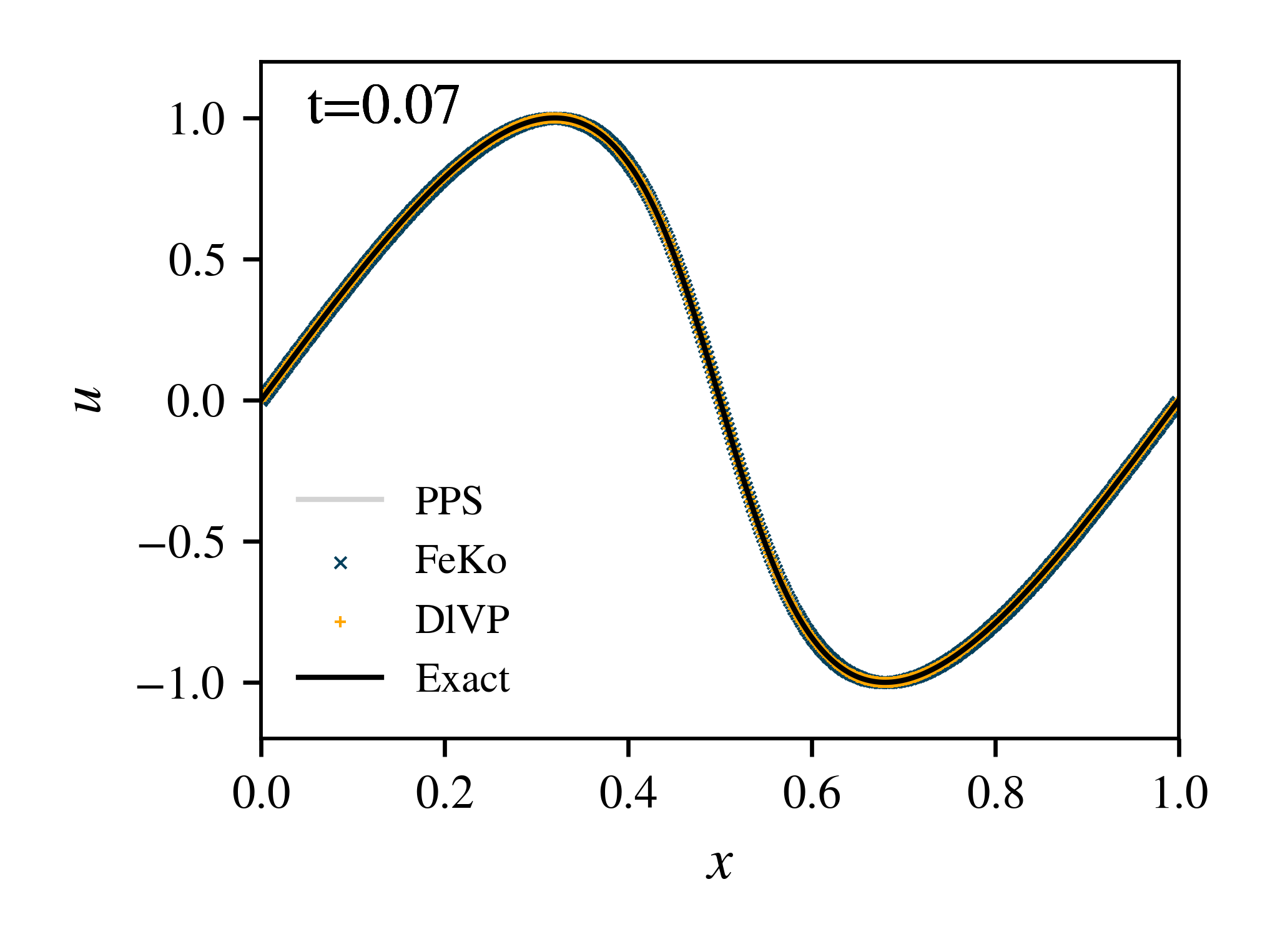}};
		\draw (1.8,-0.8) node {$(a)$};
	\end{tikzpicture}
         \begin{tikzpicture}
 	\draw (0,0) node[inner sep=0]{\includegraphics[width=0.32\linewidth]{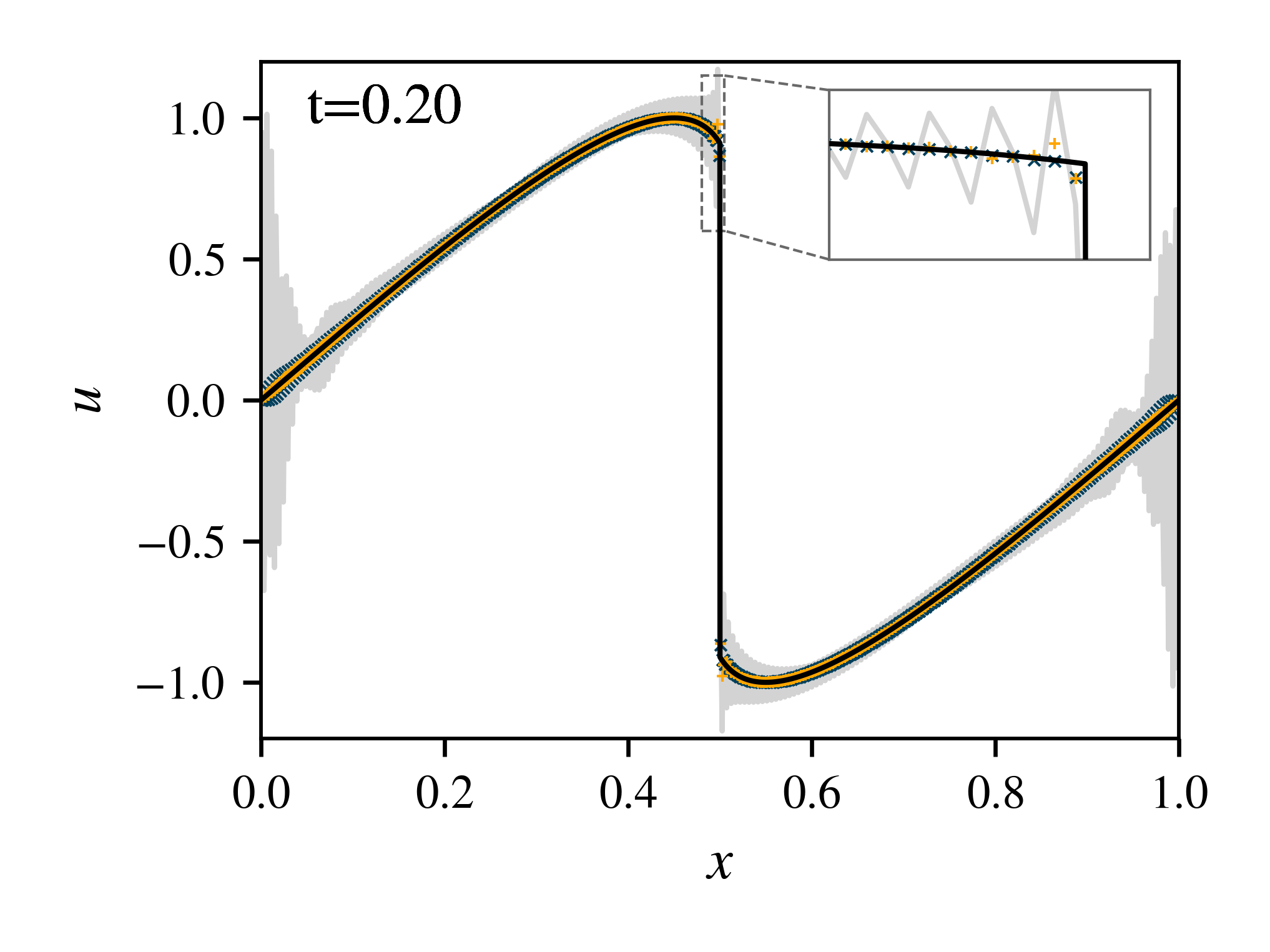}};
		\draw (1.8,-0.8) node {$(b)$};
	\end{tikzpicture}
         \begin{tikzpicture}
 	\draw (0,0) node[inner sep=0]{\includegraphics[width=0.32\linewidth]{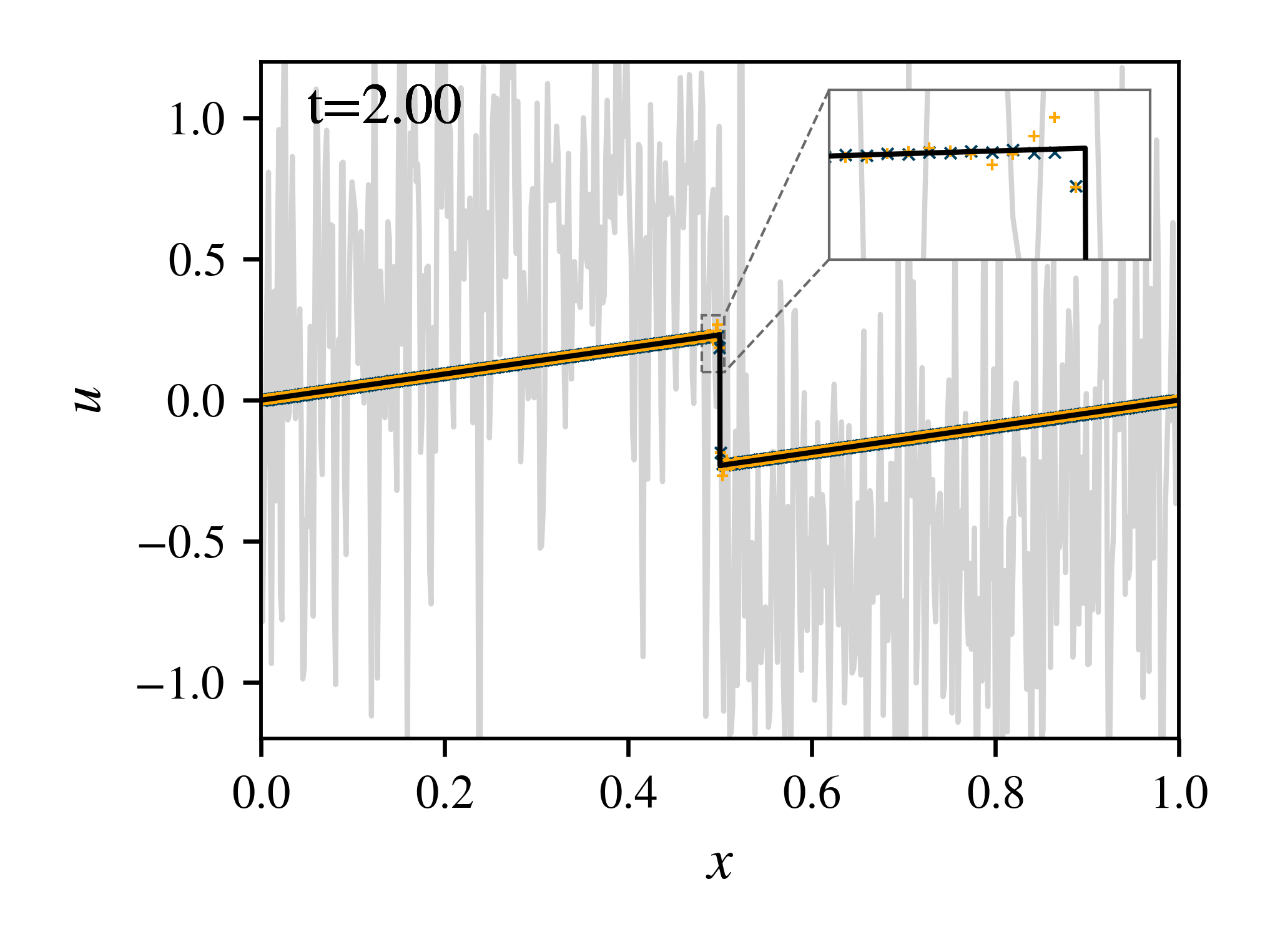}};
		\draw (1.8,-0.8) node {$(c)$};
	\end{tikzpicture}
 \caption{Plots vs. $x$ of the velocity field $u$ obtained by using SR-FeKo for $(\alpha=0.7,\gamma=0.99)$ (blue crosses) and SR-DlVP (with the de La Vall\'ee Poussin filter) for $(\alpha=0.89,\gamma=0.9)$ (orange plusses) at different representative times in panels $(a)$-$(c)$ before and after the shock formation. We compare these approximations with the exact solution (black line) and the $2/3$-dealiased PPS solution (light grey line). The zoomed insets in each panel show Gibbs-like oscillations in SR-DlVP approximation. The spatial resolution is $N_x=615$.}
 \label{fig:DlVP-time}
\end{figure}

We now use the SR scheme with the de La Vall\'ee Poussin kernel (SR-DlVP) which is a non-positive high-order approximation kernel.  The functional form and parameters for this filter are given in Appendix~\ref{app:filters}. We deduce the optimal parameters $(\alpha=0.89, \gamma=0.9)$ for SR-DlVP by the same minimisation scheme, as described earlier for SR-FeKo.  
In Fig.~\ref{fig:DlVP-time}, we plot the best SR-FeKo approximation with the best SR-DlVP approximation at different times; we set $N_x=615$. We also perform convergence analysis for the best SR-DlVP approximation using the $L^1$ and $L^2$  norms, as shown in Table.~\ref{tab:2_cs_ic0-conv-dlvp} (rows 2 and 3 respectively). For times before the shock, the SR-DlVP shows spectral convergence (see column 2 of Table~\ref{tab:2_cs_ic0-conv-dlvp}); the difference is not appreciable in panel $(a)$ of Fig.~\ref{fig:DlVP-time}. However, just after the shock formation, in panel $(b)$ of Fig.~\ref{fig:DlVP-time}, SR-DlVP approximation shows oscillations near the shock with a small, bounded amplitude which are absent for the SR-FeKo approximation. Similar oscillations were observed for the SVV method by~\cite{karamanos2000spectral}. After the shock formation, SR-FeKo and SR-DlVP schemes have the same rate of convergence in terms of the  $L^1$ norm, for approximating shock solutions which have only BV smoothness in space. If we use the $L^2$ norm, the order of convergence for SR-FeKo ($0.8$) is higher than that for the SR-DlVP ($0.5$) (see columns 3-4 of Tables~\ref{tab:1_cs_ic0-conv-feko}-\ref{tab:2_cs_ic0-conv-dlvp}). This results may seem surprising knowing that the de La Vall\'ee Poussin kernel is a high-order approximation kernel, whereas the Fej\'er--Korovkin kernel is a low-order one. In fact, as it is emphasised in~\cite{leveque2002finite}, ``The order of accuracy is not everything"; the faster global convergence of SR-FeKo, in $L^2$ norm, can be traced back to the Fej\'er--Korovkin kernel which retains the monotonicity and maximum-principle properties of the solution in the approximate one. Since this is not true for the de La Vall\'ee Poussin kernel, we see bounded oscillations and under-shoots or over-shoots of the approximation near discontinuities. Thus, high-order approximations as the sole criterion for designing accurate numerical schemes can be misleading. Even so, if we compute the order of convergence only in the smooth regions of the solution and exclude the discontinuity, then we expect the SR-DlVP schemes to perform better. We conclude that SR approximations, which use the high-order de La Vall\'ee Poussin kernel, yield better convergence when the solution is smooth, but only limited convergence for the shock capture compared to the low-order Fej\'er--Korovkin kernel.

However, the choice of kernel is special and depends on the properties of the problem at hand. For problems where we do not expect monotonicity, such as the Shu--Osher problem in Section~\ref{subsec:Results_Euler}, we find that the SR-DlVP approximation is desirable and allows us to capture shocks and fine structures in the solutions.  

\begin{table}[h!]
\centering
 \begin{tabularx}{\linewidth}{@{}|YY|YY|YY|YY|@{}}
    \hline 
    \rule{0pt}{2.5ex}    
     SR & &  \multicolumn{2}{>{\hsize=\dimexpr1\hsize+4\tabcolsep+1\arrayrulewidth\relax}X|}{\centering $t=0.07$} & \multicolumn{2}{>{\hsize=\dimexpr2\hsize+4\tabcolsep+1\arrayrulewidth\relax}X|}{\centering $t=0.20$} &
      \multicolumn{2}{>{\hsize=\dimexpr2\hsize+4\tabcolsep+1\arrayrulewidth\relax}X|}{\centering $t=2.00$} \\
      DlVP & $N_x$ & Error & Order & Error & Order & Error & Order \\ 
        \hline        
            \rule{0pt}{3ex}
\multirow{8}{*}{ $L^1$ error}  & $39$    & $1.2 \cdot 10^{-5} $ &         & $3.7 \cdot 10^{-2} $ &        &  $1.0 \cdot 10^{-2} $ &         \\ 
& $65  $  & $3.8 \cdot 10^{-7} $ &  $6.82 $  & $2.4 \cdot 10^{-2} $ &  $0.85$  &  $5.7 \cdot 10^{-3} $  &   $1.13$  \\
& $123 $  & $6.9 \cdot 10^{-10}$ &  $9.89 $  & $1.4 \cdot 10^{-2} $ &  $0.83$  &  $3.2 \cdot 10^{-3} $  &   $0.90$  \\
& $205 $  & $2.8 \cdot 10^{-13}$ &  $15.28$  & $9.0 \cdot 10^{-3} $ &  $0.88$  &  $1.9 \cdot 10^{-3} $  &   $1.03$  \\
& $615 $  & $5.0 \cdot 10^{-15}$ &  $3.66 $  & $3.2 \cdot 10^{-3} $ &  $0.93$  &  $6.3 \cdot 10^{-4} $  &   $0.99$  \\
& $1599$  & $4.1 \cdot 10^{-15}$ &  $0.21 $  & $1.3 \cdot 10^{-3} $ &  $0.96$  &  $2.6 \cdot 10^{-4} $  &   $0.92$  \\
& $2665$  & $4.0 \cdot 10^{-15}$ &  $0.05 $  & $7.9 \cdot 10^{-4} $ &  $0.95$  &  $1.6 \cdot 10^{-4} $  &   $0.91$  \\
& $7995$  & $1.0 \cdot 10^{-14}$ &  $-0.83$  & $2.7 \cdot 10^{-4} $ &  $0.95$  &  $6.0 \cdot 10^{-5} $  &   $0.93$  \\
           \hline
             \rule{0pt}{3ex}
\multirow{8}{*}{$L^2$ error}  & $39$     & $1.6 \cdot 10^{-5}$  &        &   $4.1 \cdot 10^{-2}$  &        & $1.5 \cdot 10^{-2}$  &      \\                      
 & $65  $   & $5.7\cdot10^{-7}$  & $6.54 $  &   $2.7\cdot10^{-2}$ &  $0.76 $ & $1.1\cdot 10^{-2} $& $0.65$  \\
 & $123 $   & $1.1\cdot10^{-9}$  & $9.88 $  &   $1.7\cdot10^{-2}$ &  $0.70 $ & $7.6\cdot 10^{-3} $& $0.50$  \\
 & $205 $   & $4.4\cdot10^{-13}$ & $15.23$  &   $1.2\cdot10^{-2}$ &  $0.71 $ & $6.0\cdot 10^{-3} $& $0.46$  \\
 & $615 $   & $7.0\cdot10^{-15}$ & $3.76 $  &   $5.7\cdot10^{-3}$ &  $0.70 $ & $3.6\cdot 10^{-3} $& $0.44$  \\
 & $1599$   & $5.2\cdot10^{-15}$ & $0.31 $  &   $3.1\cdot10^{-3}$ &  $0.60 $ & $2.4\cdot 10^{-3} $& $0.44$  \\  
 & $2665$   & $5.7\cdot10^{-15}$ & $-0.17$  &   $2.4\cdot10^{-3}$ &  $0.54 $ & $1.9\cdot 10^{-3} $& $0.44$  \\
 & $7995$   & $1.1\cdot10^{-14}$ & $-0.59$  &   $1.3\cdot10^{-3}$ &  $0.50 $ & $1.1\cdot 10^{-3} $& $0.44$  \\
           \hline
    \end{tabularx}
    \caption{Convergence analysis for the SR-DlVP approximations using $(\alpha=0.89,\gamma=0.9)$ at three representative times before and after the shock formation $t=0.07,0.2$ and $2.0$. Here, the second row gives convergence analysis based on the $L^1$ error and the third row, for the $L^2$ error. All sets of runs do not employ $2/3$-dealiasing.} 
    \label{tab:2_cs_ic0-conv-dlvp}
\end{table}%

For SR and SP implementations of the 1D inviscid Burgers equation, the choices of the kernel $K_m(x)$ used in the convolution operator $\mathcal{M}_m$ and governing parameters $(\alpha,\gamma)$ are guided by the convergence analysis of~\cite{bessemain}. In panel $(a)$ of Fig.~\ref{fig:SPSR-difffilter}, we show SR approximations obtained by using different positive kernels, for the optimal parameter choice $(\alpha=0.7,\gamma=0.99)$. In panel $(b)$ of the same figure, we show the corresponding SP approximations with the optimal parameter set $(\alpha=0.65,\gamma=0.99)$. There is minimal variation across different choices of the allowed filters; a closer examination near the shock reveals that Fej\'er--Korovkin kernel gives a slightly better result compared to the Jackson and Jackson--de La Vall\'{e} Poussin kernels (see zoomed insets). We present the functional forms of these filters in the subsection on positive kernels in Appendix~\ref{app:filters}.  

In Fig.~\ref{fig:SRSPSVV}, we compare the best SR-FeKo and SP-FeKo approximations with the best SVV approximation. SVV, SR and SP methods approximate the exact solution well far from the shock; however, near the shock in panels $(b,c)$ we see that small Gibbs-like oscillations still persist in the SVV approximation. The SR-FeKo and SP-FeKo methods preserve the monotonicity of the solution and approximate the shock very closely. SP-FeKo approximations still display persistent tygers in the smooth regions of the flow, due to their discontinuous-in-time nature. In order to simplify the analysis, we choose to focus our study on the continuous-in-time SR approximations.  

\begin{figure}[h!]
 \centering
 \begin{tikzpicture}
\draw (0,0) node[inner sep=0]{\includegraphics[width=0.45\linewidth]{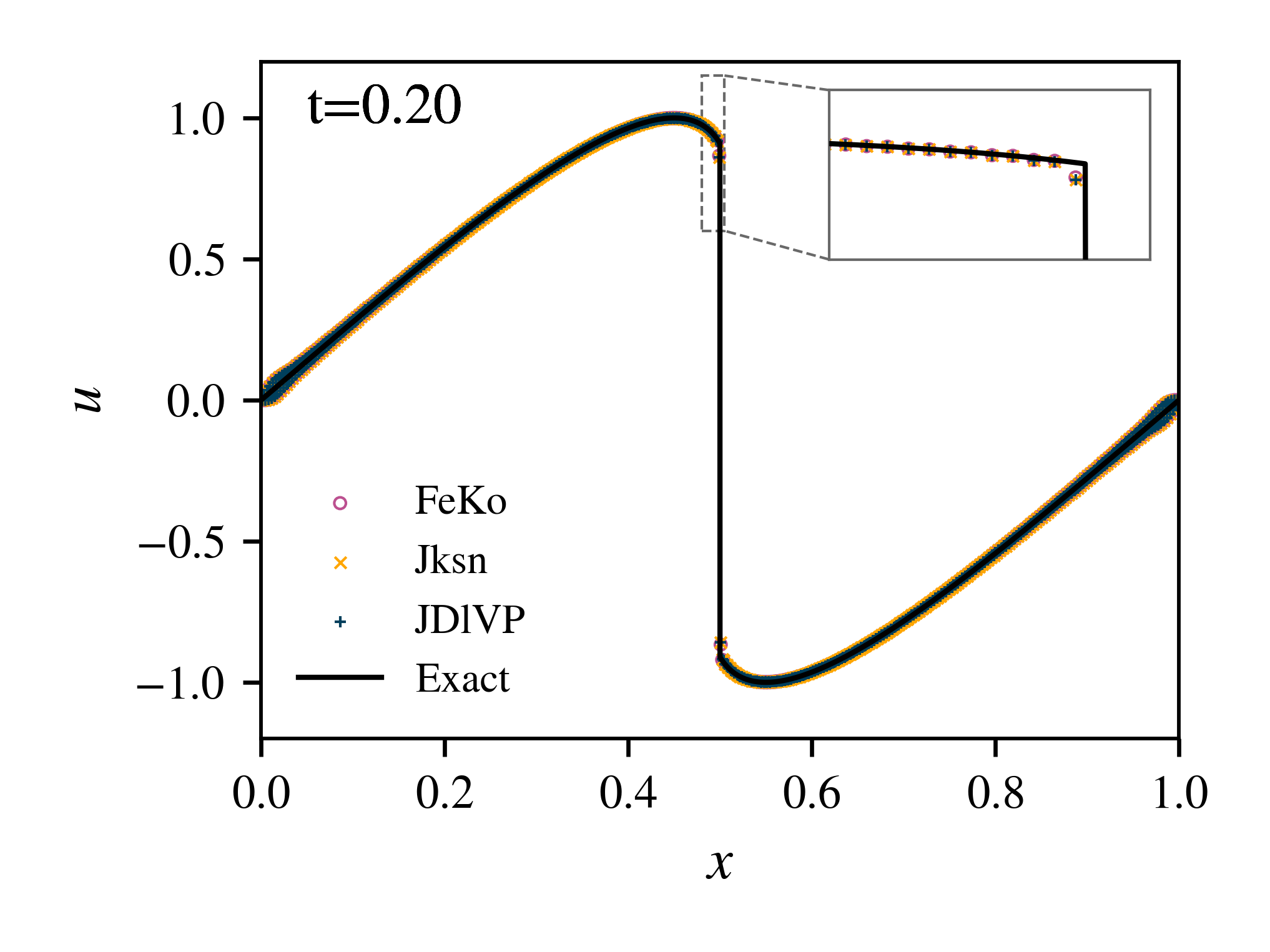}};
    \draw (2.8,-1.2) node {$(a)$};
\end{tikzpicture}
\begin{tikzpicture}
\draw (0,0) node[inner sep=0]{\includegraphics[width=0.45\linewidth]{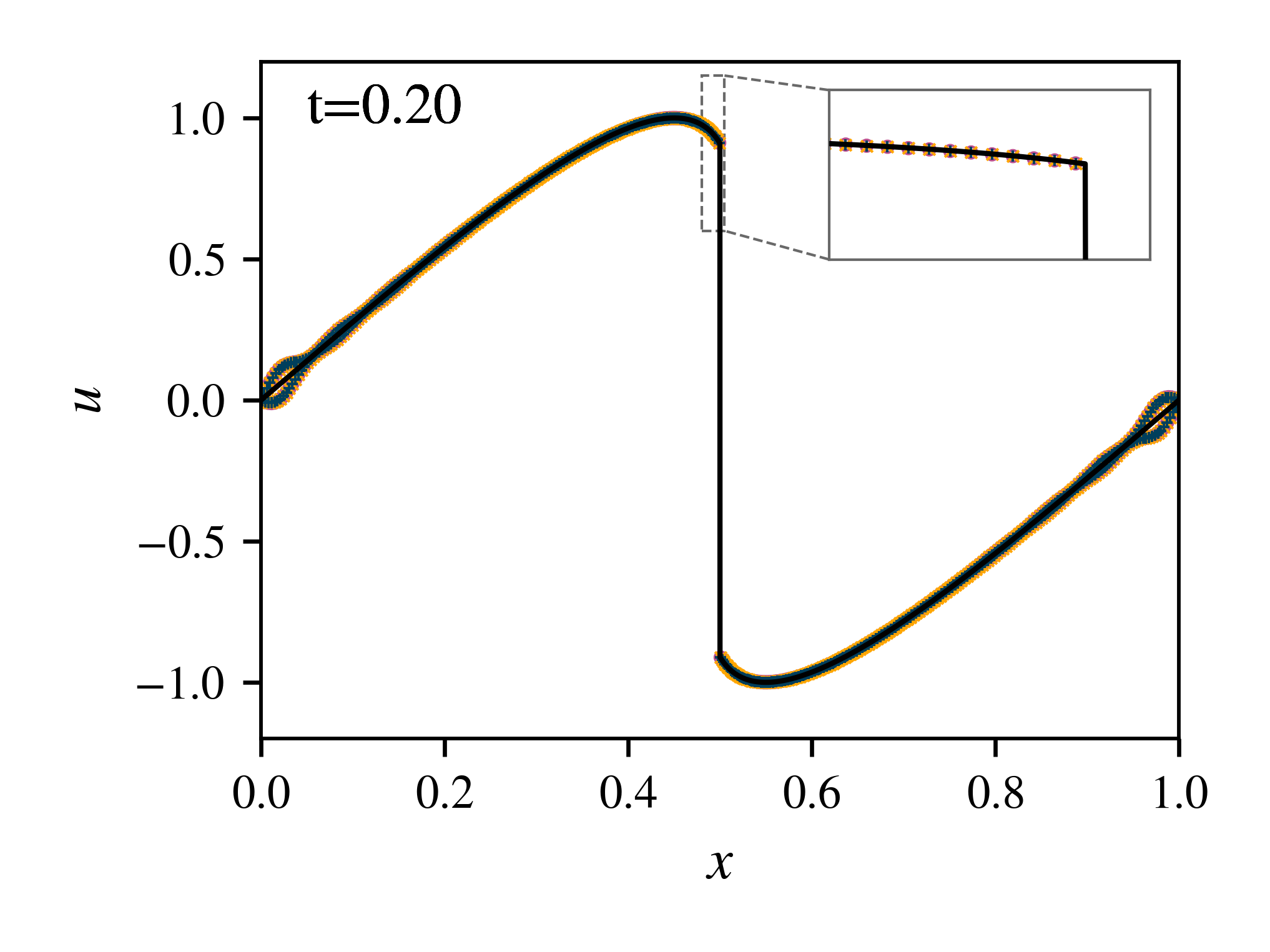}};
    \draw (2.8,-1.2) node {$(b)$};
\end{tikzpicture}
 \caption{Plots vs. $x$ of the velocity field $u$ obtained by using $(a)$ spectral relaxation and $(b)$ spectral purging methods with different kernels listed in the legend for the optimal parameter set, at time $t=0.2$, after the shock formation around $t_* \sim 0.159$. We compare these approximations with the exact solution (black line). The zoomed insets in each panel show that the behaviour of the approximated solution is the same across these kernel choices.}
 \label{fig:SPSR-difffilter}
\end{figure}
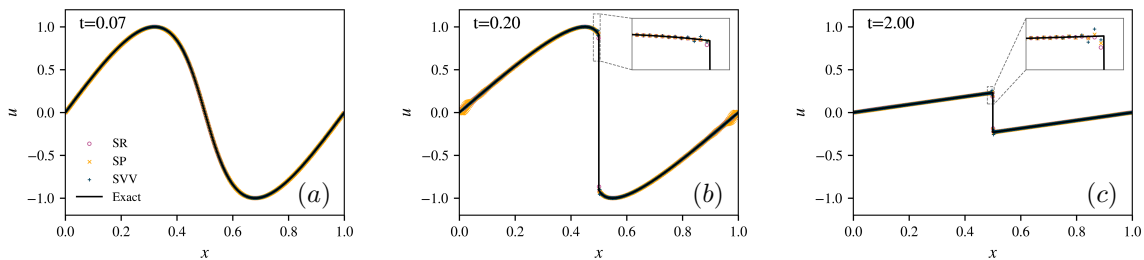
\begin{figure}[h!]
     \centering
        \begin{tikzpicture}
 	\draw (0,0) node[inner sep=0]{\includegraphics[width=0.32\linewidth]{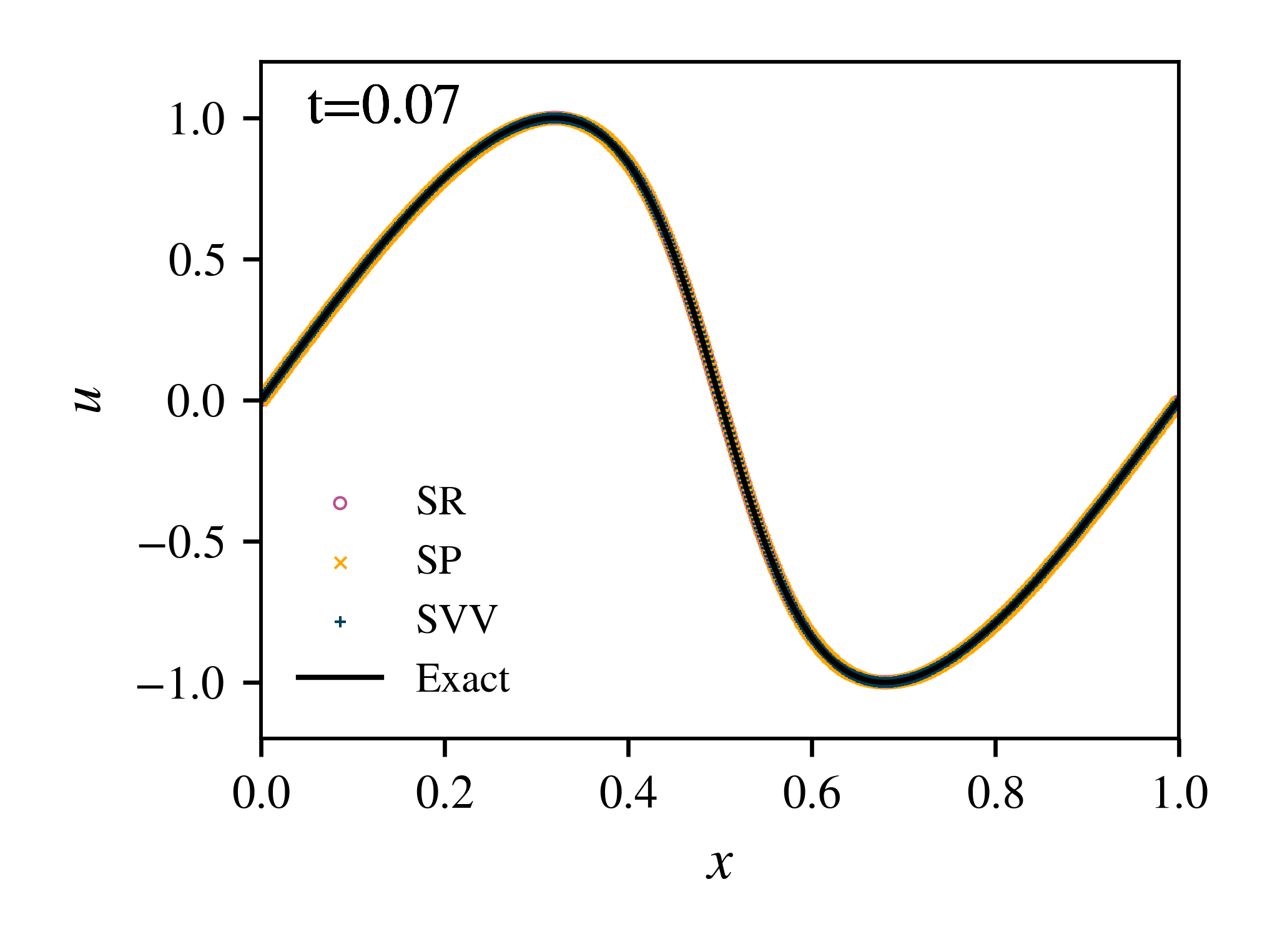}};
		\draw (1.8,-0.8) node {$(a)$};
	\end{tikzpicture}
         \begin{tikzpicture}
 	\draw (0,0) node[inner sep=0]{\includegraphics[width=0.32\linewidth]{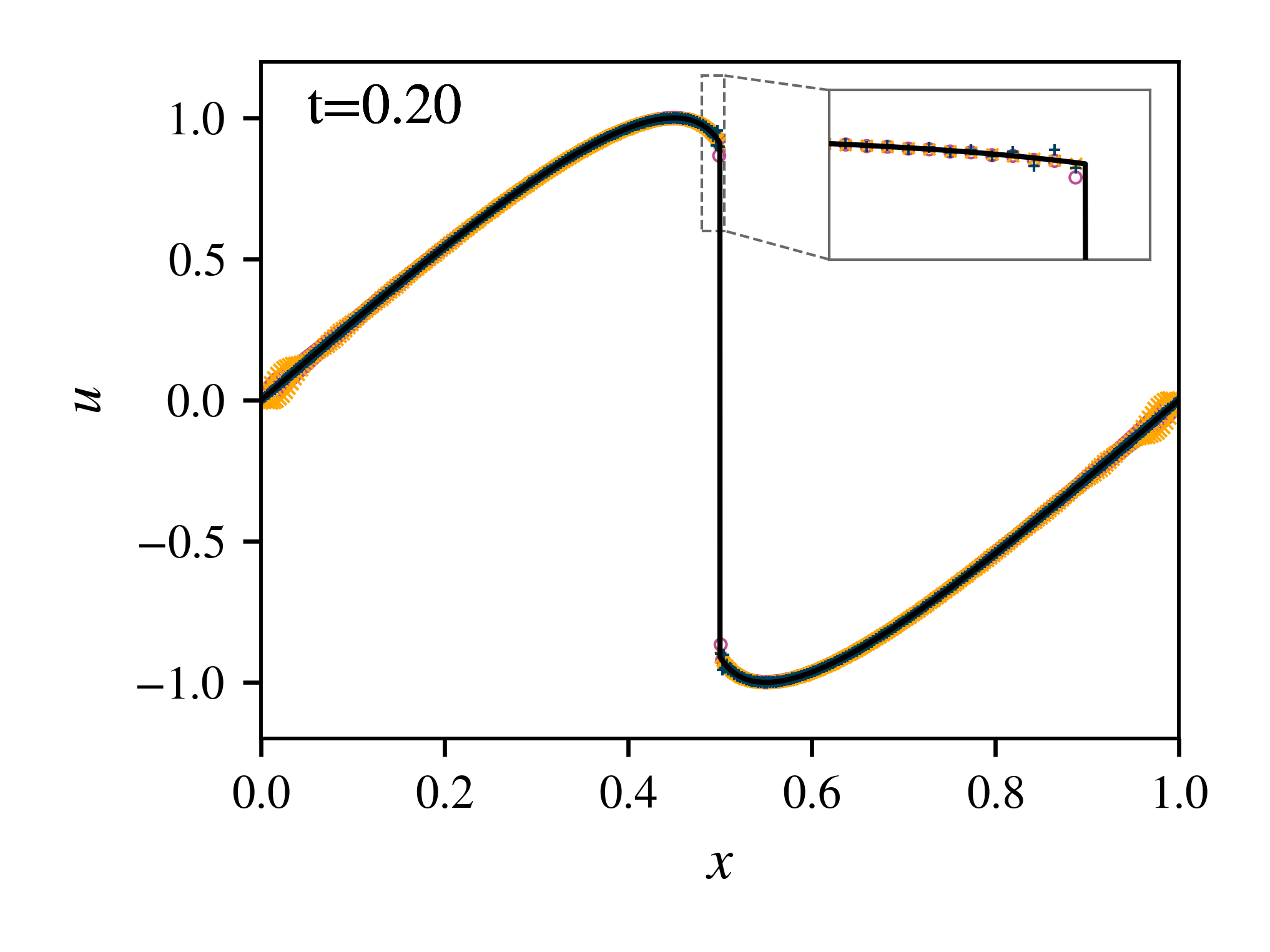}};
		\draw (1.8,-0.8) node {$(b)$};
	\end{tikzpicture}
         \begin{tikzpicture}
 	\draw (0,0) node[inner sep=0]{\includegraphics[width=0.32\linewidth]{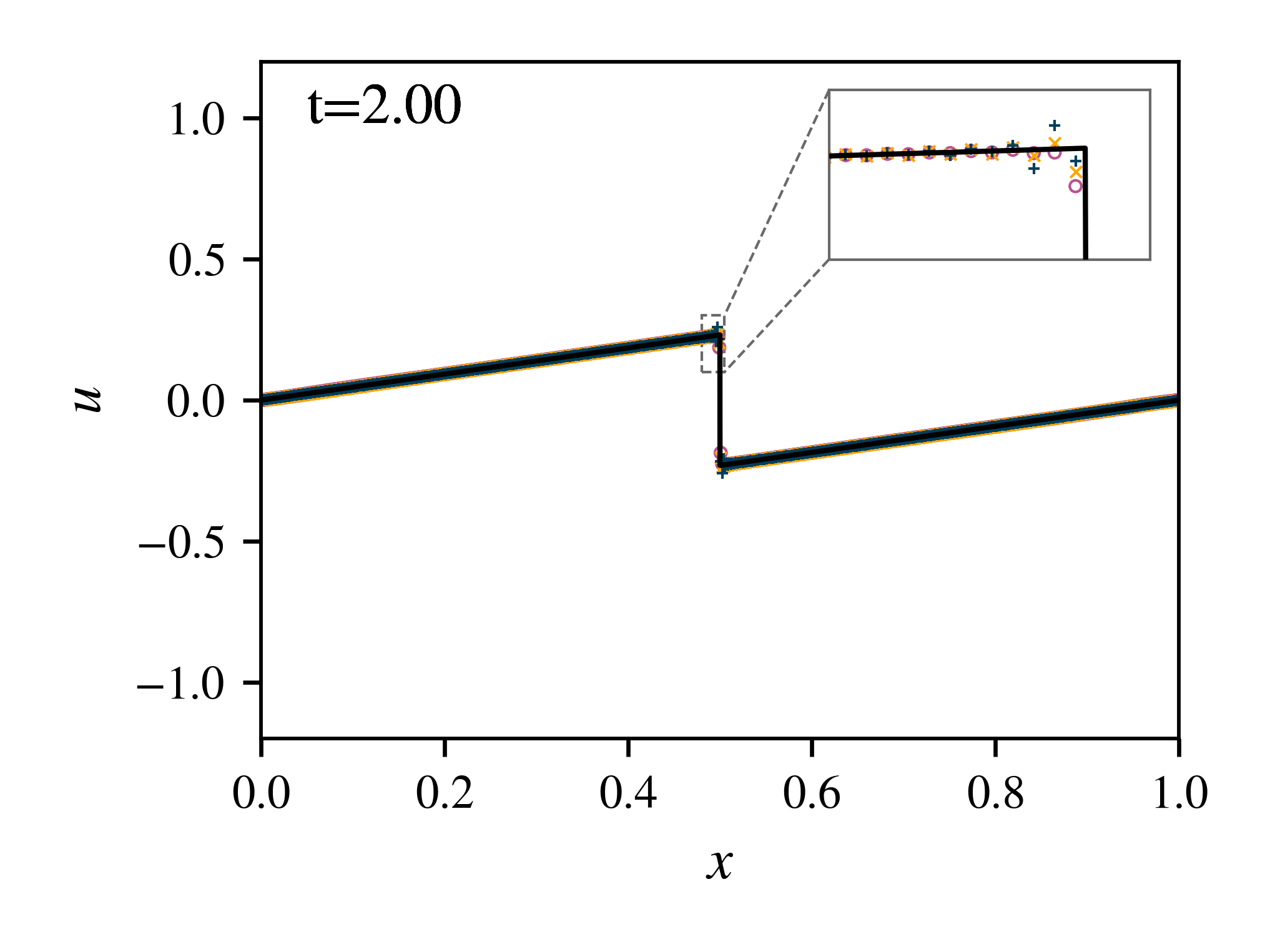}};
		\draw (1.8,-0.8) node {$(c)$};
	\end{tikzpicture}
   \caption{Plots vs. $x$ of the velocity field $u$ obtained by using the best SR-FeKo $(\alpha=0.7,\gamma=0.99)$ (pink circles) and SP-FeKo $(\alpha=0.65,\gamma=0.99)$ (orange crosses) approximations, compared with the widely used SVV approximations with $(\epsilon=N^{-1}, M=2\sqrt{N})$ (blue plusses), where $N_x =2N+1=615$. We describe the SVV filter in Appendix~\ref{app:filters}.}
     \label{fig:SRSPSVV}
 \end{figure}

\begin{figure}[h!]
     \centering
        \begin{tikzpicture}
 	\draw (0,0) node[inner sep=0]{\includegraphics[width=0.32\linewidth]{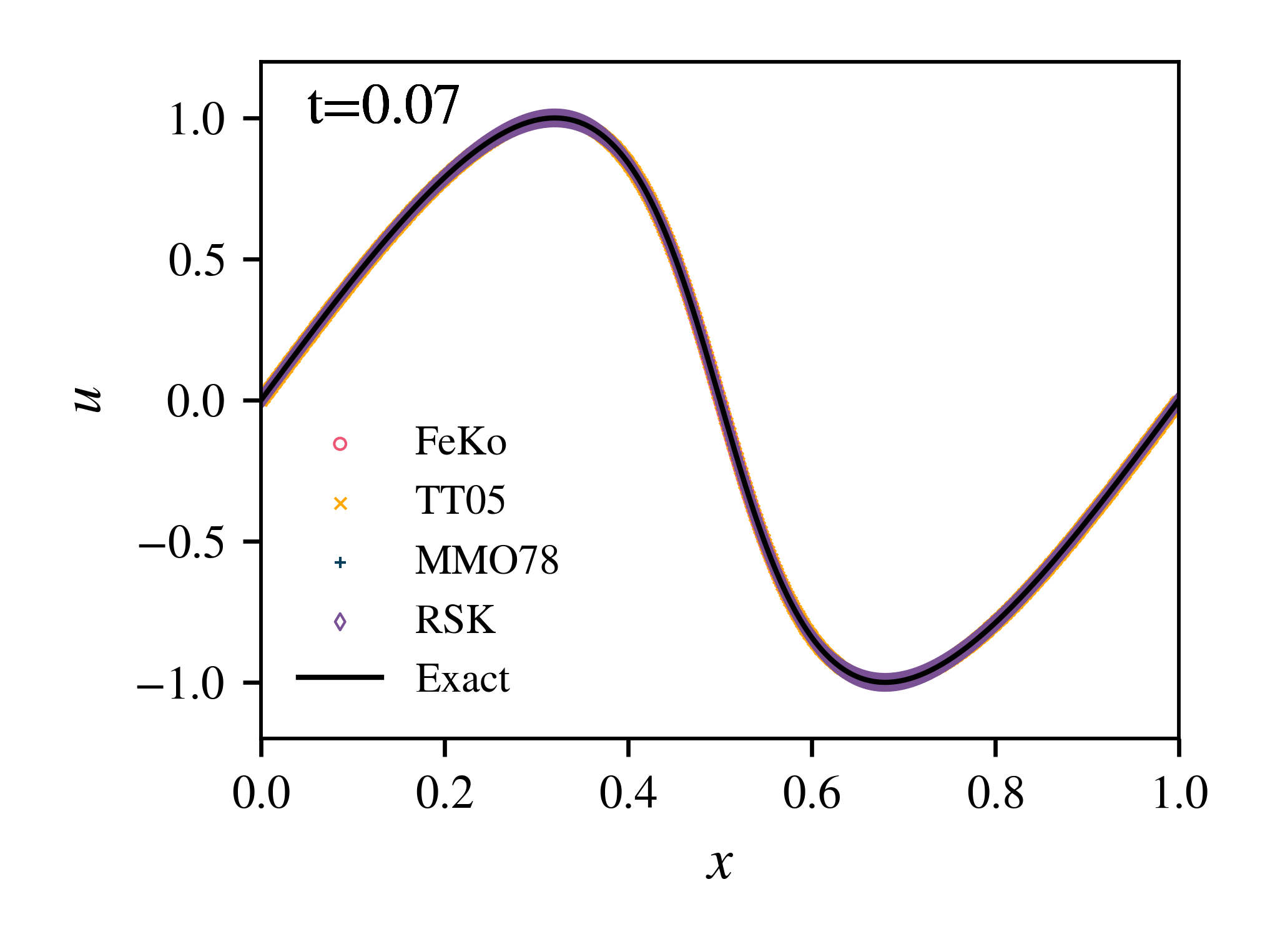}};
		\draw (1.8,-0.8) node {$(a)$};
	\end{tikzpicture}
         \begin{tikzpicture}
 	\draw (0,0) node[inner sep=0]{\includegraphics[width=0.32\linewidth]{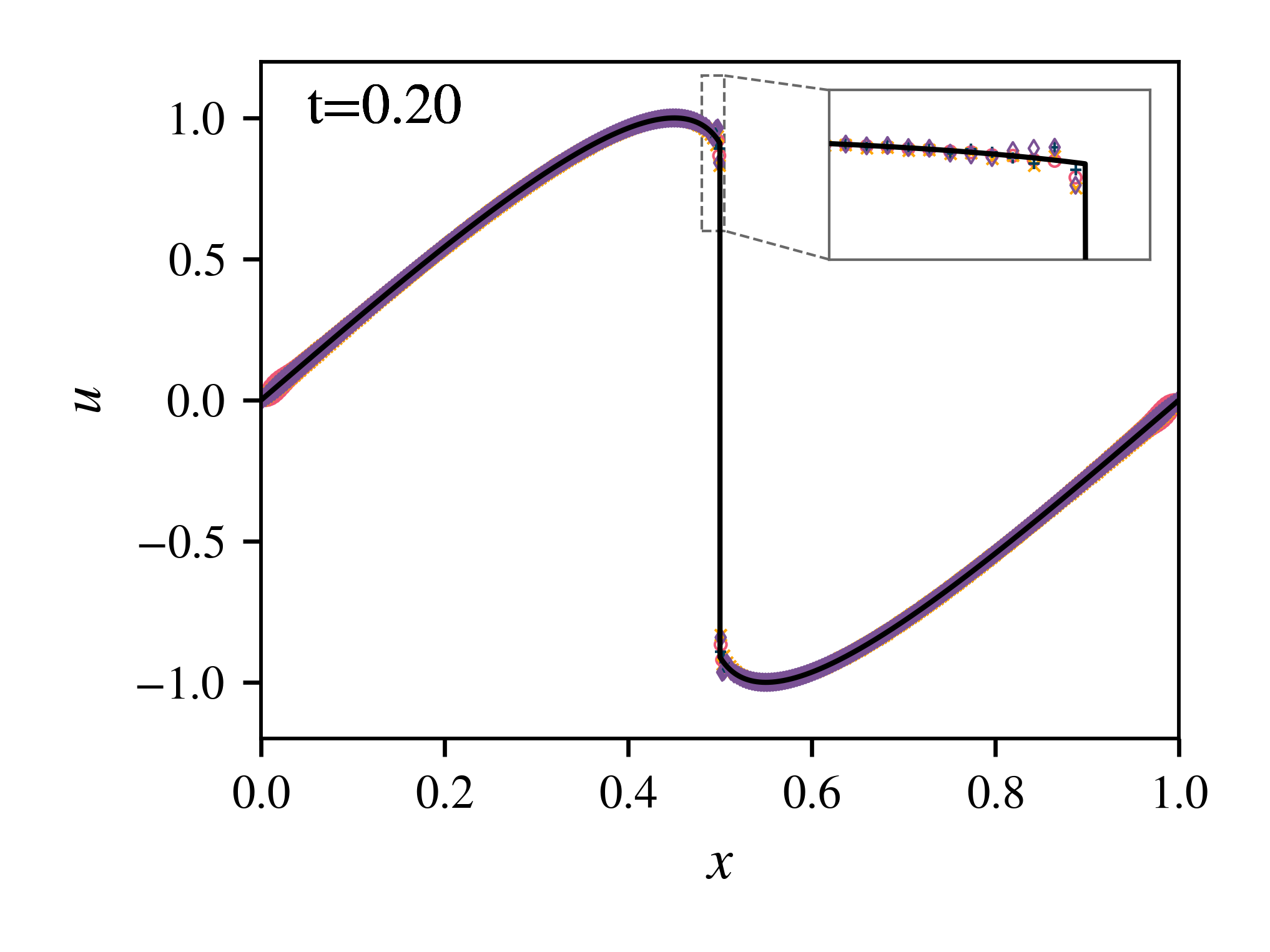}};
		\draw (1.8,-0.8) node {$(b)$};
	\end{tikzpicture}
         \begin{tikzpicture}
 	\draw (0,0) node[inner sep=0]{\includegraphics[width=0.32\linewidth]{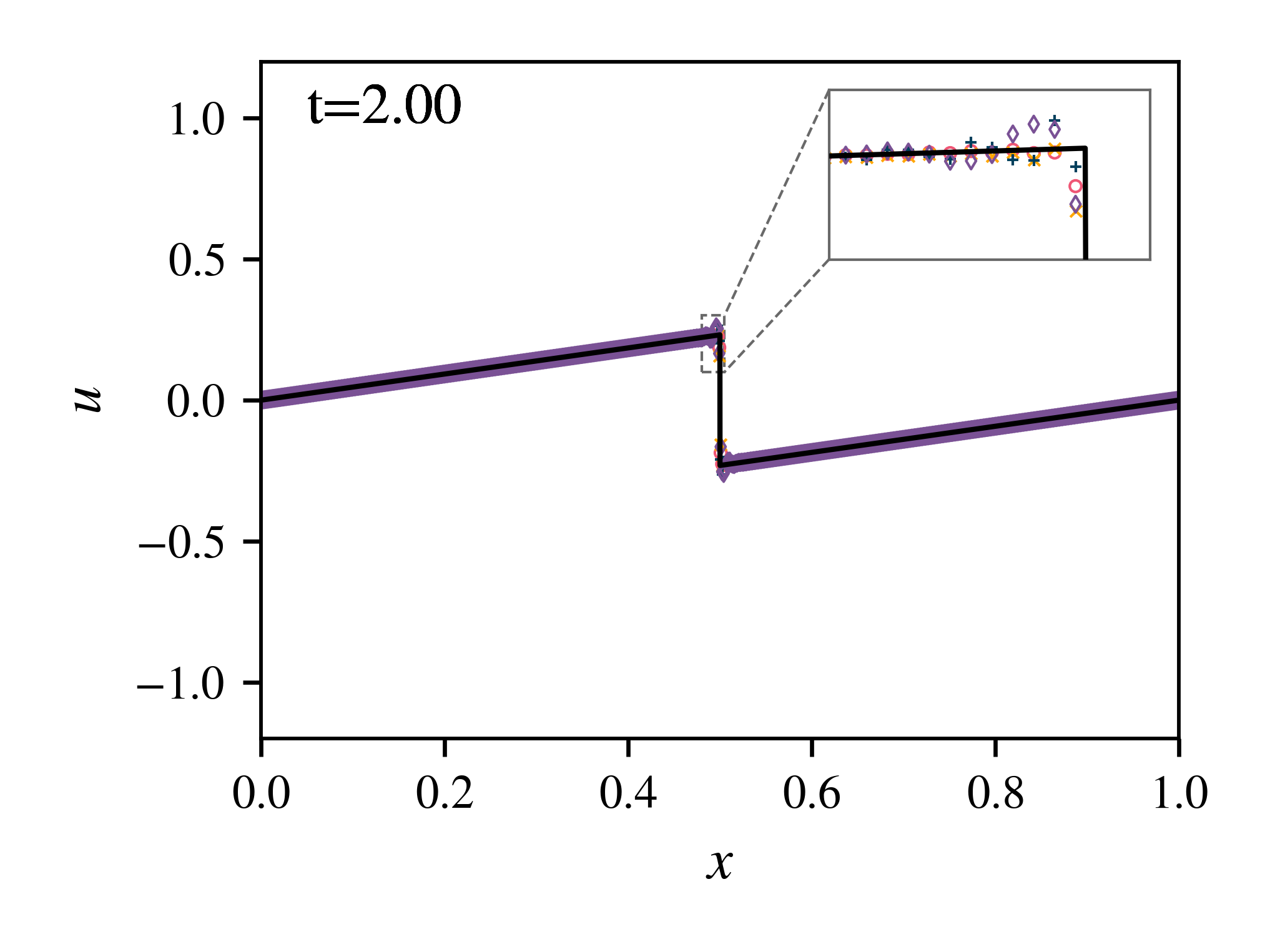}};
		\draw (1.8,-0.8) node {$(c)$};
	\end{tikzpicture}
   \caption{Plots vs. $x$ of the velocity field $u$ obtained by using SR-FeKo $(\alpha=0.70,\gamma=0.99)$ scheme, compared with the spectral relaxation scheme implemented using other filters which are not positive. We describe all filters used here in the Appendix~\ref{app:filters}.}
     \label{fig:SR-diffkernels}
 \end{figure}
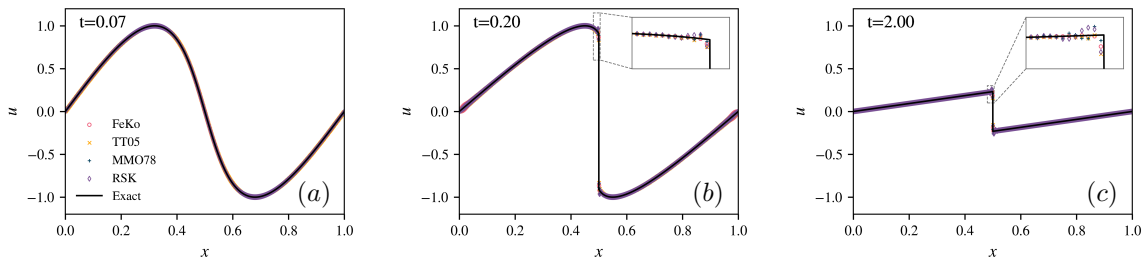

In Fig.~\ref{fig:SR-diffkernels}, we compare the SR-FeKo approximation with SR approximations that were obtained for different kernels $K_m(x)$ (used in  $\mathcal{M}_m u_N$ in Eq.~\eqref{eq:burg_x_sr}) that are prominent in the literature. In particular, we use the adaptive filters given by~\cite{tadmor2005adaptive} (TT05), the smoothing functions used by~\cite{majda1978fourier} (MMO78), and the Gaussian regularised Shannon kernel of~\cite{sun2002practical} (RSK). SR-FeKo gives the best visual convergence to the exact solution (see the zoomed inset in panels $(b,c)$) The approximations become more inaccurate as we move from SR-TT05, to SR-MMO78, and SR-RSK (in that order).

\subsection{Role of initial conditions in the choice of optimal \texorpdfstring{$\alpha$}{alpha}}
\label{subsec:Results_burg_ics}

In previous and subsequent sections, we use the novel spectral schemes without a dealiasing mechanism as the dissipation induced by these schemes usually suffices to control mild aliasing instabilities, so we can circumvent the cost of dealiasing for Eqs.~\eqref{eq:burg}. In this section, we demonstrate that, for a given PDE, different initial conditions can lead to different growth-rates for aliasing errors. Thus, the choice of kernels and the parameters $(\alpha,\gamma)$ cannot be universal; these have to be obtained through careful parameter exploration for every initial-boundary value problem.  

We investigate a shifted single-mode initial condition with SR-FeKo scheme; the implementation is as discussed in Section~\ref{sec:methods}. For clarity, we label the old and new initial conditions as below:
\begin{align}
    \texttt{IC0:    }& u_0(x) = \sin(2 \pi x) \\
    \texttt{IC1:    }& u_0(x) = \sin(2 \pi x - \pi/2) 
\end{align}
where, we have periodic boundaries for $ x \in [0,1] $.  
 
\begin{figure}[h!]
     \centering
        \begin{tikzpicture}
 	\draw (0,0) node[inner sep=0]{\includegraphics[width=0.32\linewidth]{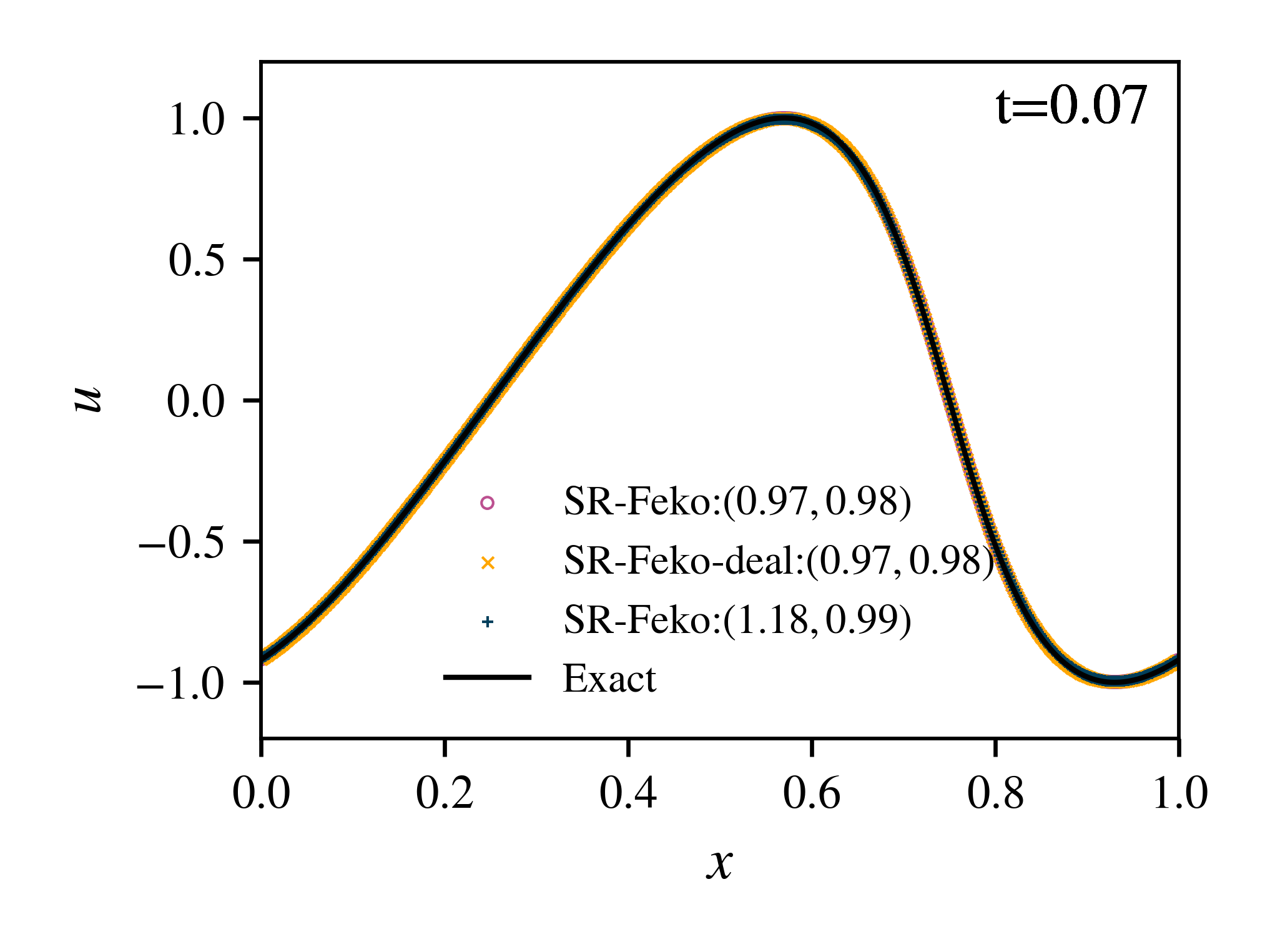}};
		\draw (-1,1.4) node {$(a)$};
	\end{tikzpicture}
         \begin{tikzpicture}
 	\draw (0,0) node[inner sep=0]{\includegraphics[width=0.32\linewidth]{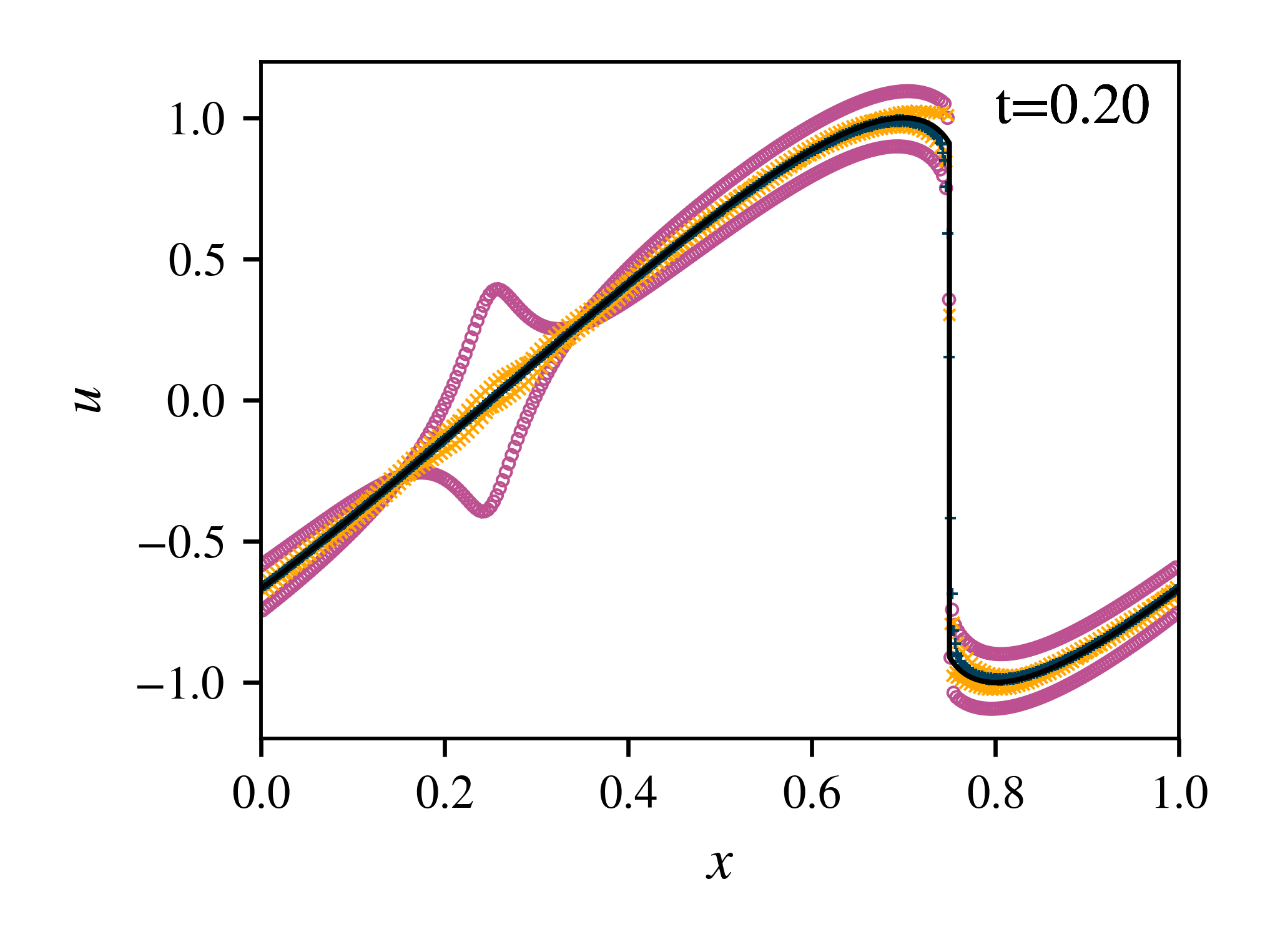}};
		\draw (-1,1.4) node {$(b)$};
	\end{tikzpicture}
         \begin{tikzpicture}
 	\draw (0,0) node[inner sep=0]{\includegraphics[width=0.32\linewidth]{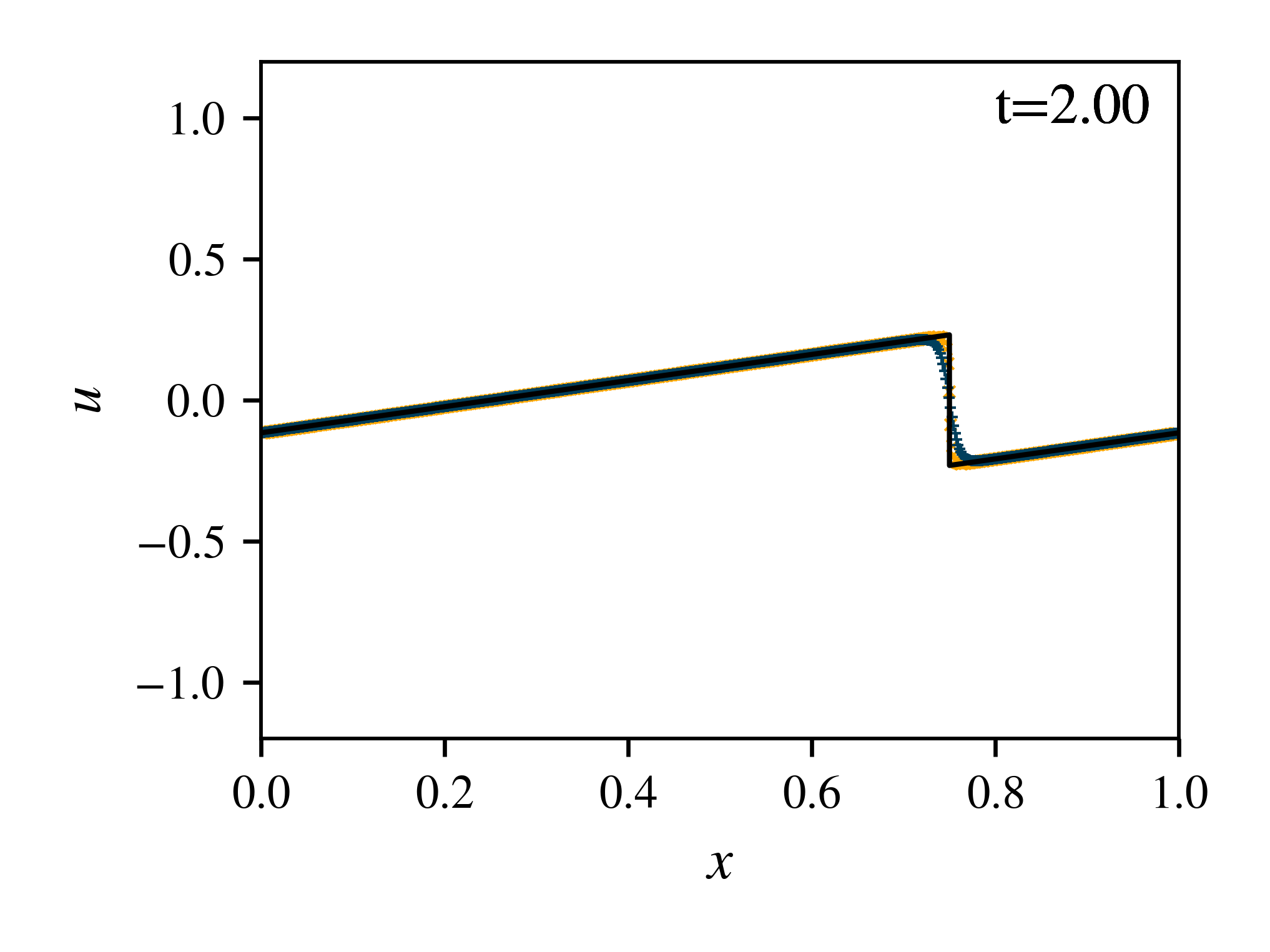}};
		\draw (-1,1.4) node {$(c)$};
	\end{tikzpicture}
   \caption{Plots vs. $x$ of the velocity field $u$ at three representative times, where the approximations were obtained using SR-FeKo and no dealiasing (plain) for $(\alpha=0.97,\gamma=0.98)$ (pink circles), $(\alpha=1.18,\gamma=0.99)$ (blue plusses) compared with SR-FeKo-deal implemented with $2/3$-dealiasing for $(\alpha=0.97,\gamma=0.98)$ (orange crosses). We show the exact solution (black line) for comparison. For stability of the plain SR approximations, we need values of $\alpha$ that are higher than the ones used with the dealiased approximations for this initial condition.}
     \label{fig:SR-IC1}
 \end{figure}
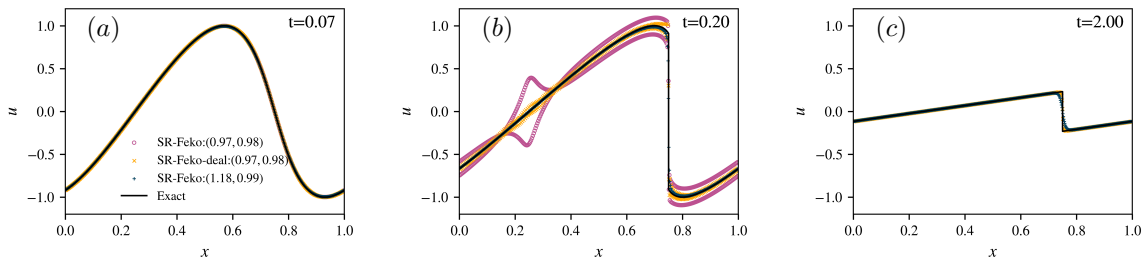

In Fig.~\ref{fig:SR-IC1}, we plot the non-dealiased (plain) SR-FeKo approximation for initial data \texttt{IC1} with $(\alpha=0.97,\gamma=0.98)$. Right after the shock formation, in panel $(b)$, the plain SR-FeKo approximation (pink circles) develops aliasing errors in the smooth region which eventually grow in time and destabilise the scheme; at long times, in panel $(c)$, the solution has numerically blown up and is no longer visible in the range of the plot. Thus, SR-FeKo is unstable for \texttt{IC1} with $\alpha<1$; this is unlike what we report for \texttt{IC0} in the previous section: SR-FeKo provides a stable and convergent approximation for \texttt{IC0} with values of $\alpha$ in $[0.7,2[$ and $\gamma \in [0.9,1)$; for $\alpha>1$, the scheme is over-dissipative but still converges. If we allow $\alpha>1$ for \texttt{IC1}, SR-FeKo yields a convergent approximation for $(\alpha=1.18,\gamma=0.99)$. The corresponding convergence analysis is given in Table~\ref{tab:4_cs_ic1-conv-nodeal}. In Fig.~\ref{fig:SR-IC1}, the plain SR-FeKo approximation (blue plusses), with the new optimal parameter set, shows visual convergence to the exact solution and does not blow-up at long times. Convergence analysis, performed in~\cite{bessemain}, proves that, for a positive kernel $K_m(x)$ and parameters $(\alpha,\gamma)$ such that ${\alpha}/{2} < \gamma<1$ and $1<\alpha<2$, there is convergence in $L^2$ norm of the approximate solution $u_N$ to the exact entropic solution of the 1D inviscid Burgers equation as $N\rightarrow\infty$. In fact, the condition $\alpha>1$ is a sufficient condition for the aliasing error to be controlled or dominated by the dissipation induced by the relaxation term. This is confirmed by Fig.~\ref{fig:SR-IC1}; for this \texttt{IC}, $\alpha>1$ is optimal for obtaining a convergent SR approximation without dealiasing; but this does not mean that for $\alpha <1$ there is never any convergence for the SR method without dealiasing. As we observe in Fig.~\ref{fig:SR-IC1}, it can depend on the initial condition. Indeed, particular initial conditions (like \texttt{IC0}) may exist for which supplementary conditions such as compensation mechanisms can make the SR method converge for a wider range of the parameter $\alpha$. The additional compensation mechanisms are not taken into account by the general convergence analysis which covers all initial conditions; they deserve to be studied in greater detail from a mathematical point of view.

\begin{figure}[h!]
 \centering
 \begin{tikzpicture}
\draw (0,0) node[inner sep=0]{\includegraphics[width=0.4\linewidth]{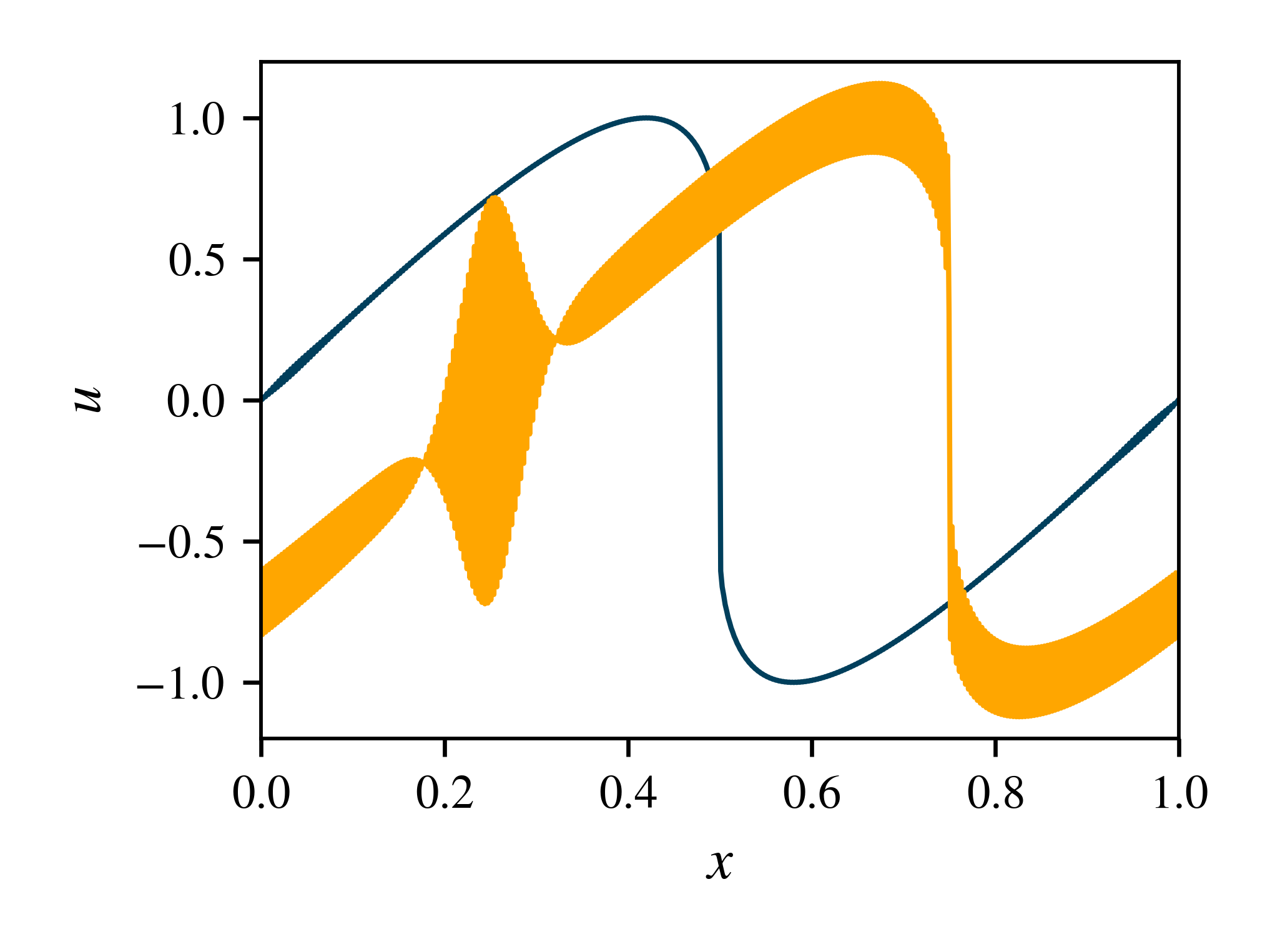}};
    \draw (2.0,1.6) node {$(a)$};
\end{tikzpicture}
\begin{tikzpicture}
\draw (0,0) node[inner sep=0]{\includegraphics[width=0.4\linewidth]{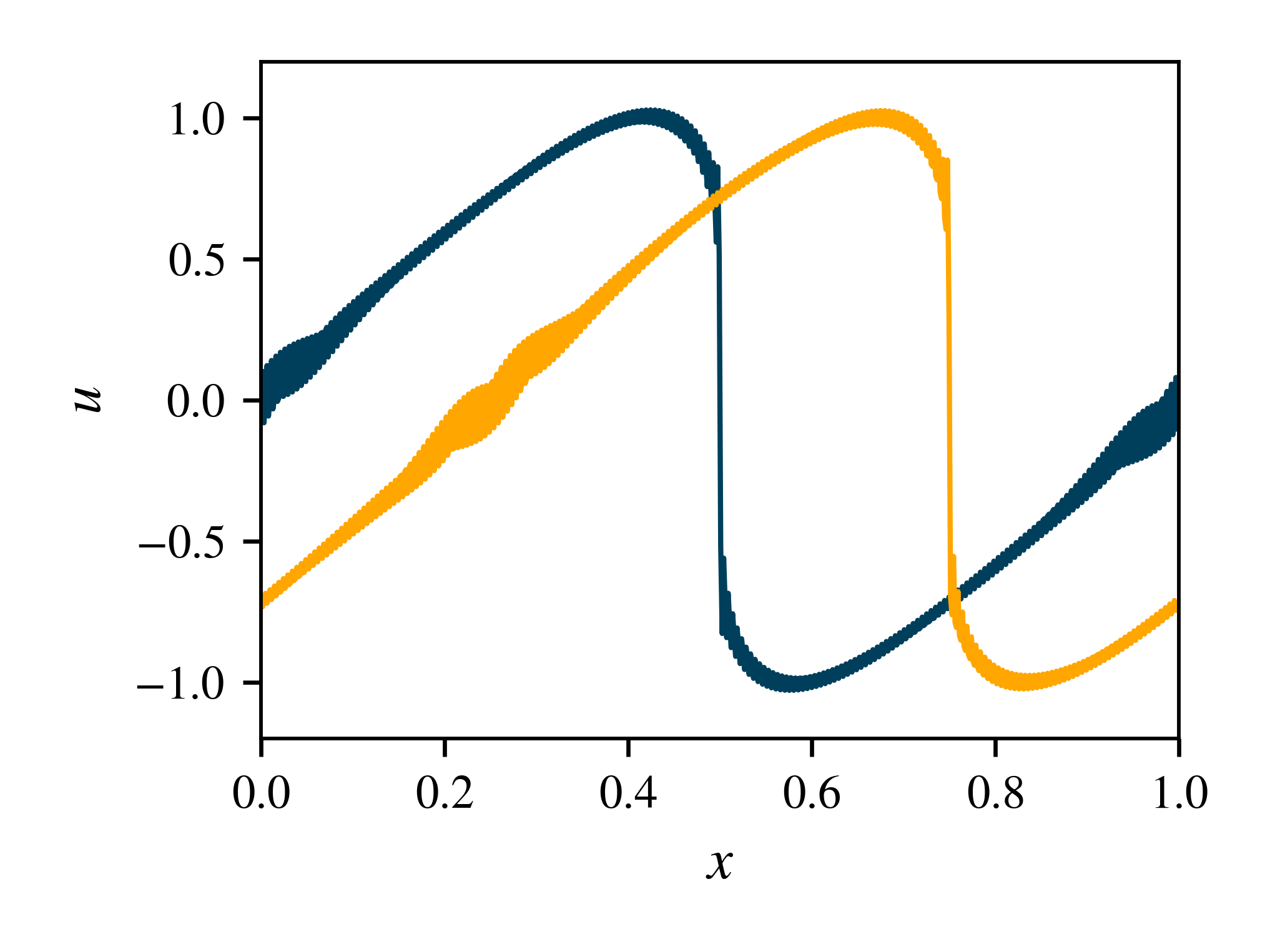}};
    \draw (2.0,1.6) node {$(b)$};
\end{tikzpicture}
 \caption{Plots vs. $x$ of the velocity field $u(x,t)$ obtained by using PPS schemes for the initial conditions \texttt{IC0} and \texttt{IC1} at time $t=0.17$ just after the time of shock formation for this model. In panel $(a)$, dealiasing is not used and aliasing errors build at different rates for the different \texttt{ICs}. In panel $(b)$ we use $2/3$-dealiasing: for \texttt{IC0}, the PPS approximation becomes worse and develops tygers. For \texttt{IC1}, the solution is more well-behaved and the aliasing instability is bypassed. We use $N_x=615$ for all these runs.}
 \label{fig:PPS-deal}
\end{figure}

To expose the role of \texttt{ICs} in producing aliasing errors of varying strengths, we look at the plain PPS approximations for both \texttt{IC0} and \texttt{IC1}. In panel $(a)$ of Fig.~\ref{fig:PPS-deal}, we show the plain PPS for \texttt{IC0} in blue, and \texttt{IC1} in orange, at $t=0.17$ after the shock formation. The growth of aliasing error at the edge of the domain is very mild in \texttt{IC0}. For \texttt{IC1}, the errors accumulate around $x\simeq 0.25$ and grow to destabilise the PPS scheme before $t=0.2$. In panel $(b)$ of the same figure, we repeat PPS for both \texttt{ICs} with $2/3$-dealiasing. Now, the results are uniform across the \texttt{ICs} and tygers emerge in both cases. Thus, $2/3$-dealiasing plays an important role in stabilising and homogenising the Fourier pseudospectral scheme. However, for \texttt{IC0}, dealiasing is expensive and disadvantageous; the dealiased PPS approximation (blue line in panel $(b)$) is far-less accurate than the non-dealiased one (blue line in panel $(a)$). For \texttt{IC1}, dealiasing improves the stability and accuracy of the PPS approximation (compare orange lines for plain in panel $(a)$ and dealiased in $(b)$). When we use $2/3$-dealiasing in SR-FeKo for \texttt{IC1}, the stability and convergence of the scheme is ameliorated. In Table~\ref{tab:3_cs_ic1-conv-deal}, we present numerical convergence analysis for the dealiased SR-FeKo approximation for $(\alpha=0.97,\gamma=0.98)$. Convergence is improved compared to the plain SR-FeKo approximation for $(\alpha=1.18,\gamma=0.99)$ discussed above (see Table~\ref{tab:4_cs_ic1-conv-nodeal}). Furthermore, in Fig.~\ref{fig:SR-IC1}, stability is restored for the dealiased SR-FeKo approximation (orange crosses) for the same parameters as the non-dealiased one (pink circles). Dealiased SR schemes require a smaller value of $\alpha$ to achieve stability, when compared to their non-dealiased counterparts. Low values of $\alpha$ imply less diffusion and higher order of convergence of the scheme since the order of convergence $\sim (2-\alpha)$ (see Fig.~\ref{fig:t1_conv_sr_alpha}). Thus, dealiasing makes it possible to obtain a stable SR scheme that is less diffusive and shows a higher rate of convergence.

\begin{figure}[h!]
     \centering
        \begin{tikzpicture}
 	\draw (0,0) node[inner sep=0]{\includegraphics[width=0.32\linewidth]{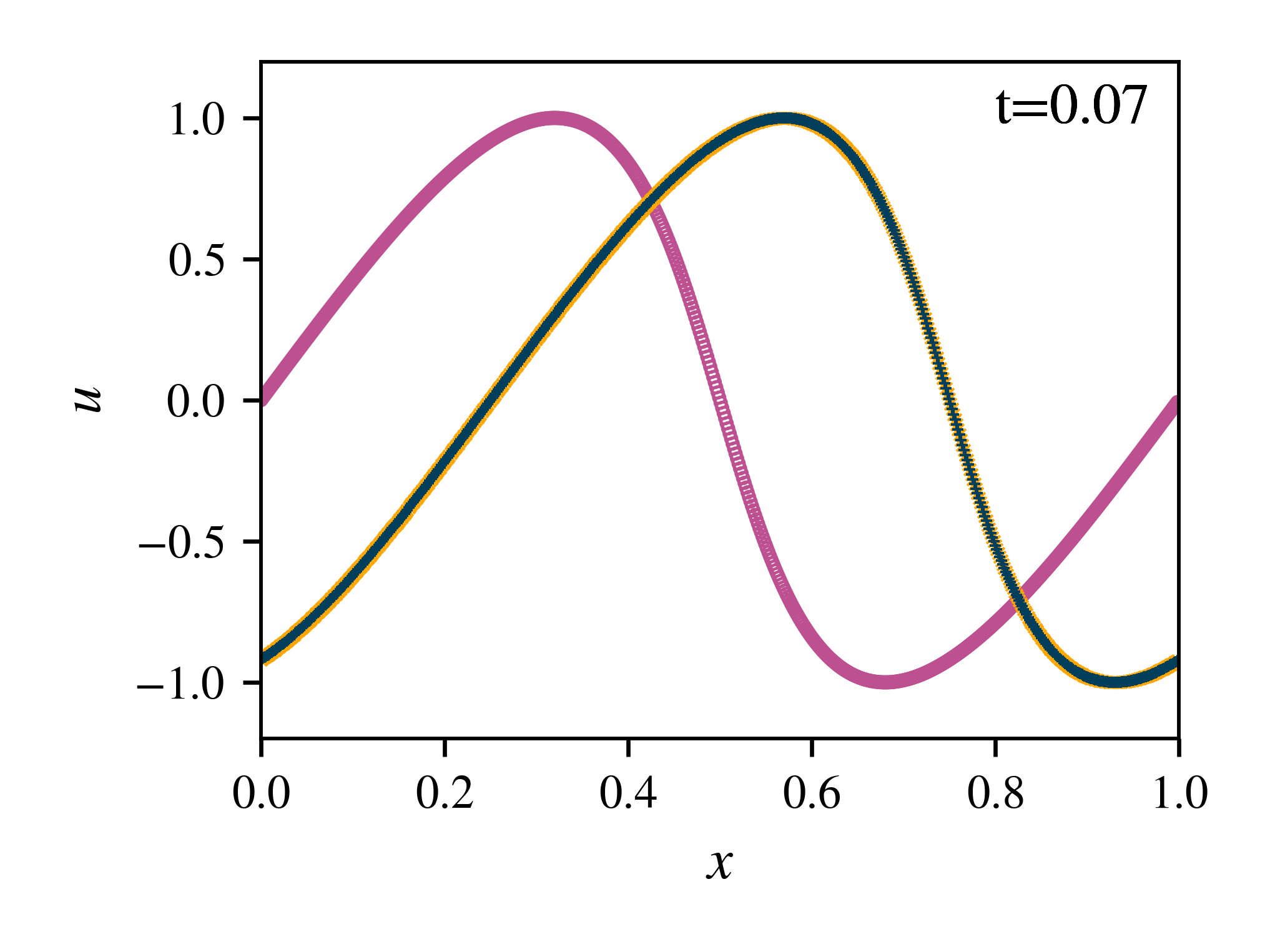}};
		\draw (-1,-0.8) node {$(a)$};
	\end{tikzpicture}
         \begin{tikzpicture}
 	\draw (0,0) node[inner sep=0]{\includegraphics[width=0.32\linewidth]{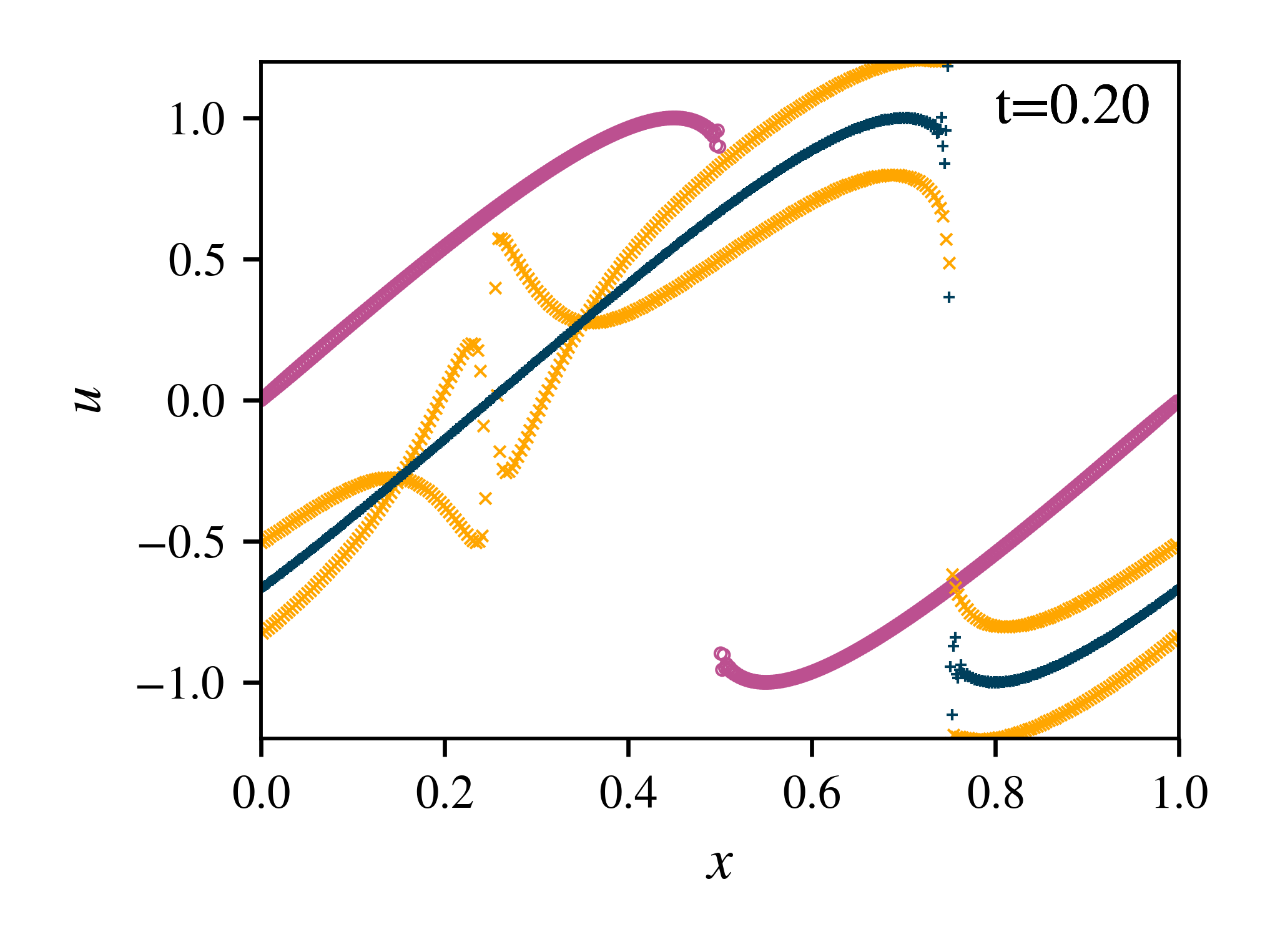}};
		\draw (-1,-0.8) node {$(b)$};
	\end{tikzpicture}
         \begin{tikzpicture}
 	\draw (0,0) node[inner sep=0]{\includegraphics[width=0.32\linewidth]{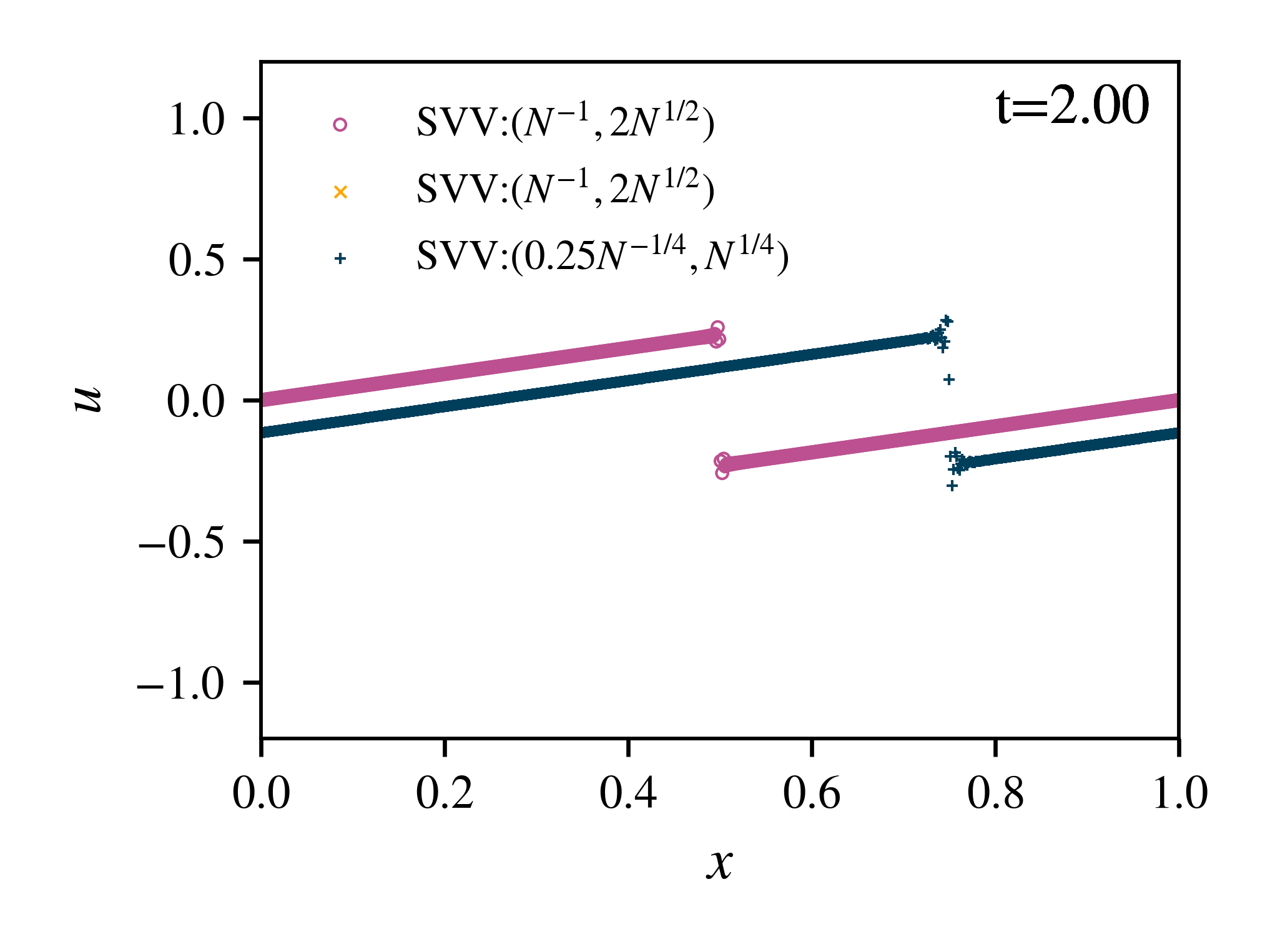}};
		\draw (-1,-0.8) node {$(c)$};
	\end{tikzpicture}
        \caption{Plots vs. $x$ of the velocity field $u (x,t)$ for three representative times, where we use the SVV method with parameters $(\epsilon,M)$ listed in the legend for different initial configurations. We do not use $2/3$-dealiasing in these runs. The spatial resolution is $N_x=615$. }
    \label{fig:svv-deal-conv}
\end{figure}
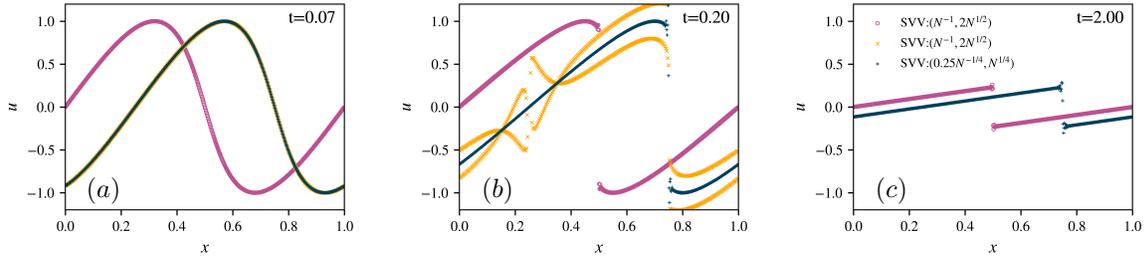

This effect is also observed for SVV schemes which circumvent dealiasing as well. In Fig.~\ref{fig:svv-deal-conv}, SVV approximation of \texttt{IC0} (pink circles) is stable and convergent for $(\epsilon=N^{-1},M=2N^{1/2})$. When we use the same parameters for \texttt{IC1} (orange crosses), the resulting SVV approximation destabilises. The more dissipative choice of $(\epsilon=0.25 N^{-1/4},M=N^{1/4})$ is optimal for this case (shown as blue plusses) and then restores the stability and convergence of the scheme. 

\begin{table}[h!]
\centering
 \begin{tabularx}{\linewidth}{@{}|YY|YY|YY|YY|@{}}
    \hline 
    \rule{0pt}{2.5ex}    
      SR & &  \multicolumn{2}{>{\hsize=\dimexpr1\hsize+4\tabcolsep+1\arrayrulewidth\relax}X|}{\centering $t=0.07$} & \multicolumn{2}{>{\hsize=\dimexpr2\hsize+4\tabcolsep+1\arrayrulewidth\relax}X|}{\centering $t=0.20$} &
      \multicolumn{2}{>{\hsize=\dimexpr2\hsize+4\tabcolsep+1\arrayrulewidth\relax}X|}{\centering $t=2.00$} \\

      FeKo & $N_x$ & Error & Order & Error & Order & Error & Order \\ 
        \hline        
            \rule{0pt}{3ex}
 \multirow{8}{*}{ $L^1$ error} &39	& $1.8\cdot 10^{-2}$ &		&	$6.0\cdot 10^{-2}$	& 	 &	$2.9\cdot 10^{-2}$	&  \\
& 65	& $1.3\cdot 10^{-2}$ &	$0.62$	&	$4.4\cdot 10^{-2}$	& $0.60$ &	$2.1\cdot 10^{-2}$	&  $0.65$  \\
& 123	& $8.3\cdot 10^{-3}$ &	$0.76$	&	$2.8\cdot 10^{-2}$	& $0.72$ &	$1.2\cdot 10^{-2}$	&  $0.79$  \\
& 205	& $5.6\cdot 10^{-3}$ &	$0.76$	&	$1.9\cdot 10^{-2}$	& $0.75$ &	$8.4\cdot 10^{-3}$	&  $0.78$  \\
& 615	& $2.4\cdot 10^{-3}$ &	$0.77$	&	$7.9\cdot 10^{-3}$	& $0.80$ &	$3.6\cdot 10^{-3}$	&  $0.78$  \\
& 1599	& $1.1\cdot 10^{-3}$ &	$0.80$	&	$3.5\cdot 10^{-3}$	& $0.84$ &	$1.7\cdot 10^{-3}$	&  $0.80$  \\
& 2665	& $7.4\cdot 10^{-4}$ &	$0.80$	&	$2.3\cdot 10^{-3}$	& $0.83$ &	$1.1\cdot 10^{-3}$	&  $0.80$  \\
& 7995	& $3.1\cdot 10^{-4}$ &	$0.80$	&	$9.6\cdot 10^{-4}$	& $0.81$ &	$4.6\cdot 10^{-4}$	&  $0.80$  \\
           \hline                                                                 
             \rule{0pt}{3ex}
 \multirow{8}{*}{$L^2$ error} &39	& $2.2\cdot 10^{-2}$ &		&	$1.3\cdot 10^{-1}$	& 	 &	$5.9\cdot 10^{-2}$	&  \\
& 65	& $1.6\cdot 10^{-2}$ &	$0.61$	&	$1.1\cdot 10^{-1}$	& $0.33$ &	$5.0\cdot 10^{-2}$	&  $0.32$ \\
& 123	& $1.0\cdot 10^{-2}$ &	$0.76$	&	$8.7\cdot 10^{-2}$	& $0.42$ &	$3.9\cdot 10^{-2}$	&  $0.39$ \\
& 205	& $6.8\cdot 10^{-3}$ &	$0.76$	&	$6.9\cdot 10^{-2}$	& $0.44$ &	$3.2\cdot 10^{-2}$	&  $0.38$ \\
& 615	& $2.9\cdot 10^{-3}$ &	$0.77$	&	$4.2\cdot 10^{-2}$	& $0.46$ &	$2.1\cdot 10^{-2}$	&  $0.39$ \\
& 1599	& $1.4\cdot 10^{-3}$ &	$0.79$	&	$2.8\cdot 10^{-2}$	& $0.43$ &	$1.4\cdot 10^{-2}$	&  $0.40$ \\
& 2665	& $9.0\cdot 10^{-4}$ &	$0.79$	&	$2.3\cdot 10^{-2}$	& $0.41$ &	$1.2\cdot 10^{-2}$	&  $0.40$ \\
& 7995	& $3.8\cdot 10^{-4}$ &	$0.80$	&	$1.5\cdot 10^{-2}$	& $0.40$ &	$7.6\cdot 10^{-3}$	&  $0.40$ \\
          \hline                                                                   
    \end{tabularx}
    
    \caption{Convergence analysis for SR-FeKo approximations using  $(\alpha=1.18,\gamma=0.99)$ at three representative times, $t=0.07,0.2$ and $2.0$. Here, the second row contains convergence analysis based on the  $L^1$ error and the third row, for the $L^2$ error. All sets of runs do not employ $2/3$-dealiasing.} 
    \label{tab:4_cs_ic1-conv-nodeal}
\end{table}%

\begin{table}[h!]
\centering
 \begin{tabularx}{\linewidth}{@{}|YY|YY|YY|YY|@{}}
    \hline 
    \rule{0pt}{2.5ex}    
      SR & &  \multicolumn{2}{>{\hsize=\dimexpr1\hsize+4\tabcolsep+1\arrayrulewidth\relax}X|}{\centering $t=0.07$} & \multicolumn{2}{>{\hsize=\dimexpr2\hsize+4\tabcolsep+1\arrayrulewidth\relax}X|}{\centering $t=0.20$} &
      \multicolumn{2}{>{\hsize=\dimexpr2\hsize+4\tabcolsep+1\arrayrulewidth\relax}X|}{\centering $t=2.00$} \\

      FeKo & $N_x$ & Error & Order & Error & Order & Error & Order \\ 
        \hline        
            \rule{0pt}{3ex}
\multirow{8}{*}{ $L^1$ error} & 39   & $ 1.1\cdot 10^{-2}$ &         &  $7.6\cdot 10^{-2}$ &       & $1.6\cdot 10^{-2}$ &     \\           
&  65   & $ 7.0\cdot 10^{-3}$ &  $0.91$  &  $5.3\cdot 10^{-2}$ &  $0.71$ & $1.0\cdot 10^{-2}$ & $0.89$ \\
&  123  & $ 3.8\cdot 10^{-3}$ &  $0.97$  &  $4.1\cdot 10^{-2}$ &  $0.41$ & $5.9\cdot 10^{-3}$ & $0.91$ \\
&  205  & $ 2.4\cdot 10^{-3}$ &  $0.91$  &  $3.3\cdot 10^{-2}$ &  $0.43$ & $3.7\cdot 10^{-3}$ & $0.89$ \\
&  615  & $ 8.1\cdot 10^{-4}$ &  $0.98$  &  $2.0\cdot 10^{-2}$ &  $0.45$ & $1.3\cdot 10^{-3}$ & $0.94$ \\
&  1599 & $ 3.1\cdot 10^{-4}$ &  $0.99$  &  $1.2\cdot 10^{-2}$ &  $0.52$ & $5.3\cdot 10^{-4}$ & $0.96$ \\
&  2665 & $ 1.9\cdot 10^{-4}$ &  $0.98$  &  $9.3\cdot 10^{-3}$ &  $0.54$ & $3.2\cdot 10^{-4}$ & $0.96$ \\
&  7995 & $ 6.4\cdot 10^{-5}$ &  $0.99$  &  $5.1\cdot 10^{-3}$ &  $0.54$ & $1.1\cdot 10^{-4}$ & $0.95$ \\
           \hline                                                                   
             \rule{0pt}{3ex}
\multirow{8}{*}{$L^2$ error}  & 39	& $1.3 \cdot 10^{-2}$ & 		&$1.3 \cdot 10^{-1}$& 		&	$4.3 \cdot 10^{-2}$ & 	\\
 & 65	& $8.4 \cdot 10^{-3}$ & 	$0.90$	&	$9.9 \cdot 10^{-2}$ & 	$0.59$	&	$3.4 \cdot 10^{-2}$ & 	$0.47$ \\
 & 123	& $4.5 \cdot 10^{-3}$ & 	$0.97$	&	$7.3 \cdot 10^{-2}$ & 	$0.49$	&	$2.5 \cdot 10^{-2}$ & 	$0.49$ \\
 & 205	& $2.9 \cdot 10^{-3}$ & 	$0.91$	&	$5.7 \cdot 10^{-2}$ & 	$0.47$	&	$2.0 \cdot 10^{-2}$ & 	$0.45$ \\
 & 615	& $9.8 \cdot 10^{-4}$ & 	$0.98$	&	$3.5 \cdot 10^{-2}$ & 	$0.46$	&	$1.1 \cdot 10^{-2}$ & 	$0.49$ \\
 & 1599	& $3.8 \cdot 10^{-4}$ & 	$0.99$	&	$2.2 \cdot 10^{-2}$ & 	$0.49$	&	$7.2 \cdot 10^{-3}$ & 	$0.49$ \\
 & 2665	& $2.3 \cdot 10^{-4}$ & 	$0.98$	&	$1.7 \cdot 10^{-2}$ & 	$0.49$	&	$5.6 \cdot 10^{-3}$ & 	$0.49$ \\
 & 7995	& $7.8 \cdot 10^{-5}$ & 	$0.99$	&	$9.8 \cdot 10^{-3}$ & 	$0.49$	&	$3.2 \cdot 10^{-3}$ & 	$0.50$ \\
           \hline                                                                   
    \end{tabularx}
    \caption{Convergence analysis for SR-FeKo approximations using  $(\alpha=0.97,\gamma=0.98)$ at three representative times $t=0.07,0.2$ and $2.0$. Here, the second row contains convergence analysis based on the  $L^1$ error and the third row, for the $L^2$ error. We employ $2/3$-dealiasing which makes this computation stable for the choice of kernel and parameters quoted above.} 
    \label{tab:3_cs_ic1-conv-deal}
\end{table}%

\subsection{2x2 system: 1D shallow water equation}
\label{subsec:Results_shalwat}

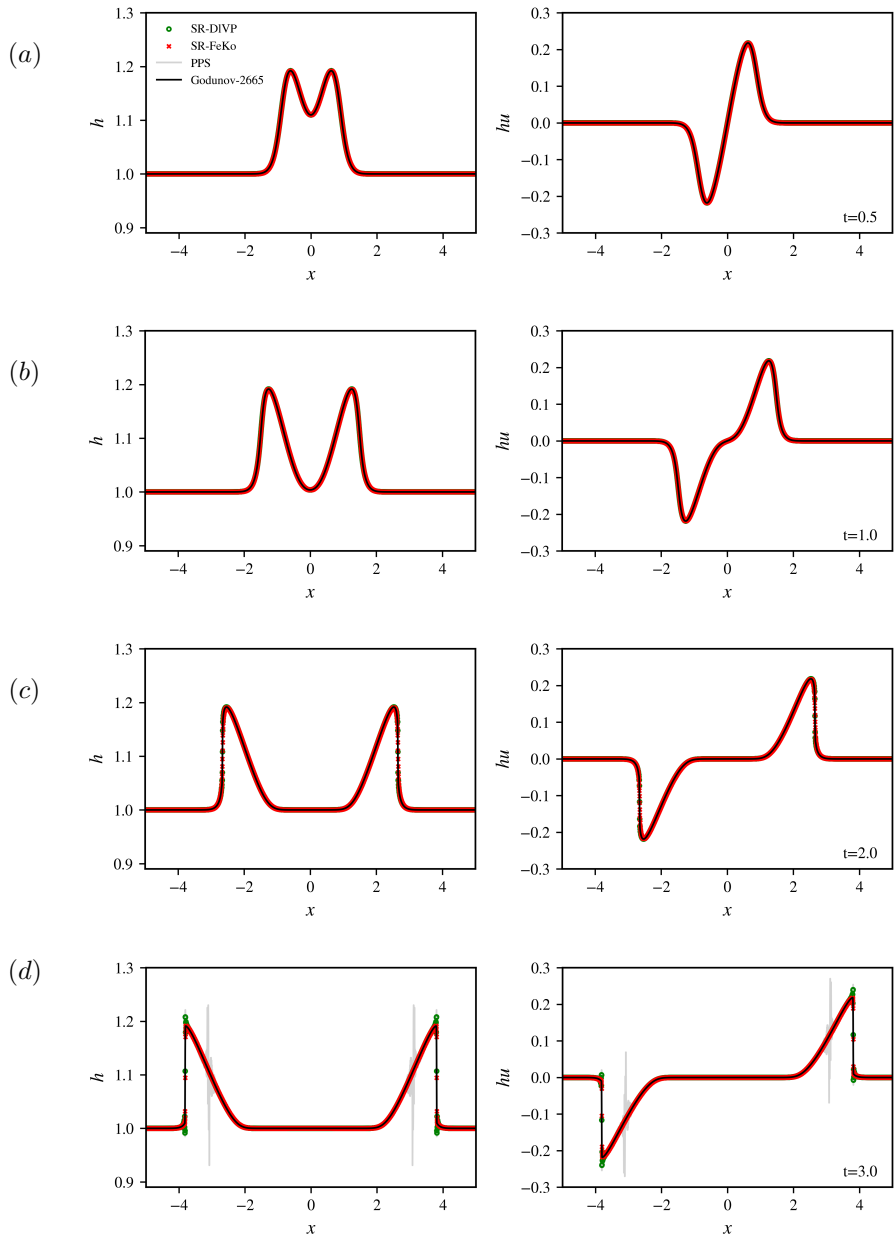
\begin{figure}[h!]
    \centering
    \begin{tikzpicture}
	\draw (0,0) node[inner sep=0]{\includegraphics[width=0.7\linewidth]{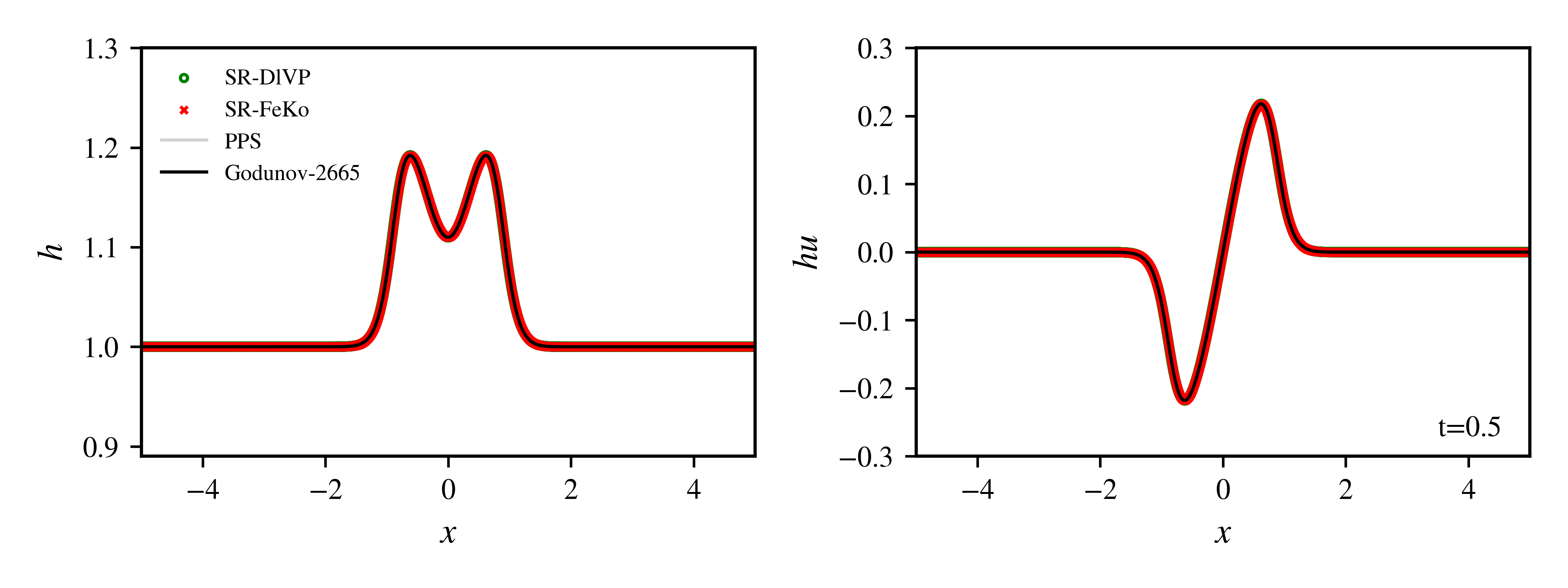}};
	\draw (-6.2,1.2) node {$(a)$};
	\end{tikzpicture}
    \begin{tikzpicture}
 	\draw (0,0) node[inner sep=0]{\includegraphics[width=0.7\linewidth]{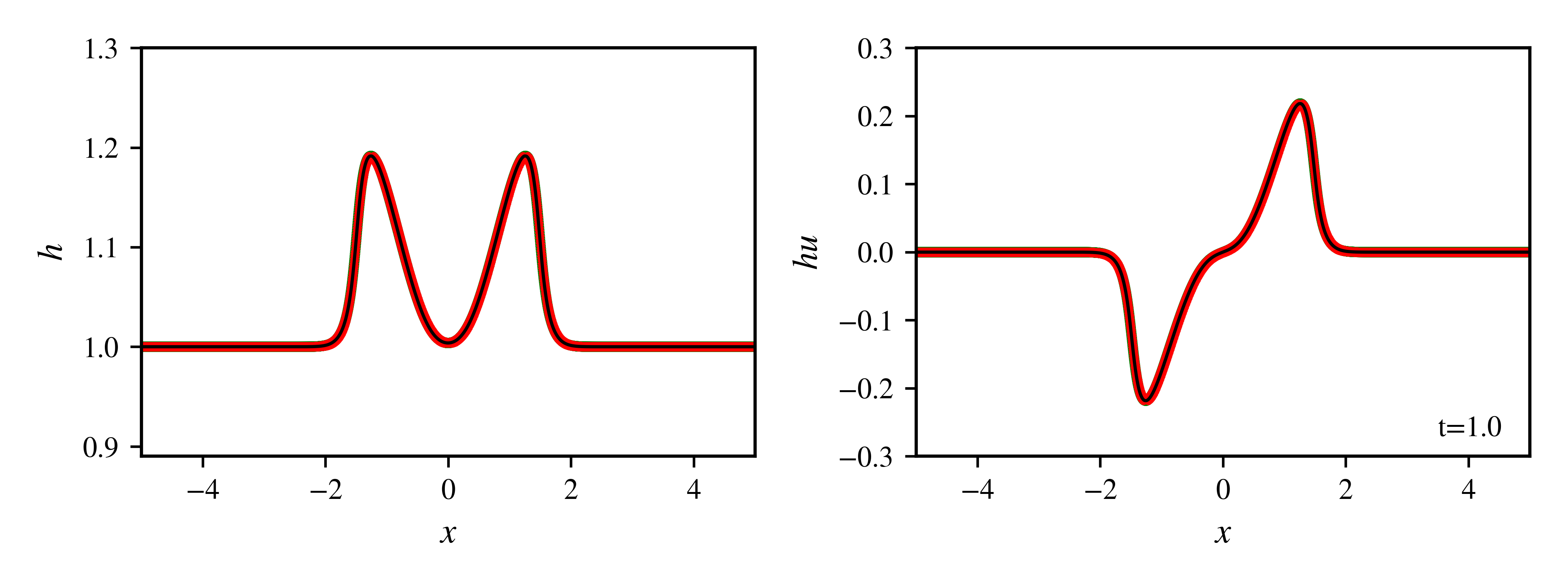}};
	\draw  (-6.2,1.2) node {$(b)$};
	\end{tikzpicture}   
    \begin{tikzpicture}
 	\draw (0,0) node[inner sep=0]{\includegraphics[width=0.7\linewidth]{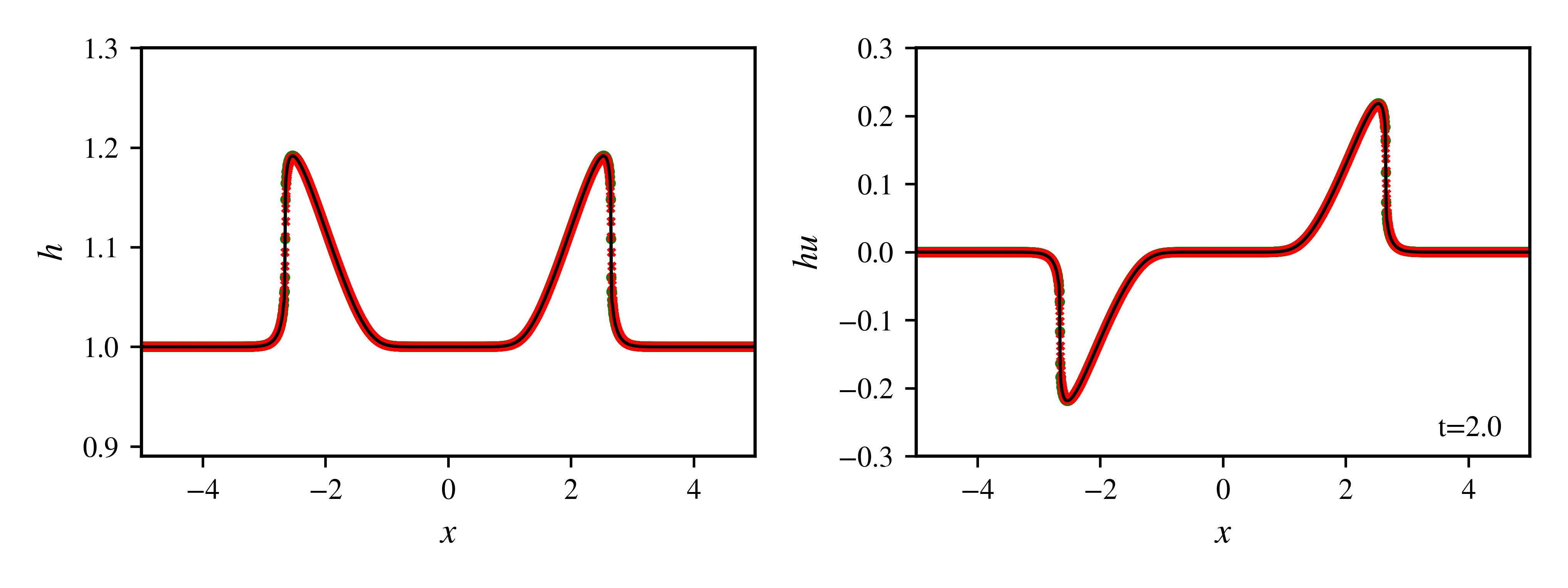}};
	\draw (-6.2,1.2) node {$(c)$};
	\end{tikzpicture}
    \begin{tikzpicture}
 	\draw (0,0) node[inner sep=0]{\includegraphics[width=0.7\linewidth]{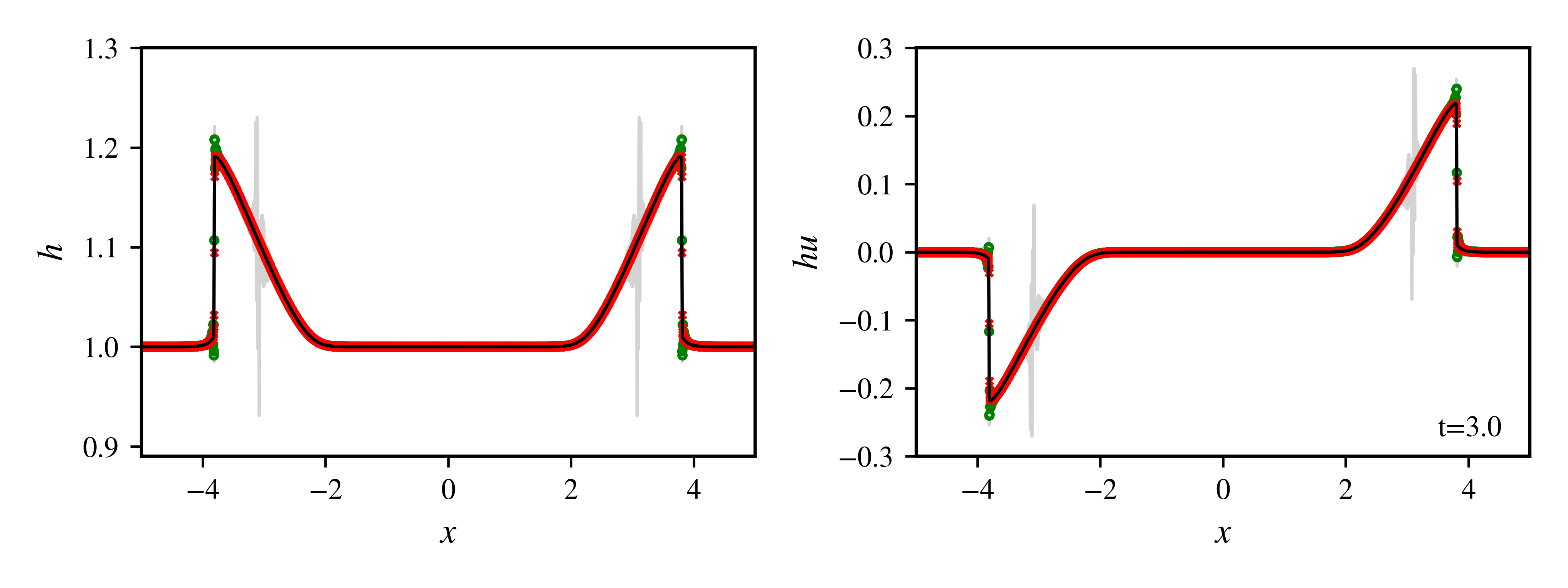}};
	\draw  (-6.2,1.7) node {$(d)$};
	\end{tikzpicture}    
 \caption{ Hump of water: Plot vs. $x$ of $h(x,t)$ (left) and $hu(x,t)$ (right) fields. Different panels give the fields at different instants of time. We show the following approximations in the plot: $2/3$-dealiased PPS (grey line), SR-FeKo with $(\alpha=0.5, \gamma=0.99)$ (red crosses) and SR-DlVP with $(\alpha=0.8225, \gamma=0.98,)$ (green circles). Dealiasing is not performed for the SR runs. The reference solution from \texttt{CLAWPACK} is given by Godunov-2665 (black line) for comparison. All runs have a spatial resolution of $N_x = 2665$. }
    \label{fig:swhump}
\end{figure}

\begin{figure}[h!]
    \centering
    \begin{tikzpicture}
    \draw (0,0) node[inner sep=0]{\includegraphics[width=0.7\linewidth]{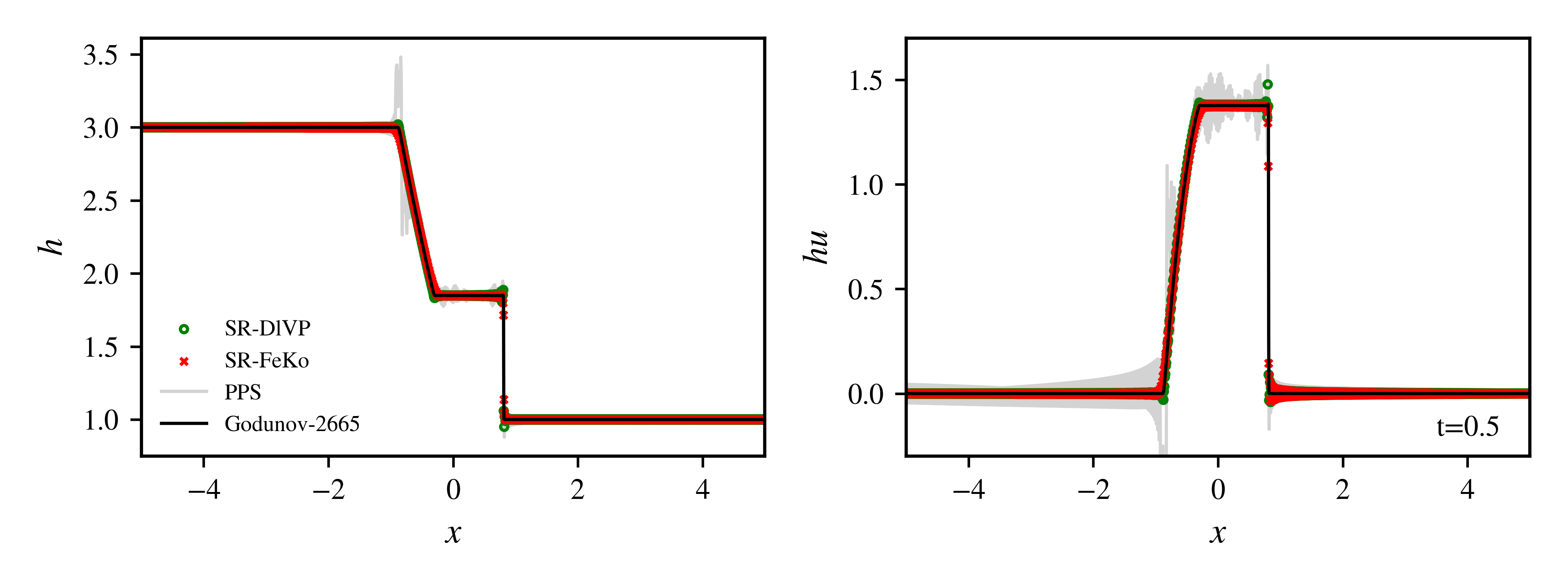}};
	\draw (-6.2,1.2) node {$(a)$};
	\end{tikzpicture}
    \begin{tikzpicture}
 	\draw (0,0) node[inner sep=0]{\includegraphics[width=0.7\linewidth]{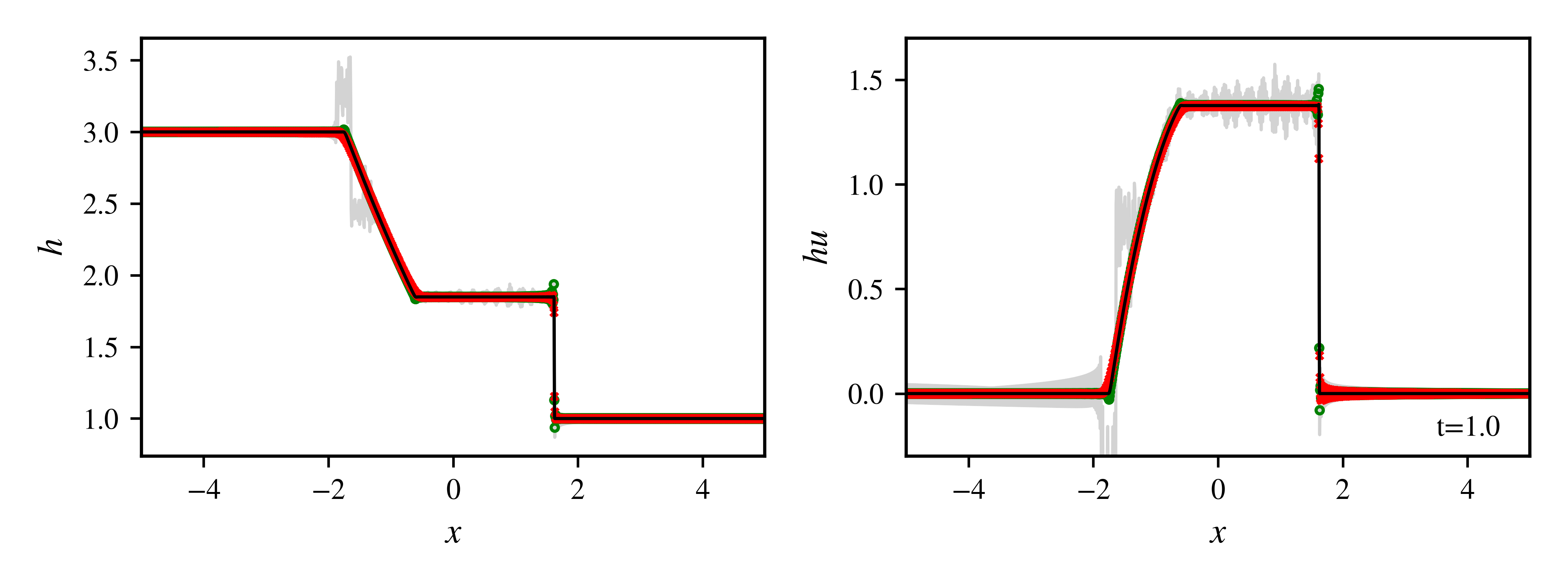}};
	\draw (-6.2,1.2) node {$(b)$};
	\end{tikzpicture}   
    \begin{tikzpicture}
 	\draw (0,0) node[inner sep=0]{\includegraphics[width=0.7\linewidth]{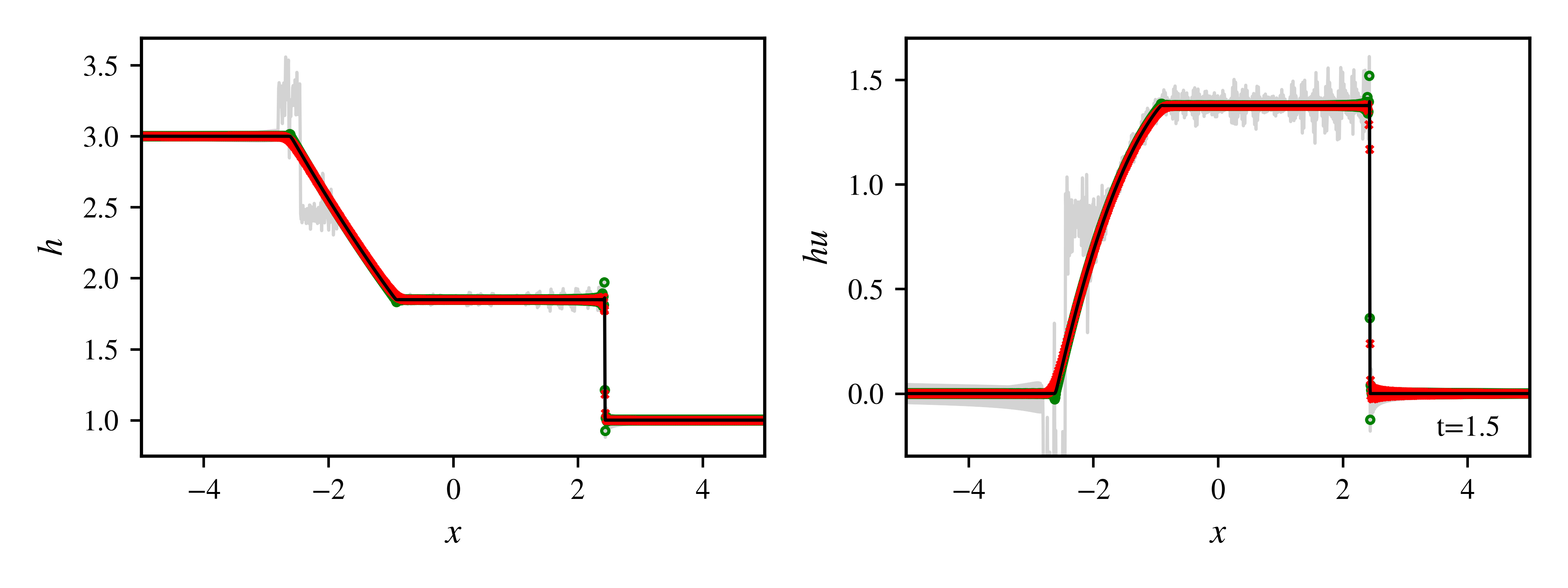}};
	\draw (-6.2,1.2) node {$(c)$};
	\end{tikzpicture}
    \begin{tikzpicture}
	\draw (0,0) node[inner sep=0]{\includegraphics[width=0.7\linewidth]{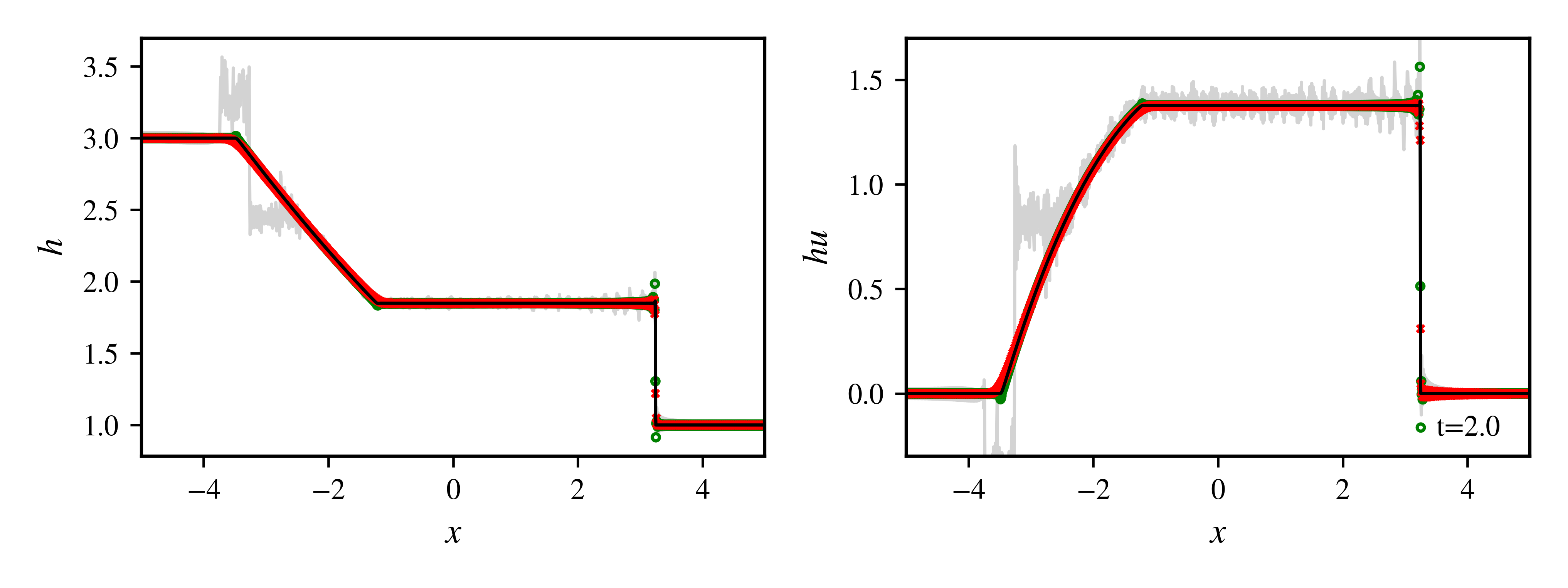}};
	\draw (-6.2,1.2) node {$(d)$};
	\end{tikzpicture}    
  \caption{ Dam break: Plot vs. $x$ of $h(x,t)$ (left) and $hu(x,t)$ (right) fields. Different panels give the fields at different instants of time. We show the following approximations in the plot: $2/3$-dealiased PPS (grey line), SR-FeKo with $(\alpha=0.55, \gamma=0.99)$ (red crosses), and SR-DlVP with $(\alpha=0.7, \gamma=0.98,)$ (green circles). Dealiasing is not performed for the SR runs. The reference solution from \texttt{CLAWPACK} is given by Godunov-2665 (black line) for comparison. All runs have a spatial resolution of $N_x = 2665$. }
    \label{fig:swdambreak}
\end{figure}
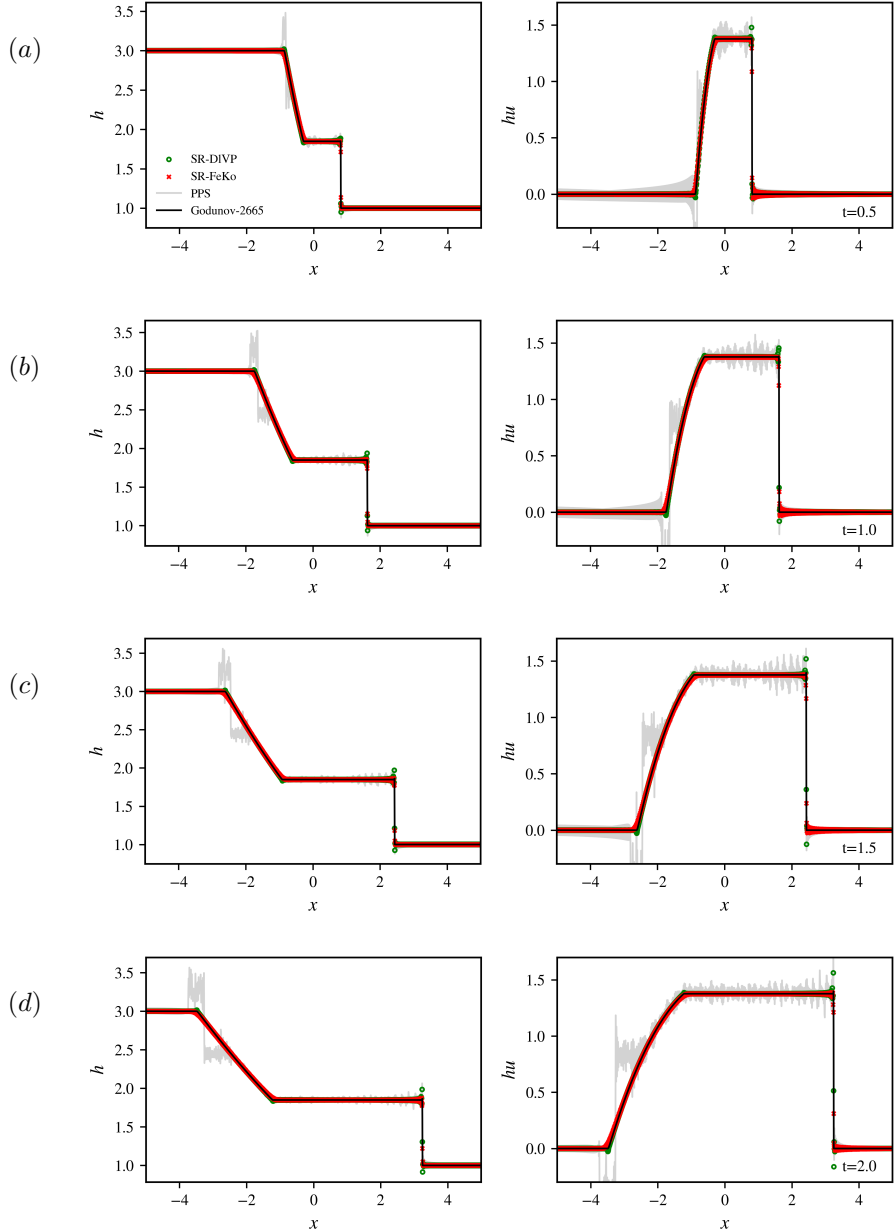

We now extend the spectral relaxation scheme to systems of nonlinear hyperbolic conservation laws; this gives us an opportunity to check the efficiency of the SR method beyond the simplest case of the 1D inviscid Burgers equation. In this section, we consider the 1D shallow water equations, a $2 \times 2$ system of scalar conservation laws. This model governs the flow of a thin layer of fluid of constant density, which is bounded on one side by a pressure surface and on the other side by a fixed surface. In conservative form, the 1D shallow water equations are given by
\begin{subequations}
    \begin{align}
\left[\begin{array}{c}
h \\
h u
\end{array}\right]_t+\left[\begin{array}{c}
u h \\
h u^2+\tfrac{1}{2} g h^2
\end{array}\right]_x=0,
\label{eq:shalwat}
    \end{align}
where $h(x,t)$ is height of the fluid surface, $u(x,t)$ is horizontal velocity of the fluid, and $g$ is the constant of acceleration due to gravity. Thus, the height $h(x,t)$ and the momentum $hu$ are conserved for this system.  

We employ spectral relaxation schemes, with the low-order Fej\'er--Korovkin (SR-FeKo) and the high-order de La Vall\'ee Poussin (SR-DlVP) kernels, to approximate the discontinuous solutions of Eq.~\eqref{eq:shalwat}. We test the schemes on two sets of initial conditions (\texttt{ICs}): $(a)$ the hump of water \texttt{IC} is initially smooth and $(b)$ the dam break \texttt{IC} is initialised as a piecewise constant function. Note that, for a given \texttt{IC}, relaxation is applied to both $h$ and $hu$ fields with the same kernel $K_m(x)$ and parameters $(\alpha,\gamma)$. Reference solutions are the high-resolution ($N_x=2665$ cells) approximations given by Godunov finite-volume schemes of \texttt{CLAWPACK}. For SR approximations that follow, the best parameter sets $(\alpha,\gamma)$ are chosen manually to give the maximum visual convergence to the reference solutions. We report numerical evidence that the SR schemes closely model the reference solution, well beyond the time at which shocks form, for both \texttt{ICs} discussed below. Furthermore, they are capable of following the shock dynamics for the Riemann problem in Eq.~\eqref{eq:shalwat_bw}. 

\begin{itemize}
        \item \textbf{Hump of water:} The water-hump \texttt{IC} is given below for periodic boundary conditions on $x \in [-5,5]$:
                  \begin{align}
          h(x,0) &= 1+0.4 e^{-\beta x^2}, \qquad \qquad      u(x,0) = 0.              \label{eq:shalwat_hump}
          \end{align}
        The hump of water, initialised with zero velocity, gives rise to two waves that move in opposite directions and steepen due to the effects of nonlinearity. The steepening leads to the formation of a shock at time $t \simeq 3.0$. In Fig.~\ref{fig:swhump}, we plot the conserved variables $h(x,t)$ and $hu(x,t)$ at $t=0.5, 1.0, 2.0, 3.0$ in panels $(a)-(d)$, respectively. The dealiased PPS approximation on $N_x=2665$ points (grey line) and the reference solution from \texttt{CLAWPACK} (black line) are shown for comparison. For $t<3$ in panels $(a)-(c)$, the PPS approximation converges to the exact solution which is smooth. After the shock formation in panel $(d)$, PPS approximations of both $h$ and $hu$ fields show tyger-like structures~\cite{ray2011resonance} that develop spontaneously in smooth parts of the solutions, away from the shocks. SR-FeKo approximations with $(\alpha=0.5,\gamma=0.99)$ (red crosses) and SR-DlVP with $(\alpha=0.8225,\gamma=0.98)$ (green circles) are plotted in the same figure. We observe that SR-FeKo and SR-DlVP are highly accurate and overlap with the reference solution for early times $t<3$. In panel $(d)$ of Fig.~\ref{fig:swhump}, SR-FeKo closely follows the reference solution and thus provides a more accurate approximation of the shock. The SR-DlVP approximation is also very good; however, it leads to some bounded oscillations at the shock.  

        \item \textbf{Dam break:} This \texttt{IC} models the breaking of a dam that separates two levels of water at $t=0.$ The water is initially at rest on both sides of the dam. We note that the \texttt{IC} used here is discontinuous and presents a greater challenge for pure pseudospectral methods. Furthermore, the solution is not periodic; this necessitates that we use mirror symmetrisation of the domain $x \in [-5,5]$ about $x=-5$ to continue applying Fourier schemes. The initial configuration is given by
        \begin{equation}
         h(x, 0)=\left\{\begin{array}{ll}
h_l & \text { if } x<0, \\
h_r & \text { if } x>0,
\end{array} \quad u(x, 0)=0,\right.   \label{eq:shalwat_bw}
        \end{equation}
        where $h_l > h_r \ge 0$.

        In Fig.~\ref{fig:swdambreak}, we plot the conserved variables $h(x,t)$ and $hu(x,t)$ at $t=0.5, 1.0, 1.5, 2.0$ in panels $(a)-(d)$ respectively. The dealiased PPS approximation on $N_x=2665$ points (grey line) and the reference solution from \texttt{CLAWPACK} (black line) are shown for comparison. At early times $t=0.5$, shown in panel $(a)$, the PPS approximations of both $h(x,t)$ and $hu(x,t)$ fields develop oscillations near the discontinuity, which leach into the smooth regions of the solutions at later times. We show results from the SR-FeKo approximation using $(\alpha=0.55,\gamma=0.99)$ (red crosses) and the SR-DlVP approximation using $(\alpha=0.70,\gamma=0.98)$ (green circles) in Fig.~\ref{fig:swdambreak}. We notice that the SR-DlVP approximation displays bounded oscillations near the shock, but follows the reference solution more closely. As it approaches the shocks, the SR-FeKo approximation is slightly overdamped and displays bounded oscillations on the right corners. However, far from the discontinuity, it provides a considerably accurate approximation of the flow fields. Both SR approximations are able to follow the movement of discontinuities in this case.   
        
\end{itemize}        
\end{subequations}
   
\subsection{3x3 system: 1D compressible Euler equations}
\label{subsec:Results_Euler}

The 1D compressible Euler equations of gas dynamics constitute a $3 \times 3$ system of nonlinear hyperbolic conservation laws. In conservative form, it is given by 
\begin{subequations}
    \begin{align}
\left[\begin{array}{c}
\rho \\
E \\
\rho u
\end{array}\right]_t+\left[\begin{array}{c}
\rho u \\
(E+p) u \\
\rho u^2+p 
\end{array}\right]_x=0 \,,
\label{eq:1deulereq}
    \end{align}
where $\rho(x,t)$ is the density, $u(x,t)$ is the velocity, $E(x,t)$ is the total energy, and $p=(\gamma-1)\rho$ is the pressure of the gas. We set $\gamma =1.4$ for the ratio of specific heats of diatomic gases (major composition of air).  

We consider the initial-boundary value problem for Eq.~\eqref{eq:1deulereq} on the non-periodic, bounded domain $x \in [x_L,x_R]$. We use Chebyshev polynomials of the first kind $T_k(x)= \cos(k \cos^{-1} (x))$ as the basis for pseudospectral projection. Discrete Chebyshev coefficients are determined by quadratures, analogous to Eq.~\eqref{eq:fourcoeffs}, using the discrete cosine transform (see \texttt{REDFT00} in \texttt{FFTW3} library). To implement boundary conditions at $x_L$ and $x_R$, these points should be present in the collocation grid. We thus use Chebyshev extrema points given by $\{ X_j = \cos(\pi j/N), j = 0, 1,..., N \}$ where $N$ is order of the highest-order polynomial in the pseudospectral basis. Since the Chebyshev polynomials are defined on $X \in [-1,1]$, we use a simple rescaling $x=(X+1)/2$, for the blast-waves \texttt{ICs} in Eq.~\eqref{eq:1deul_bw} where the domain is $x \in [0,1]$.  

Since Chebyshev extrema points cluster and form a fine grid at the boundaries, we are constrained to evolve the system with a very small time step. For $N$-point Chebyshev discretisation, the CFL stability criterion gives $\Delta t \sim \mathcal{O}(1/N^2)$. 
We can alleviate this restriction by once more rescaling the grid to redistribute the higher density of points at the walls. 
We follow the procedure of~\cite{kosloff1993modified} and apply the transformation given below to Chebyshev extrema points $\{ x_j , j=0, 1,...,  N \}$:
\begin{align}
    \chi_j = \frac{\sin^{-1} (\beta \ x_j)}{\sin^{-1} (\beta)}\,,
    \label{eq:kosloff}
\end{align}

where we use $\beta=0.999$. Now, $\{ \chi_j, j=0, 1,..., N \}$ is the computational grid that we finally employ. The detailed implementation of boundary conditions, which we use in the test problems, is presented in Appendix~\ref{app:cbc}.  

The relaxation kernel is now applied to the Chebyshev spectral approximations; note that $N_x=N$ and the relations $m(N)$ and $\tau(N)$ remain as in Eq.~\eqref{eq:sr_params}. We test the SR-FeKo and SR-DlVP schemes on standard test problems for shock-fitting schemes and discuss their performances below. The parameters $(\alpha,\gamma)$ are chosen manually to give the best visual convergence for these models.  
\begin{itemize}
    \item \textbf{Sod Shock tube} This Riemann-type problem is solved on the domain $ x \in [-1,1]$, where the boundaries at $x_L=-1$ and $x_R=1$ are treated as \textit{reflecting solid walls}~\cite{sod1978survey}. The initial data are given by
        \begin{equation}
            \begin{array}{lr}
\left(\rho_l, u_l, p_l\right)=(1,0,1) & -1 \leqslant x \leqslant 0 ,\label{eq:1deul_ss} \\
\left(\rho_r, u_r, p_r\right)=(0.125,0,0.10) & 0 \leqslant x \leqslant 1.
\end{array}
        \end{equation}

    In Fig.~\ref{fig:eulsod}, we plot the fields $\rho$ (left), $u$ (center), and $p$ (right) at consecutive times in panels $(a)-(d)$.
    The reference solution from \texttt{HyPar} (black line) is shown for comparison --- it is computed using a Characteristic-based 5th order WENO method for a spatial resolution of $N_x=8000$ points. The dealiased Chebyshev-PPS approximation on $N_x=615$ points destabilises quickly and cannot be seen in the plots; for consistency, it is shown in the legend. Despite the challenge posed by this test case for PPS schemes, the SR schemes are able to approximate accurately the reference solution. In particular, SR-FeKo with $(\alpha=0.785,\gamma=0.99)$ produces an oscillation-free approximation of the solution; it is however slightly overdamped near shocks and corners. The SR-DlVP approximation with $(\alpha=0.94,\gamma=0.95)$ models the piecewise profile more closely, at the cost of bounded and small oscillations near corners. We do not implement $2/3$-dealiasing for these SR runs. The spatial resolution used for the spectral runs is $N_x=615$.
    
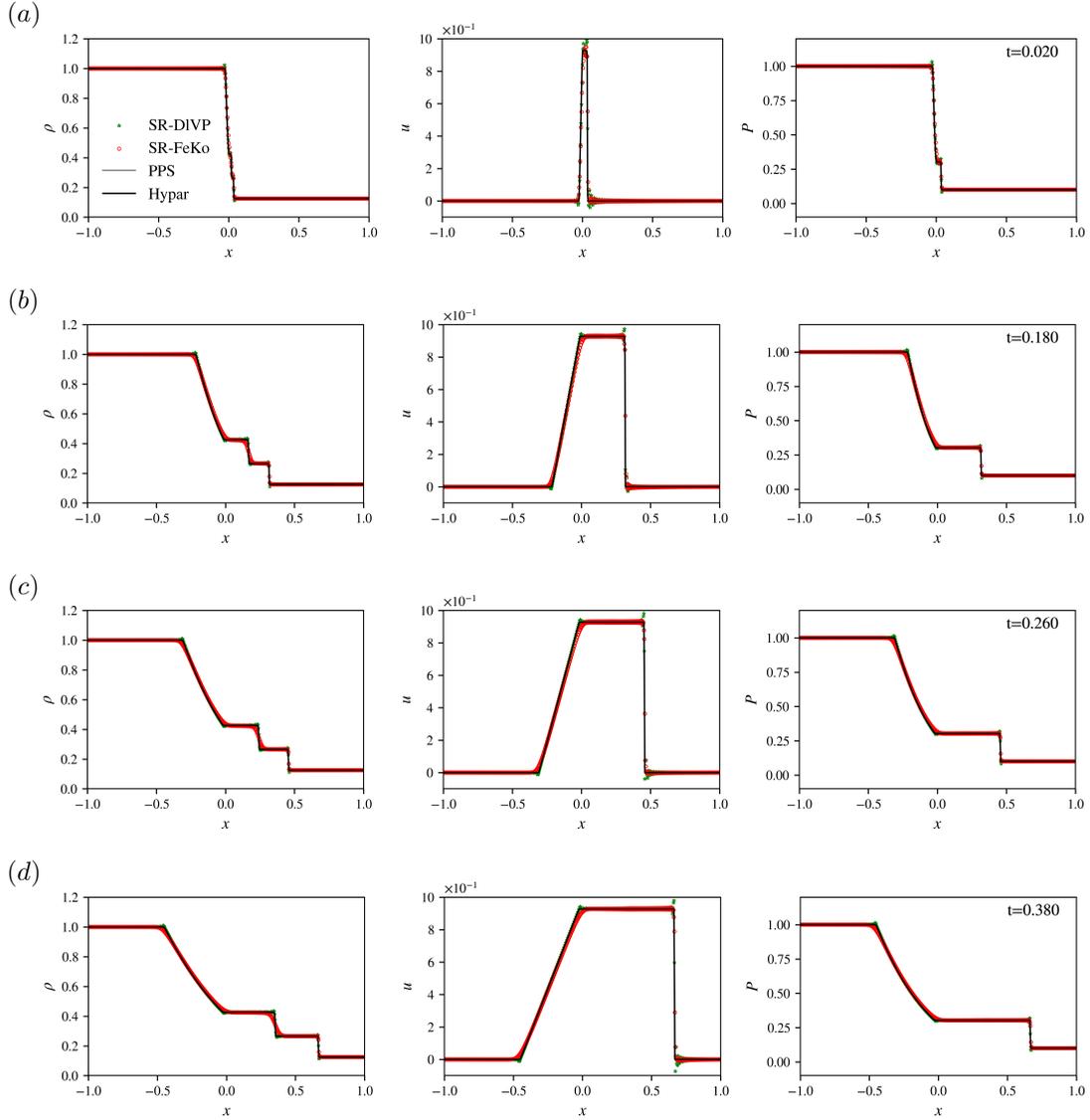
\begin{figure}[h!]
    \centering
        \begin{tikzpicture}
    	\draw (0,0) node[inner sep=0]{\includegraphics[width=0.9\linewidth]{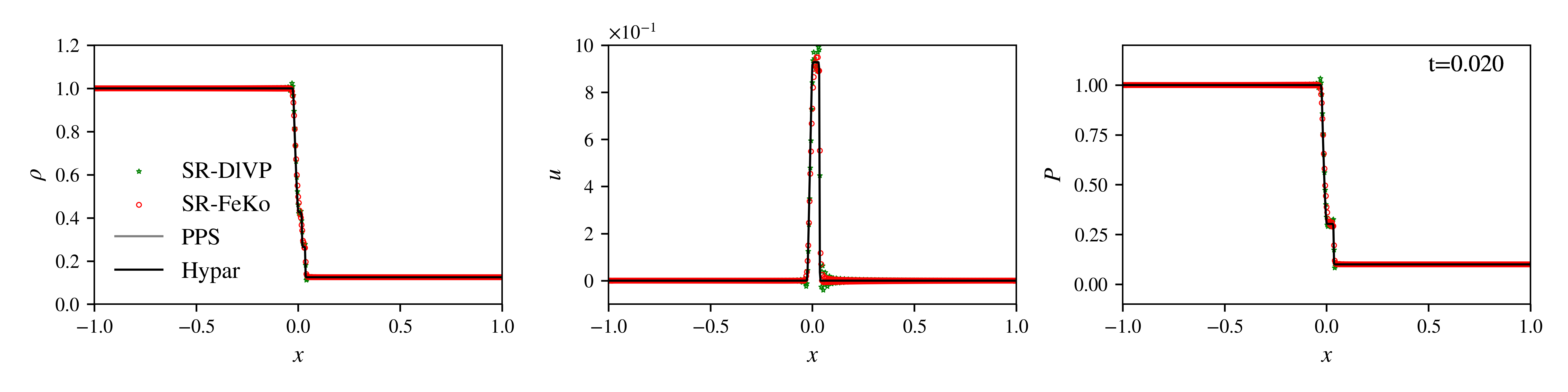}};
        \draw (-7.2,1.7) node {$(a)$};
    	\end{tikzpicture}
        \begin{tikzpicture}
     	\draw (0,0) node[inner sep=0]{\includegraphics[width=0.9\linewidth]{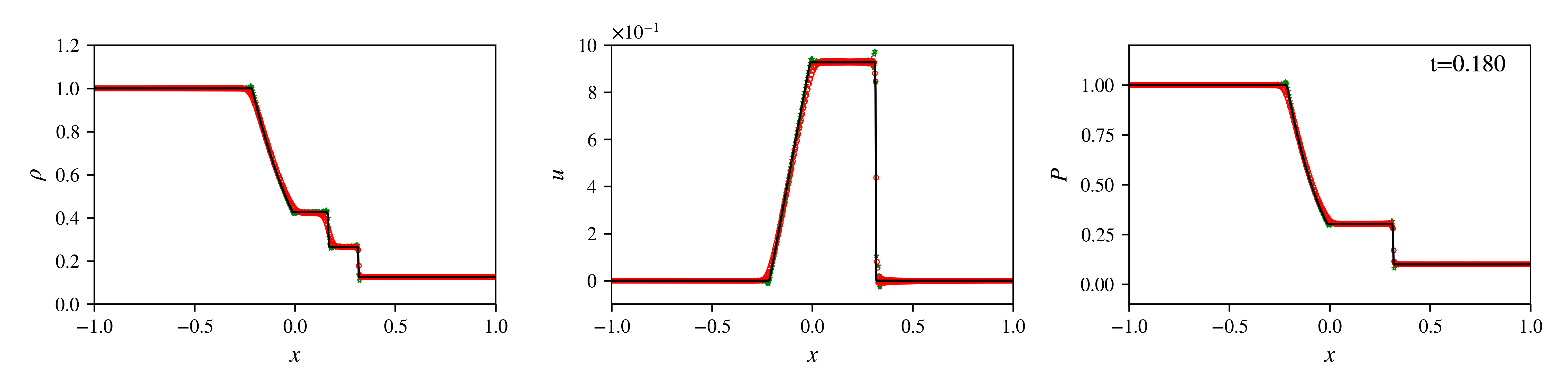}};
        \draw  (-7.2,1.7) node {$(b)$};
    	\end{tikzpicture}   
        \begin{tikzpicture}
     	\draw (0,0) node[inner sep=0]{\includegraphics[width=0.9\linewidth]{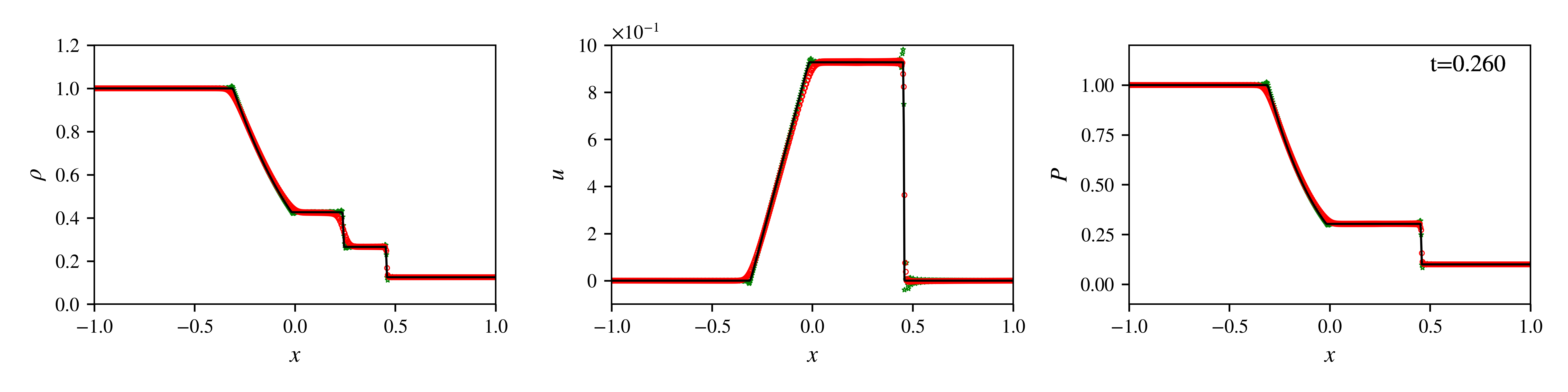}};
        \draw (-7.2,1.7) node {$(c)$};
    	\end{tikzpicture}
        \begin{tikzpicture}
    	\draw (0,0) node[inner sep=0]{\includegraphics[width=0.9\linewidth]{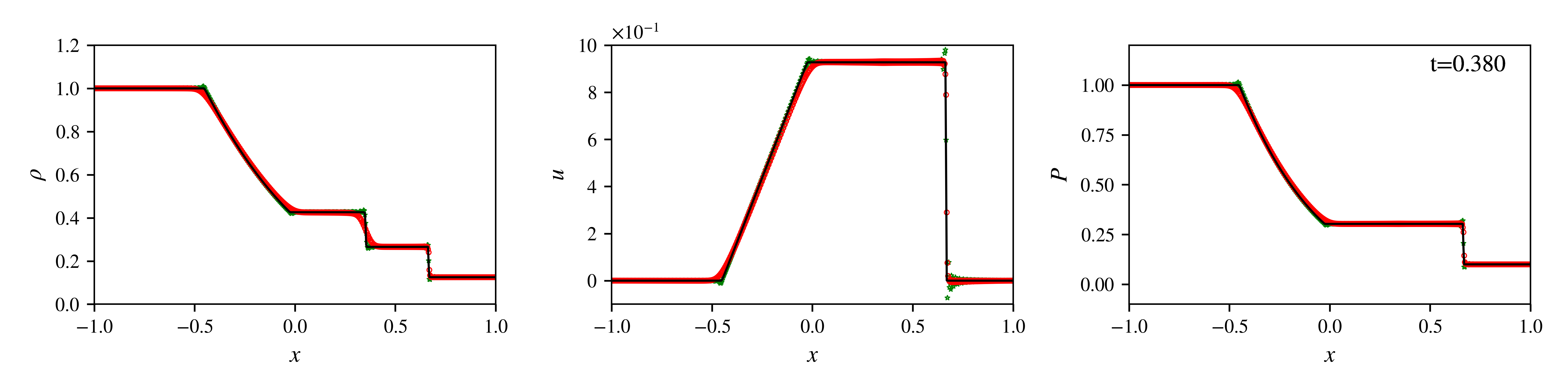}};
        \draw  (-7.2,1.7) node {$(d)$};
    	\end{tikzpicture}    
    \caption{ Sod Shock tube: Plots vs. $x$ of $\rho,u$, and $p$ fields. Different panels (rows) give the fields at different instants of time. The reference solution from \texttt{HyPar} is shown for comparison (black line). PPS approximation (grey solid line) destabilises early for this \texttt{IC} and is largely absent. SR-FeKo with $(\alpha=0.785, \gamma=0.99)$ (red circles) and SR-DlVP with $(\alpha=0.94, \gamma=0.95)$ (green stars) approximate the reference solution closely. All the runs presented here have a resolution of $N_x = 615$. Dealiasing is not performed for SR runs. }
    \label{fig:eulsod}
\end{figure}

  \item \textbf{Lax Shock tube} The domain is given by $ x \in [-1,1] $ where the boundaries at $x_L=-1$ and $x_R=1$ are treated as \textit{reflecting solid walls}~\cite{lax1954weak}. The initial data are given by
            \begin{equation}
            \begin{array}{lr}
\left(\rho_l, u_l, p_l\right)=(0.445,0.311,8.928) & -1 \leqslant x \leqslant 0, \label{eq:1deul_lx}\\
\left(\rho_r, u_r, p_r\right)=(0.5,0,1.4275) & 0 \leqslant x \leqslant 1.
\end{array}
        \end{equation}

For this test problem, we plot the fields $\rho$ (left), $u$ (center), and $p$ (right) at consecutive times in panels $(a)-(d)$ of Fig.~\ref{fig:eullax}. The reference solution is obtained from \texttt{HyPar} with $N_x=8000$ as detailed for the Sod shock-tube problem above. In panel $(a)$ at time $t=0.02$, we see the dealiased Chebyshev-PPS approximation on $N_x=615$ points. It destabilises quickly and cannot be seen in panels $(b)-(d)$; for consistency, we show it in the legend. The SR-FeKo scheme with $(\alpha=0.91,\gamma=0.99)$ produces an oscillation-free approximation of the solution; it is, however, moderately overdamped near discontinuities and corners. The SR-DlVP approximation with $(\alpha=1.10,\gamma=0.98)$ models the piecewise profile more closely, at the cost of small, bounded oscillations near corners. We do not implement $2/3$-dealiasing for these runs. The spatial resolution used for the spectral runs is $N_x=615$. 

\begin{figure}[h!]
    \centering
        \begin{tikzpicture}
        \draw (0,0) node[inner sep=0]{\includegraphics[width=0.9\linewidth]{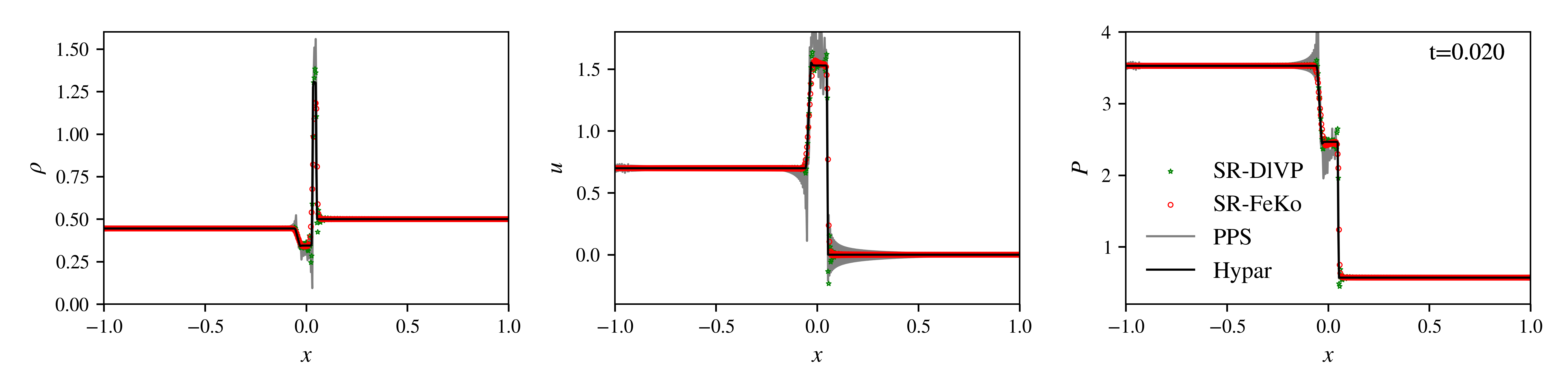}};
        \draw (-7.2,1.7) node {$(a)$};
        \end{tikzpicture}
        \begin{tikzpicture}
        \draw (0,0) node[inner sep=0]{\includegraphics[width=0.9\linewidth]{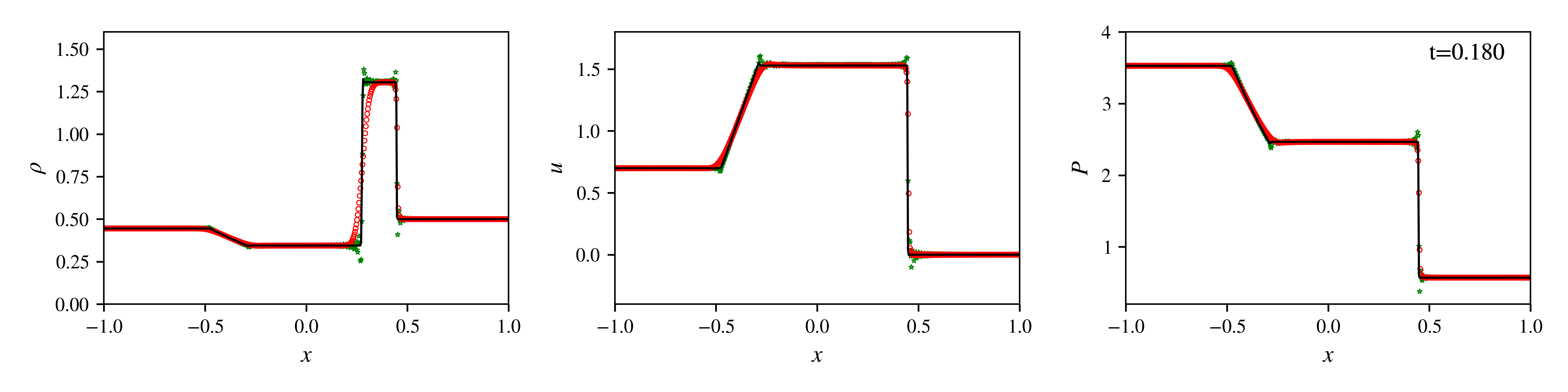}};
        \draw  (-7.2,1.7) node {$(b)$};
        \end{tikzpicture}   
        \begin{tikzpicture}
        \draw (0,0) node[inner sep=0]{\includegraphics[width=0.9\linewidth]{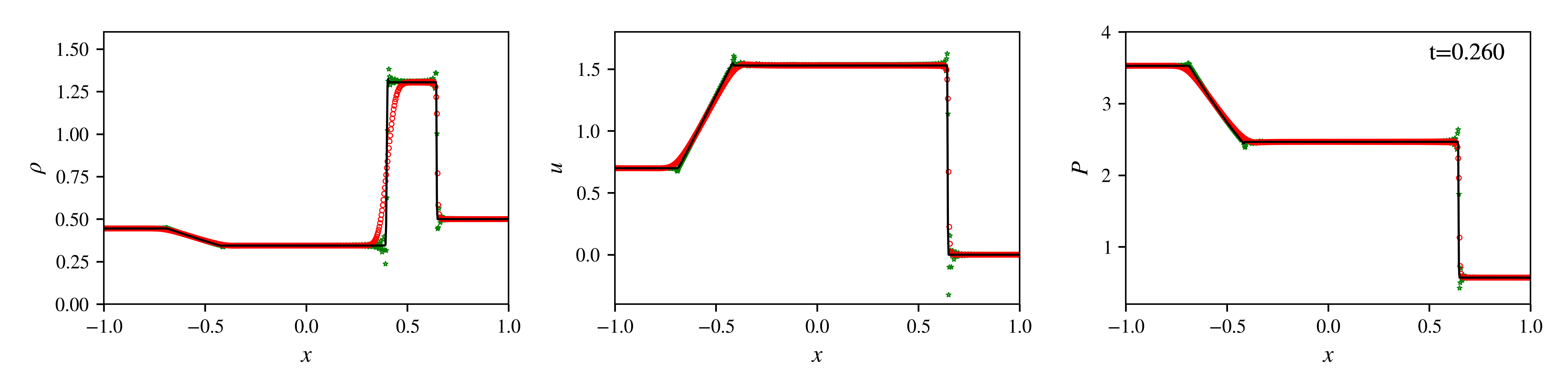}};
        \draw (-7.2,1.7) node {$(c)$};
        \end{tikzpicture}
        \begin{tikzpicture}
        \draw (0,0) node[inner sep=0]{\includegraphics[width=0.9\linewidth]{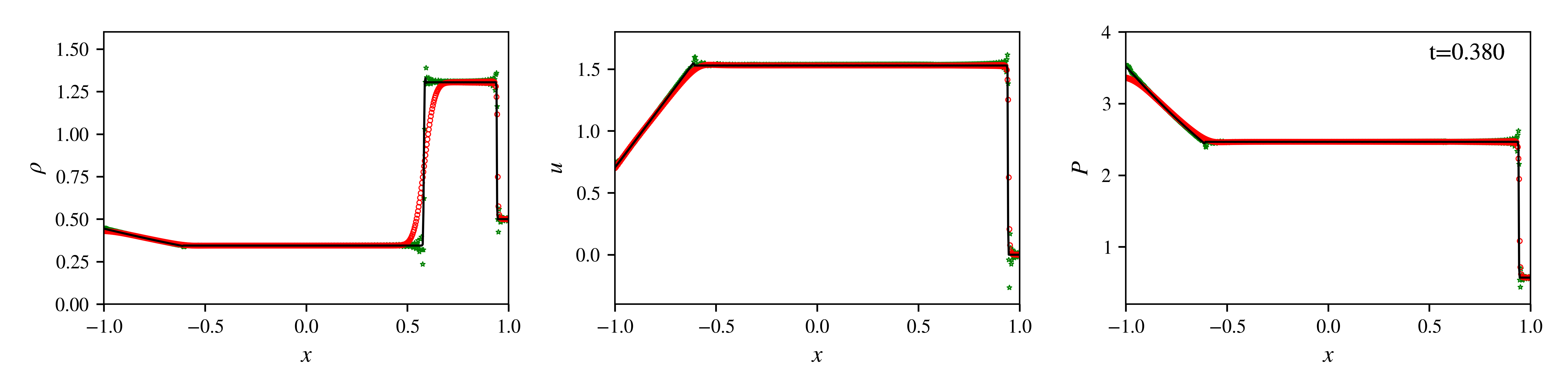}};
        \draw  (-7.2,1.7) node {$(d)$};
        \end{tikzpicture}    
    \caption{ Lax Shock tube: Plots vs. $x$ of $\rho,u$, and $p$ fields. Different panels (rows) give the fields at different instants of time. The reference solution from \texttt{HyPar} is shown for comparison (black line). PPS approximation (grey solid line) destabilises early for this \texttt{IC} and is largely absent. SR-FeKo with $(\alpha=0.91, \gamma=0.99)$ (red circles) is moderately overdamped but SR-DlVP with $(\alpha=1.10, \gamma=0.98)$ (green stars) approximate the reference solution closely, despite small bounded oscillations near corners. All the runs presented here have a resolution of $N_x = 615$. Dealiasing is not performed for SR runs.}
    \label{fig:eullax}
    \end{figure}
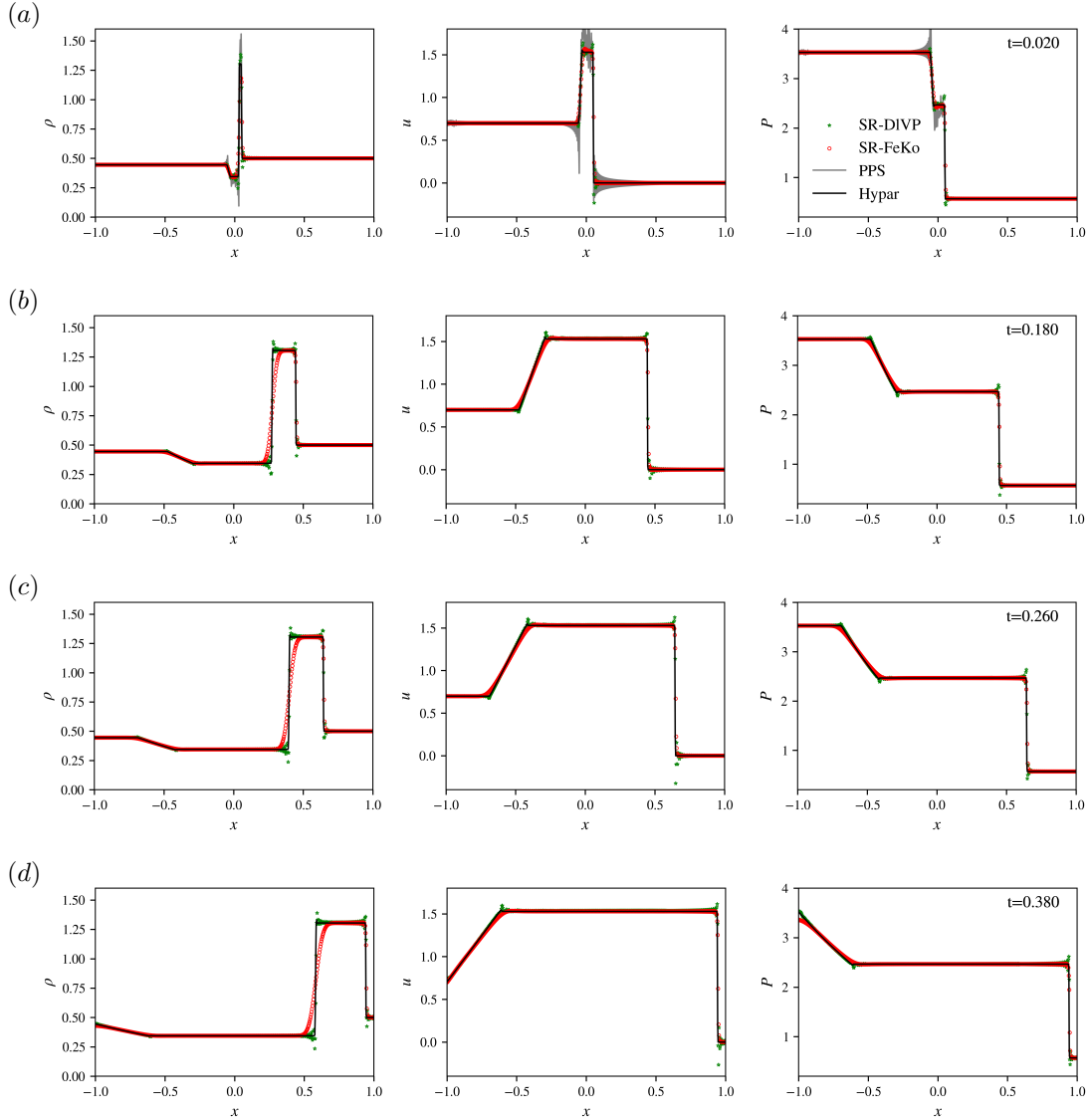
   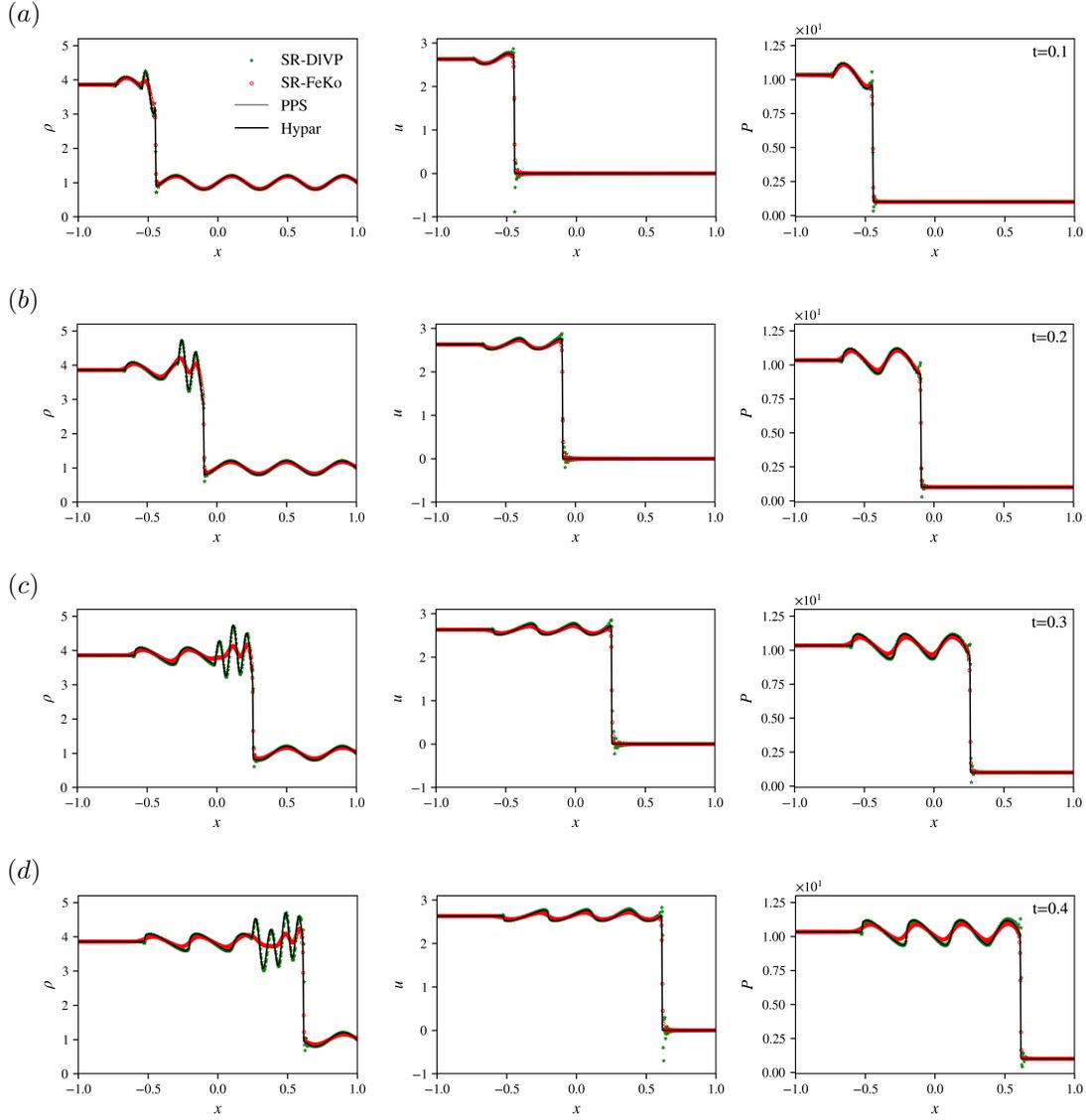
\begin{figure}[h!]
    \centering
    \begin{tikzpicture}
	\draw (0,0) node[inner sep=0]{\includegraphics[width=0.9\linewidth]{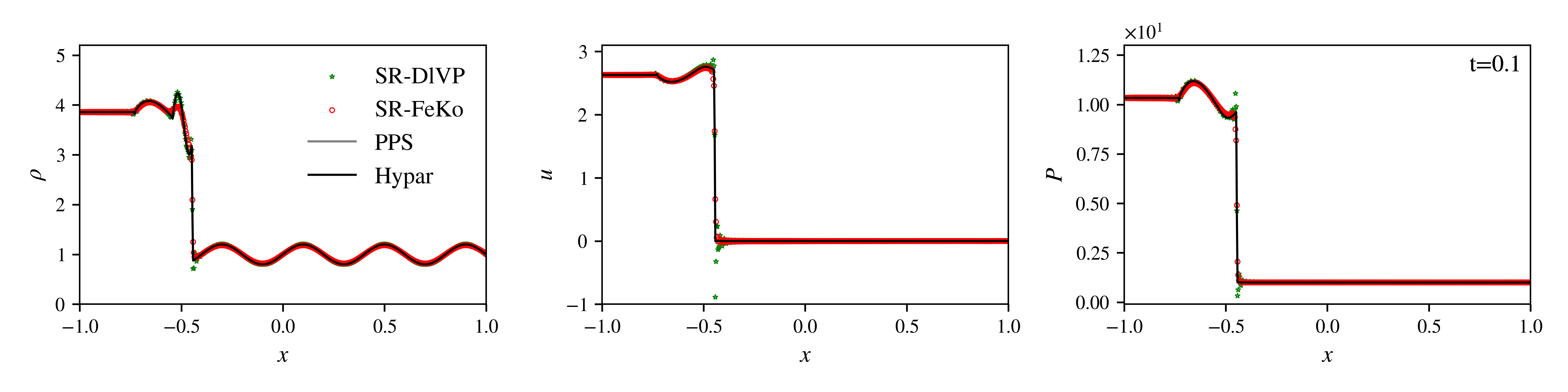}};
		\draw (-7.2,1.7) node {$(a)$};
	\end{tikzpicture}
    \begin{tikzpicture}
 	\draw (0,0) node[inner sep=0]{\includegraphics[width=0.9\linewidth]{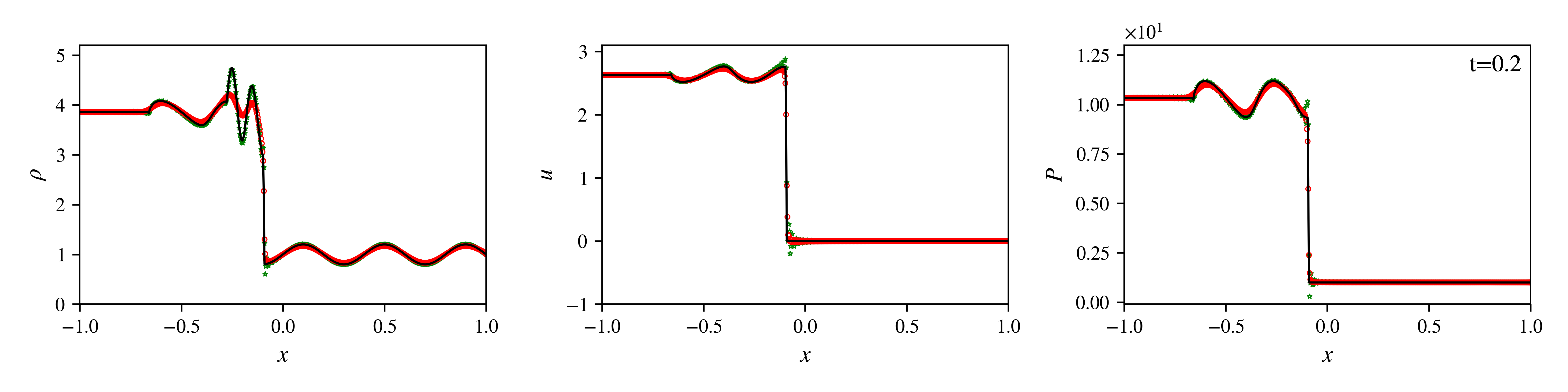}};
		\draw  (-7.2,1.7) node {$(b)$};
	\end{tikzpicture}   
    \begin{tikzpicture}
 	\draw (0,0) node[inner sep=0]{\includegraphics[width=0.9\linewidth]{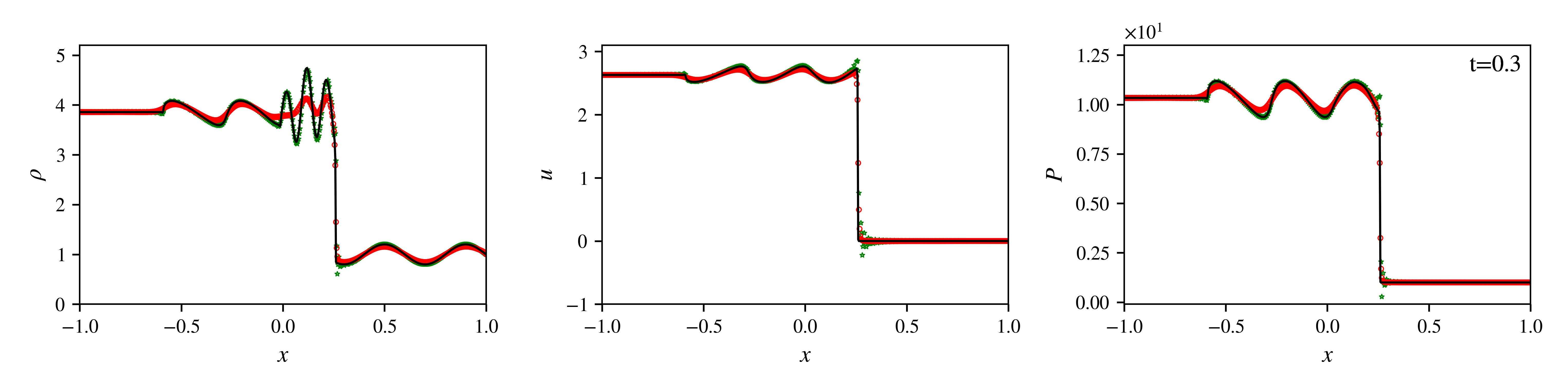}};
		\draw (-7.2,1.7) node {$(c)$};
	\end{tikzpicture}
    \begin{tikzpicture}
	\draw (0,0) node[inner sep=0]{\includegraphics[width=0.9\linewidth]{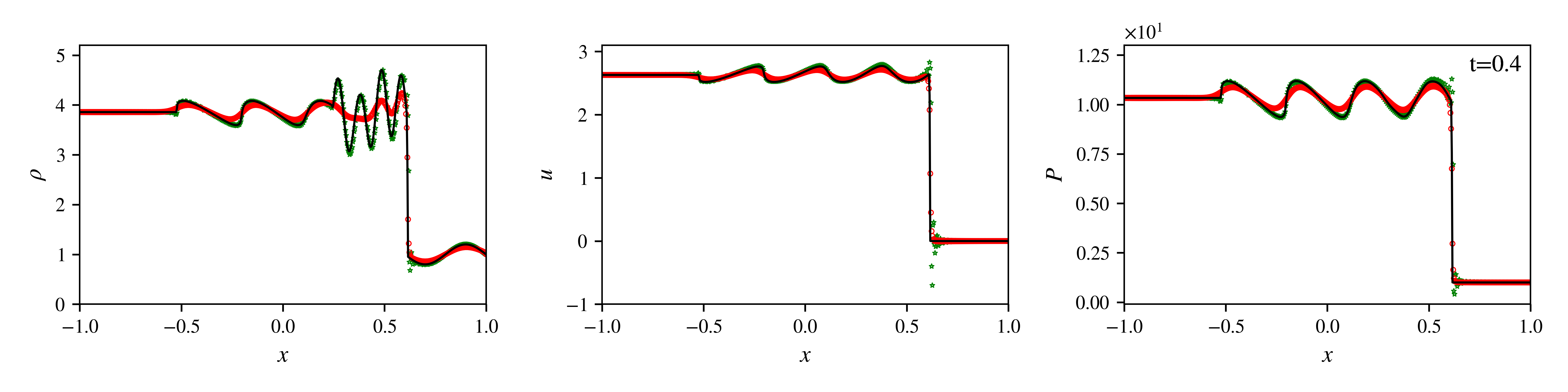}};
		\draw  (-7.2,1.7) node {$(d)$};
	\end{tikzpicture}
     \caption{ Osher--Shu problem: Plots vs. $x$ of $\rho,u$, and $p$ fields. Different panels (rows) give the fields at different instants of time. The reference solution from \texttt{HyPar} is shown for comparison (black line). The PPS approximation (grey solid line) destabilises early for this \texttt{IC} and is absent except in legend. SR-FeKo with $(\alpha=0.95, \gamma=0.97)$ (red circles) does not converge to the correct solution. SR-DlVP with $(\alpha=1.14, \gamma=0.95)$ (green stars) is able to approximate the shock as well as the sinusoidal variation in the reference solution closely, despite small, bounded oscillations near corners. All the runs presented here have a resolution of $N_x = 615$. Dealiasing is not performed for SR runs.}
    
    \label{fig:eulos}
\end{figure}

\begin{figure}[h!]
    \centering
    \begin{tikzpicture}
	\draw (0,0) node[inner sep=0]{\includegraphics[width=0.9\linewidth]{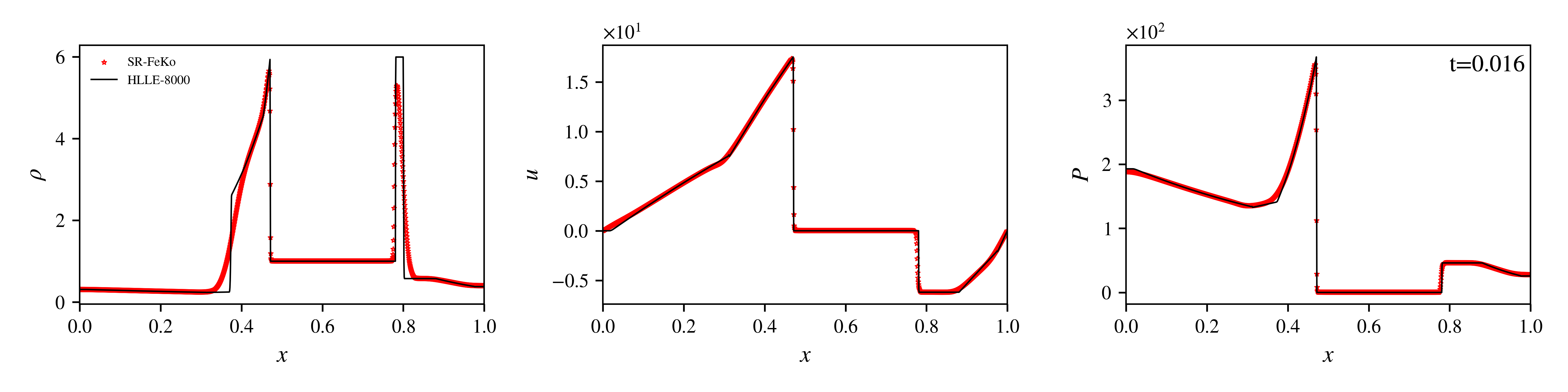}};
		\draw (-7.2,1.7) node {$(a)$};
	\end{tikzpicture}
    \begin{tikzpicture}
 	\draw (0,0) node[inner sep=0]{\includegraphics[width=0.9\linewidth]{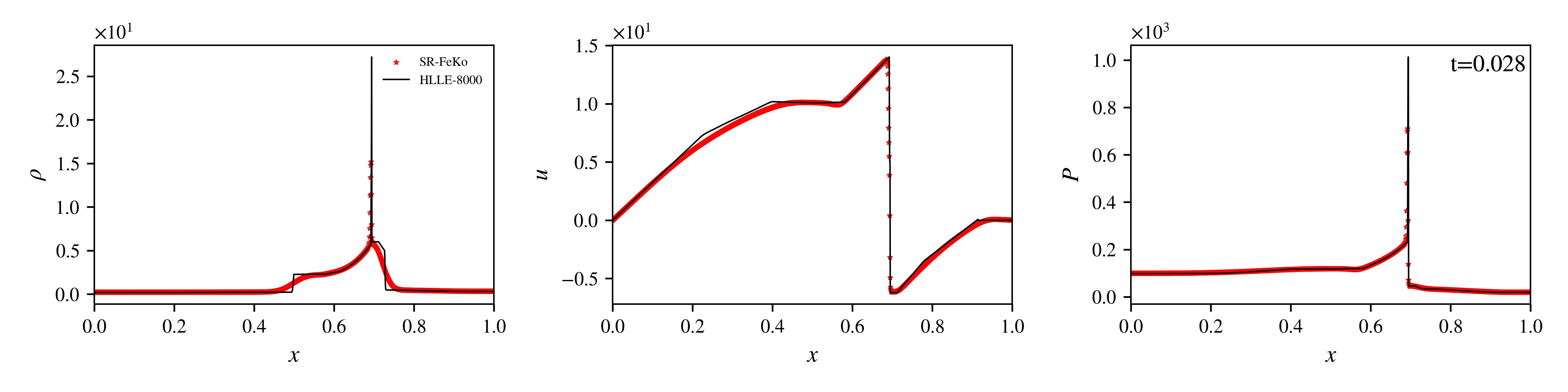}};
		\draw  (-7.2,1.7) node {$(b)$};
	\end{tikzpicture}   
    \begin{tikzpicture}
 	\draw (0,0) node[inner sep=0]{\includegraphics[width=0.9\linewidth]{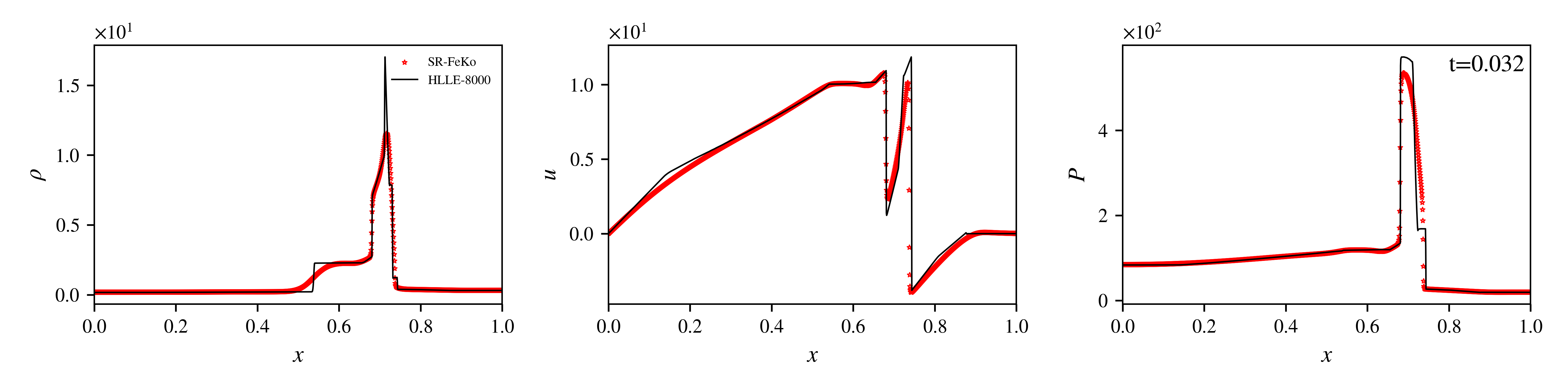}};
		\draw (-7.2,1.7) node {$(c)$};
  	\end{tikzpicture}  
      \begin{tikzpicture}
 	\draw (0,0) node[inner sep=0]{\includegraphics[width=0.9\linewidth]{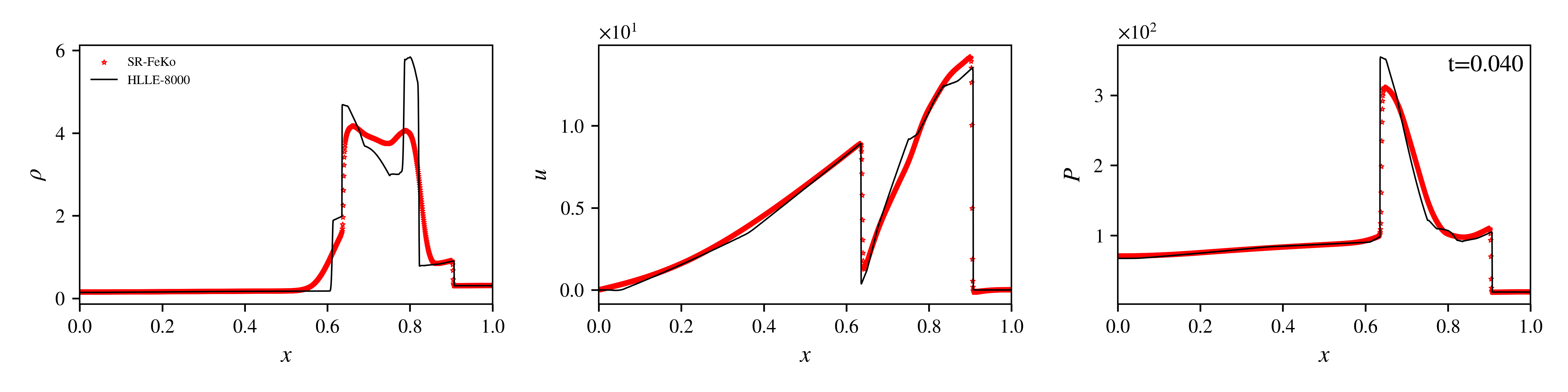}};
		\draw (-7.2,1.7) node {$(d)$};
	\end{tikzpicture}  
 \caption{ Blast-waves interaction:  Plots vs. $x$ of $\rho,u,$ and $p$ fields. Different panels (rows) give the fields at different instants of time. The reference solution from \texttt{CLAWPACK} is shown for comparison (black line).
 SR-FeKo with $(\alpha=1.355, \gamma=0.999)$ (red crosses) provides a low-order approximation of the solution. It fails to resolve the fine structures emanating in the wake of the blast-waves interaction in panel $(d)$. The SR-FeKo run presented here has a spatial resolution of $N_x = 1599$. Dealiasing is not performed for SR runs.}
    \label{fig:eulbw}
\end{figure}

    \item \textbf{Shu--Osher problem} This test case models the interaction between a Mach-3 shock moving towards the right, and an entropy wave with a sinusoidal profile~\cite{shu1988efficient}. The domain is $x \in [-1,1]$ and the boundary conditions are \textit{supersonic inflow} at $x_L=-1$ and \textit{outflow} at $x_R=1$.  The initial configuration is given by
        \begin{align}
            (\rho, u, p)_{t=0}=\left\{\begin{array}{lcrc}(3.85714, & 2.629369,& 10.33333) & x \leq-0.8, \\ (1.0+\epsilon \sin (\kappa \pi x), & 0, & 1.0) & x>-0.8,\end{array}\right. \label{eq:1deul_os}
        \end{align}
        where $\epsilon=0.2$ and $\kappa=5$.
     Since the true solution contains fine structure, we expect that a low-order scheme cannot provide an accurate approximation. This is reflected in the results presented in Fig.~\ref{fig:eulos}. Here, we plot the fields $\rho$ (left), $u$ (center), and $p$ (right) at consecutive times in panels $(a)-(d)$. The reference solution (black line) is obtained from \texttt{HyPar} with $N_x=8000$ as discussed above for the Sod shock-tube problem. This test case is challenging for the dealiased-PPS scheme; it destabilises quickly and cannot be seen in panels $(b)-(d)$. However, for consistency, we show it in the legend. We show the SR-FeKo approximations on $N_x=615$ points using $(\alpha=0.95,\gamma=0.97)$ (red circles). Since the Fej\'er--Korovkin kernel is of low-order, the SR-FeKo scheme fails to capture the oscillatory profile of the solution. By using the high-order de La Vall\'ee Poussin kernel, the SR-DlVP approximation on  the same grid, with $(\alpha=1.14,\gamma=0.95)$ (green stars) is able to resolve the shock as well as the fine oscillating structure in this problem. The Gibbs oscillations at the shock for SR-DlVP approximation can be further controlled by choosing better parameters $(\alpha,\gamma)$ found through a systematic search.

    \item \textbf{Interaction of blast waves} Finally, we model the interaction of two blast waves~\cite{woodward1984numerical}, starting from the initial data,
        \begin{equation}
(\rho,u,p)(x, t=0)= \begin{cases}(\rho_l,u_l,p_l) & 0 \leqslant x \leqslant 0.1, \\ (\rho_m,u_m,p_m) & 0.1 \leqslant x \leqslant 0.9, \\ (\rho_r,u_r,p_r) & 0.9 \leqslant x \leqslant 1,\end{cases} \label{eq:1deul_bw}
\end{equation}
where $\rho_l=\rho_m=\rho_r=1, \ u_l=u_m=u_r=0, \ p_l=10^3, \ p_m=10^{-2}, \ p_r=10^2$. The domain is $x \in [0,1]$ and the boundary conditions are \textit{reflecting solid walls }at $x_L=0$ and $x_R=1$. This problem is a challenging testing ground for any shock-capture scheme. The interaction of blast waves results in fine structures which can only be resolved accurately using a high-order scheme. However, the ensuing fluctuations in the flow fields can challenge the scheme's stability. In Fig.~\ref{fig:eulbw}, we plot the fields $\rho$ (left), $u$ (center), and $p$ (right) at consecutive times in panels $(a)-(d)$. The reference solution (black line) in this case, is obtained using HLLE scheme from \texttt{CLAWPACK} library on $N_x=8000$ cells. We perform all the spectral computations on a higher-resolution grid of $N_x=1599$ points. Despite the increase in resolution, the dealiased PPS approximation destabilises instantly and is thus not shown here. The high-order SR-DlVP approximations also destabilise early and render the scheme unusable for virtually every choice of $(\alpha,\gamma)$. Indeed, at early times, the SR-DlVP approximation of the blast-waves solution contains oscillations near the discontinuities. For the density $\rho$ and the pressure $p$, these oscillations produce negative values which are unphysical. Then, we can no longer compute $c=\sqrt{{\gamma p}/{\rho}}$ in Eq.~\eqref{eq:eigvals}, which leads to a breakdown of the spectral scheme. Here, we note that SR-DlVP approximations for the other test cases presented above do not destabilise in this manner because the oscillations in the fields $\rho$ and $p$  do not become negative-valued. In fact, the amplitude of oscillations near a shock is related to the amplitude of the shock (a fraction of the shock amplitude). For the blast-waves initial configuration, at the location $x=0.1$, there is an initial pressure discontinuity of amplitude of the order of $10^3$ as we go from $p_l= 10^{3}$ to $p_m= 10^{-2}$ in Eq.~\eqref{eq:1deul_bw}, while oscillations of the order of $10^{-2}$ are observed for $p(x\simeq 0.1,t)$ for $0<t\ll 1$. For small times $0<t\ll 1$, we obtain $p(0.1<x<0.9,t)\sim p_m = 10^{-2}$ by time continuity of the solution; then oscillations which develop in this region occasionally result in negative values of the pressure $p$ that destabilise the scheme. For other test cases, the discontinuities are not as strong and the resulting oscillations remain bounded and positive in time. In Fig~\ref{fig:eulbw}, we present SR-FeKo approximations with $(\alpha=1.35,\gamma=0.99)$. The approximate solution provided by the low-order Fej\'er-Korovkin kernel is non-oscillatory and positive. Since the SR-FeKo scheme preserves the monotonicity and in particular the positivity of the approximate solution, this scheme remains stable and noticeably accurate until the time $t=0.032$ (in panels $(a)-(c)$). This low-order approximation fails,  however, to resolve the fine structures following the interaction, especially for the density $\rho(x,t)$ in panel $(d)$ of Fig.~\ref{fig:eulbw}. We do not use $2/3$-dealiasing for the SR schemes presented above. 

We conclude that the blast-waves interaction problem presents a challenging test case for purely spectral methods (such as the methods SR, SP, SVV, ...), which employ global spectral interpolations such as the trigonometric interpolation. To our knowledge, shock-capturing spectral schemes in literature, do not perform tests with this standard problem. In~\cite{cai1993uniform}, the authors report a good approximation of the blast-waves solution, but this numerical scheme is not a purely spectral method since a piecewise polynomial (essentially non-oscillatory Lagrange polynomial interpolation) is used  to interpolate the function near discontinuities, while spectral interpolation is used only in the smooth parts of the flow.
\end{itemize}   
\end{subequations}


\subsection{1D wall approximation of the 3D-axisymmetric wall-bounded incompressible Euler equation}
\label{subsec:Results_lhmodel}
 
In this section, we apply spectral relaxation schemes to the 1D HL model, which is the wall approximation of the 3D-axisymmetric wall-bounded incompressible Euler equation, investigated in~\cite{luo2014potentially,1d_1,choi2017finite}. This model was proposed to capture the singular dynamics, at the wall, of the 3D-axisymmetric incompressible Euler equation and is given by
\begin{subequations}
\begin{align}
        \partial_t u +v \partial_x u &= 0, \label{eq:lh1du} \\
        \partial_t \omega +v \partial_x \omega &= \partial_x u, \label{eq:lh1dw} \\
        \partial_x v &= \mathcal{H} (\omega), \label{eq:lh1dv}
\end{align}
\label{eq:lh1d_full}
\end{subequations}
where $\mathcal{H}(\cdot)$ is the Hilbert transform. Equation~\eqref{eq:lh1dv} represents the Biot-Savart law for this model. For more details on the formulation of this model, see~\cite{choi2017finite,1d_1,kolluru2022insights}. Reference~\cite{kolluru2022insights} used a $2/3$-dealiased Fourier pseudospectral scheme to study the singular dynamics of this model for initial data $u(x,t=0)=10^4 \sin^2(2\pi x/L)$ and $\omega(x,t=0)=0$. Periodic boundary conditions were used on the domain $x \in [0,L]$ where $L=1/6$. The authors reported the emergence of tygers akin to those observed for the 1D inviscid Burgers~\cite{ray2011resonance} and the 3D-axisymmetric wall-bounded incompressible Euler flow~\cite{kolluru2022insights}.  

\begin{figure}[h!]
    \centering
    \begin{tikzpicture}
	\draw (0,0) node[inner sep=0]{\includegraphics[width=0.9\linewidth]{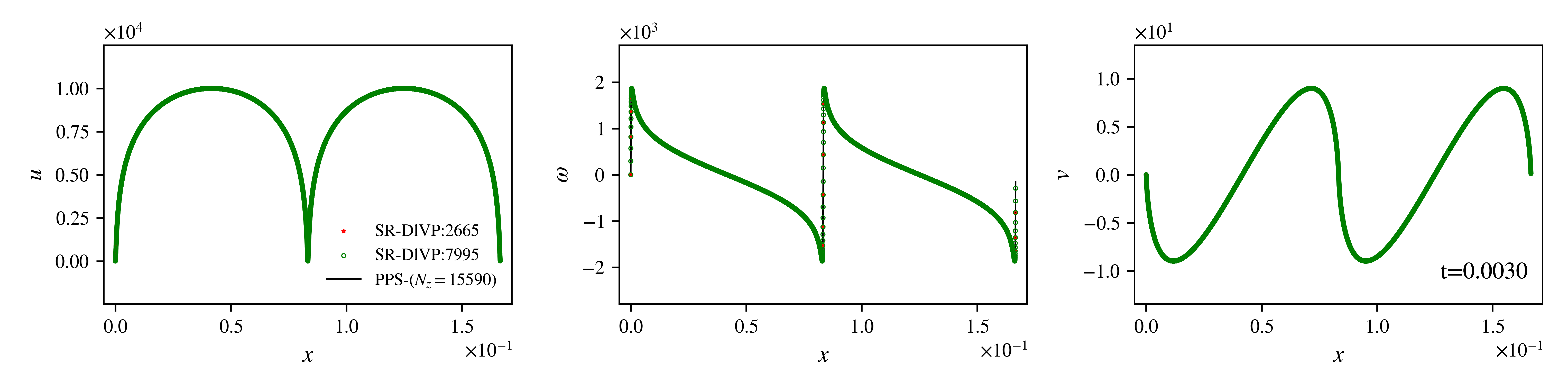}};
		\draw (-7.2,1.7) node {$(a)$};
	\end{tikzpicture}
    \begin{tikzpicture}
 	\draw (0,0) node[inner sep=0]{\includegraphics[width=0.9\linewidth]{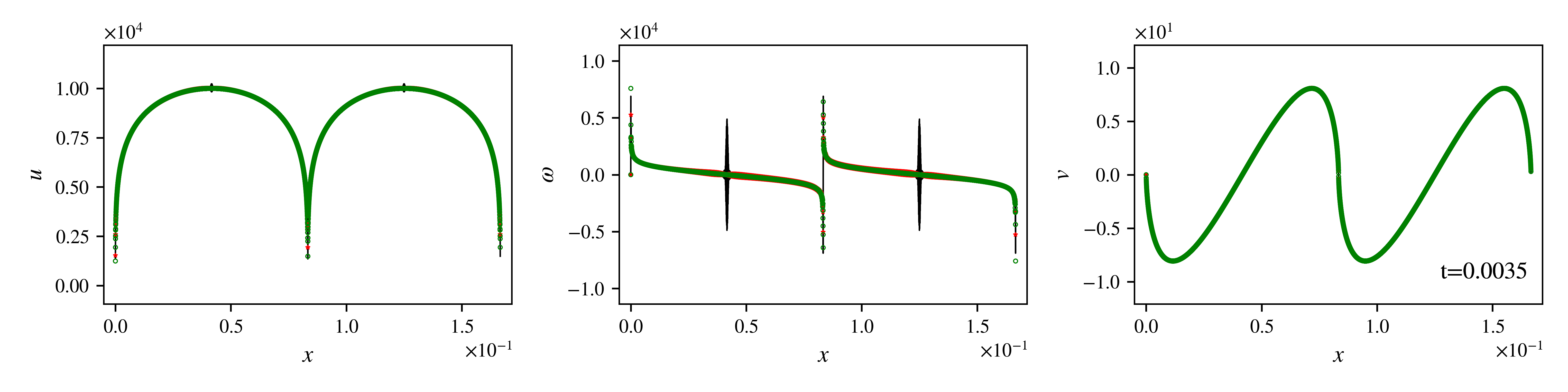}};
		\draw  (-7.2,1.7) node {$(b)$};
	\end{tikzpicture}   
    \begin{tikzpicture}
 	\draw (0,0) node[inner sep=0]{\includegraphics[width=0.9\linewidth]{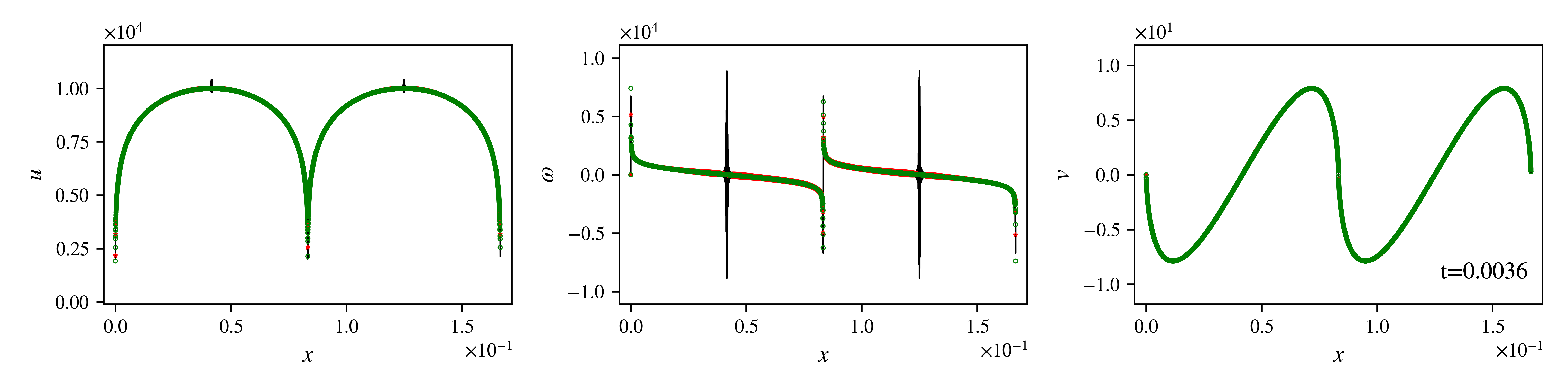}};
		\draw (-7.2,1.7) node {$(c)$};
  	\end{tikzpicture}  
       \begin{tikzpicture}
 	\draw (0,0) node[inner sep=0]{\includegraphics[width=0.9\linewidth]{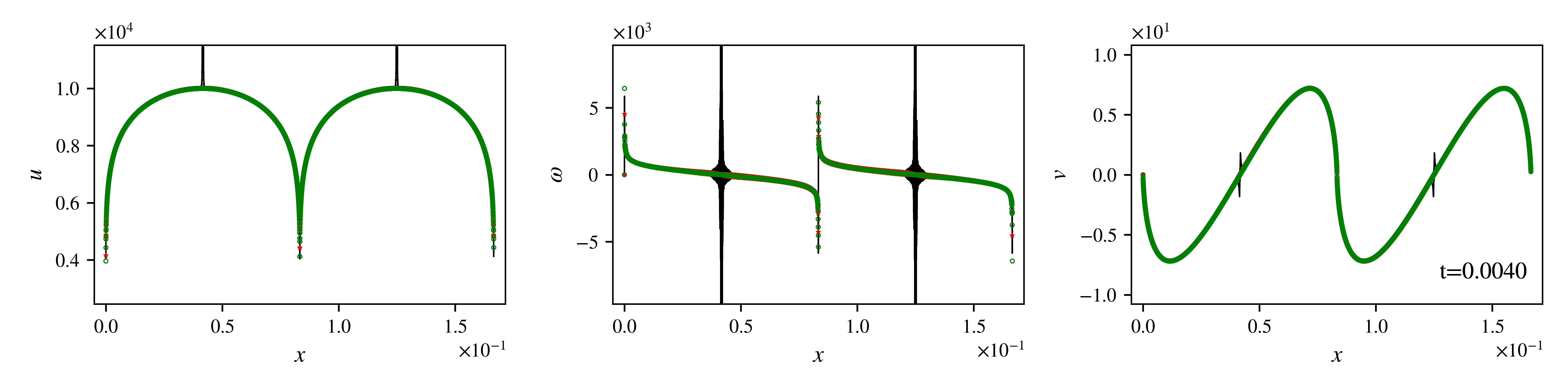}};
		\draw (-7.2,1.7) node {$(d)$};
  	\end{tikzpicture}  
 \caption{ HL model: Plot vs. $x$ of $u,\omega$, and $v$ fields. Different panels give the fields at different instants of time $t=0.003,0.0035,0.0036,0.004$ from panel $(a)-(d)$ respectively. The SR-DlVP runs are given for $N_x=2665$ (red stars) and $N_x=7995$ (green circles) resolutions. The best value of the parameters used in the spectral relaxation scheme are $(\alpha=1.6, \gamma=0.99, n=0.92)$. The black line for comparison represents the PPS approximation of a high resolution $N_x=15590$. The time of singularity $t_*=0.0035056$. Dealiasing is not performed for SR runs.}
    \label{fig:axeul}
\end{figure}
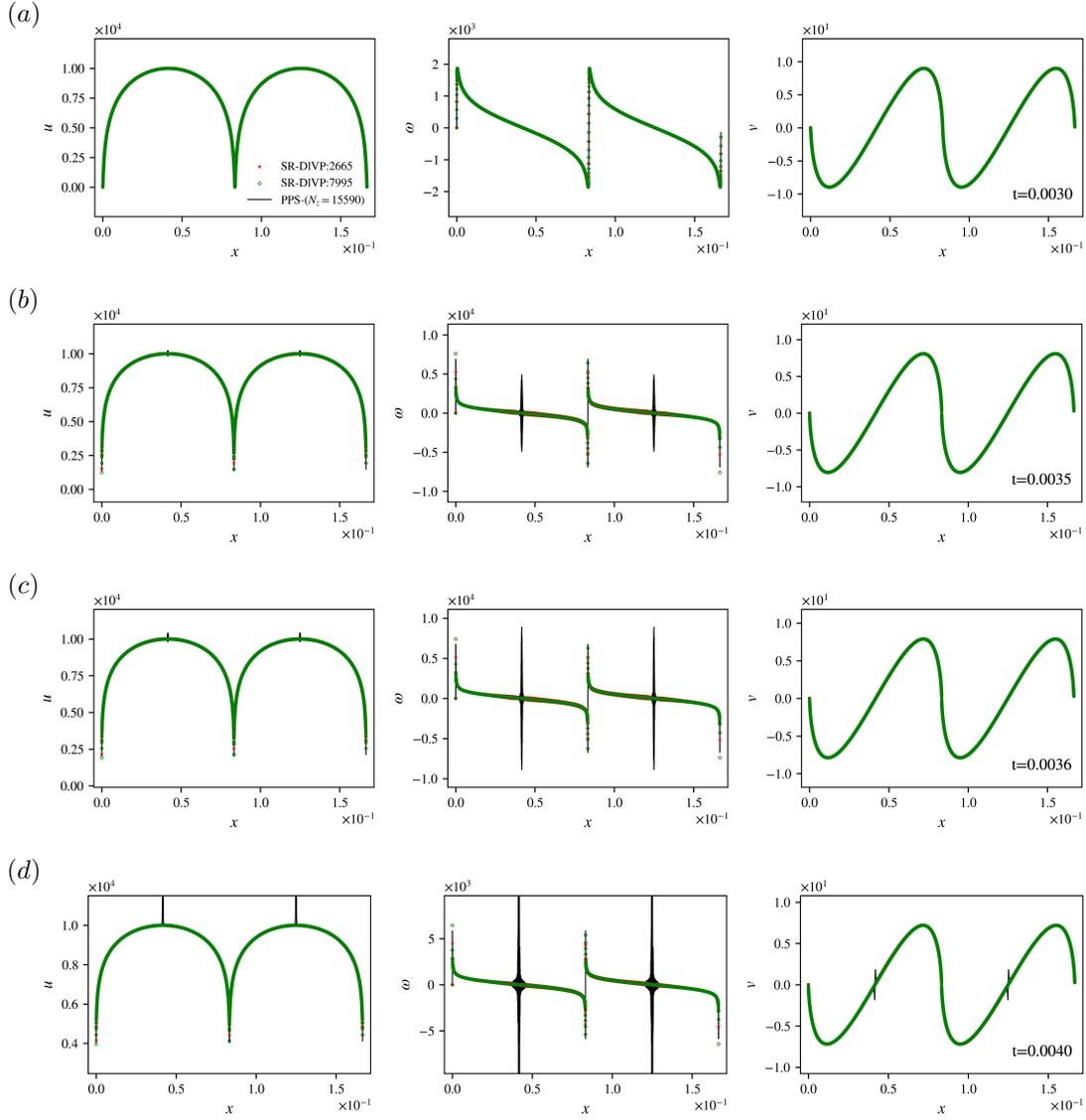

In Fig.~\ref{fig:axeul}, we plot spectral approximations of the singular solution of Eqs.~\eqref{eq:lh1d_full}; we now discuss them in detail. Since we do not have the exact solution for this model, we use a dealiased PPS approximation computed on a very high resolution grid with $N_x=15990$ (solid black line) as the reference solution. We do not show the PPS approximations of lower resolutions for greater readability. For the high-resolution PPS approximation, we see the gradual emergence of tygers in smooth parts of the solution; first for the field $\omega(x,t)$ (center column of panel $(b)$), followed by $u(x,t)$ (left column of panel $(c)$), and finally $v(x,t)$ (right column of panel $(d)$). In PPS approximations of lower resolutions, tygers develop earlier in time but at the same locations. The link between tyger-birth time and the spatial resolution of PPS scheme is discussed in~\cite{ray2011resonance,kolluru2022insights}. Briefly, tygers emerge in the system when the precipitating complex singularity moves into a strip of width $\Delta x$ which flanks   the real-line, in complex space. Here $\Delta x \sim 1/N_x$ is the smallest grid spacing in the PPS discretisation. Increasing $N_x$ reduces the width of this strip; thus, it takes longer for the singularity to encounter the strip and for tyger-birth to occur.\par

We plot SR-DlVP approximations with $(\alpha=1.6,\gamma=0.99)$ for two spatial resolutions: $N_x=2665$ (SR-DlVP:2665 in red) and $N_x=7995$ (SR-DlVP:7995 in green). In panel $(a)$ at time $t=0.0030$, both SR-DlVP approximations closely approximate the reference solution (black line). As we get close to the time of singularity in this model $t_*=0.0035056$, in panels $(b)$ and $(c)$ of Fig.~\ref{fig:axeul}, the SR-DlVP:7995 (green) approximation maintains a smooth tyger-free profile while also capturing the discontinuity at $x=0,L/2,L$ in $\omega(x,t)$ without Gibbs oscillations. The lower resolution approximation, SR-DlVP:2665 (red) stays close to the reference solution but does show attenuated tygers. 

As discussed in Section~\ref{sec:intro}, by using the analyticity strip technique, we can track the nearest complex singularity, which can ultimately precipitate finite-time blowup for Eq.~\eqref{eq:lh1d_full}. If a field $u(z,t)$ contains a pole, of order $\mu$, close to the real-line at $z_*(t) =( x_* + i \delta) \in \mathbb{C}$, the Fourier spectrum displays a strong exponentially-decaying tail in the high-$k$ region. It can be shown that $\hat{u}(k) \sim |k|^{-(\mu+1)}e^{-k\delta}e^{ikx_*}$ as $k \rightarrow \infty$~\cite{carrier2005functions,sulem1983tracing}. We can thus extract $\delta(t)$ from the slope of the exponentially-decaying tail of the spectrum. However, at times close to the time of singularity, the loss of spectral convergence in the PPS approximation and subsequent thermalisation prevents us from getting reliable estimates for $\delta(t)$. We can delay the loss of spectral convergence, and thermalisation, by using the dissipative SR schemes. Furthermore, SR schemes allow us to retain almost all of the wavenumbers at a given resolution, as they avoid $2/3$-dealiasing.\par

\begin{figure}[h!]
    \centering
    \begin{tikzpicture}
	\draw (0,0) node[inner sep=0]{\includegraphics[width=0.9\linewidth]{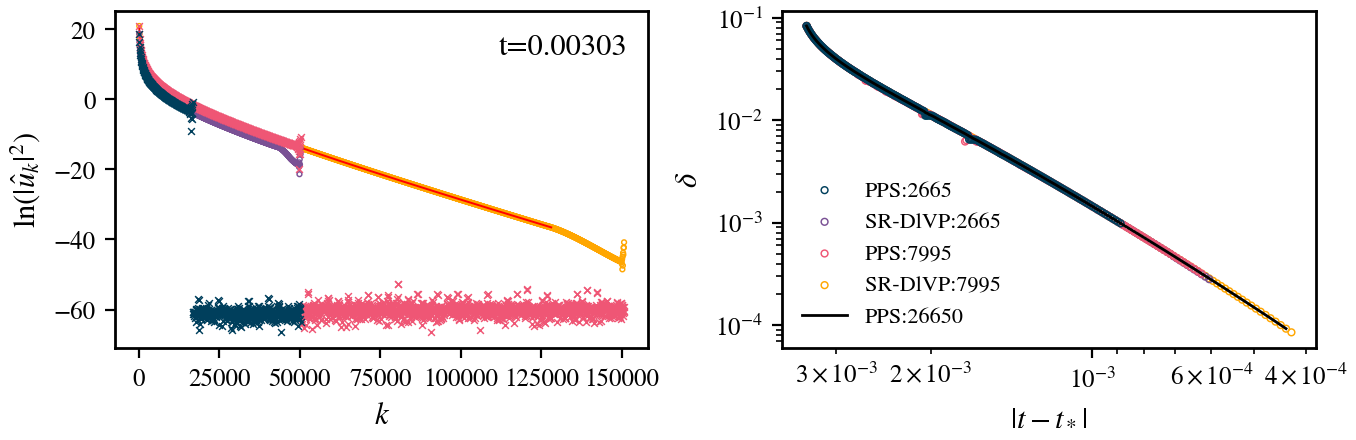}};
 		\draw (-1.1,1.3) node {\large$(a)$};
		\draw (7.0,1.3) node {\large$(b)$};
	\end{tikzpicture}
 \caption{ HL model: (Left) Plots vs. $k$ of $\ln |\hat{u}_k|^2$ and (Right) plots vs. $|t-t_*|$ of the width of the analyticity strip $\delta$ extracted from the exponential slope of the corresponding spectra, displayed on the left. In panel $(a)$, the spectra of SR-DlVP approximations with $(\alpha=1.6,\gamma=0.99)$ extend over a much larger range of wavenumbers $k$ compared to their PPS counterparts. In panel $(b)$, $\delta(t)$ estimates from SR-DlVP and PPS approximations are superposed. The reference $\delta$ is obtained using PPS approximation on a very high spatial resolution $N_x=26650$ points. Here, $t_*=0.0035056$ is the time of singularity in this model. }
    \label{fig:axeul_spectra}
\end{figure}

In panel $(a)$ of Fig.~\ref{fig:axeul_spectra}, we show the logarithm of the Fourier spectra, $\ln |\widehat{u}_k|^2$, at time $t=0.00303$ before the time of singularity. We compare the decaying tails of the spectra of SR-DlVP and PPS schemes for two different resolutions. PPS:2665 (blue) and PPS:7995 (pink) are $2/3$-dealiased. Thus, they do not use the full range of wavenumbers available to them. The corresponding SR-DlVP approximations with $N_x=2665$ (violet) and $N_x=7995$ (orange) do not have this constraint; the exponential tail is registered over a much larger range of wavenumbers for the SR approximations. Having a large range of fitting is key to a reliable estimation of $\delta(t)$. In panel $(a)$, the spectrum of SR-DlVP:2665 follows that of the higher resolution PPS:7995 approximation. However, there is a slope change at the end of the exponential tail; this coincides with the window over which DlVP filter is dissipative (negative slope of the filter in Fourier space). For a given resolution $N_x$, the DlVP kernel with $r=0.92$, acts as an almost low-pass filter with a cutoff at $k= \tfrac{2\pi}{L} \times 0.92 \times (N/2 -1)^\gamma \approx 43000 $ for $N_x=2665$. In the spectra of SR-DlVP:2665 the slope change occurs around $k \sim 40000$ as read from the Fig.~\ref{fig:axeul_spectra} and is consistent with our estimate. We plot the SR-DlVP:7995 (orange) spectrum to once again demonstrate the expansion of the range of fitting for $\delta(t)$ compared to PPS:7995 (pink). 

 We now compare the estimate for the width of the analyticity-strip $\delta(t)$ extracted from SR-DlVP approximations and PPS-approximations. In panel $(b)$ of Fig.~\ref{fig:axeul_spectra}, we show the $\delta(t)$ extracted from the corresponding spectra in panel $(a)$. In SR-DlVP spectra, the fit-window used for determination of $\delta$ is restricted to exclude the range of Fourier modes over which the de La Vall\'ee Poussin kernel has a negative slope; the fit ends before the slope change (red line over SR-DlVP:7995 shown in orange). Despite the restriction imposed above, we still retain a much larger range of available wavenumbers when compared to PPS approximations of the same resolutions. 
We use the $\delta$ obtained from a very high-resolution PPS approximation using $N_x=26650$ as a reference (solid black line). 
The estimates of $\delta$ extracted from the SR-DlVP spectra for the resolution of $N_x=2665$ (violet) and the higher resolution of $N_x=7995$ (orange) overlie the black line and are consistently accurate for longer times than analogous PPS estimates for the same resolution. Therefore, we conclude that SR-DlVP schemes are effective in increasing the time until which analyticity-strip techniques can reliably be used to track singularities in a given model. 


\section{Conclusions}
\label{sec:conclusions}
We have presented an extensive numerical study of the novel spectral relaxation and spectral purging schemes for the 1D inviscid Burgers equation in a periodic domain. Here, we have investigated the role of different kinds of kernels when used in the spectral relaxation and purging schemes. We have carried out numerical convergence analysis using a low-order kernel and a high-order kernel; we observe increased order of convergence after the shock formation for the novel schemes in both cases when compared to pure pseudospectral methods. The SR approximations obtained using high-order de La Vall\'ee Poussin kernel is similar to the one obtained using SVV methods; both approximations show bounded oscillations near the shock. SR schemes that use the low-order Fej\'e-Korovkin kernel provide a non-oscillatory and maximum-principle preserving approximation to the entropic solution for the 1D inviscid Burgers equation. We have performed a detailed study of the aliasing instability in the evolution of a single-mode initial condition for the 1D inviscid Burgers equation; we conclude that dealiasing homogenises the growth-rate of aliasing errors across initial conditions, but can be expensive. Dealiasing can stabilise SR schemes at smaller values of $\alpha$; the resulting scheme is less-diffusive, with a higher order of convergence.  

We have numerically extended the application of the novel schemes to systems of 1D nonlinear hyperbolic systems of conservation laws --- namely, the 1D shallow water equations and the 1D compressible Euler equations of gas dynamics. Here, we demonstrate the performance of SR approximations using F\'ejer--Korovkin and de La Vall\'ee Poussin kernels for a variety of test problems. Except for the interaction of blast waves for the 1D compressible Euler equations, the SR approximations were capable of accurately resolving the location and profile of the shock and other sharp structures in the flow fields. The ability to use a range of kernels that best fit the attributes and properties of the solutions at hand, allows for greater flexibility and power in the SR schemes. This is, in particular, best observed for the case of the Shu--Osher initial data evolved with the 1D Euler equation. The use of a higher-order kernel helps retain the oscillatory structures of the flow that do not respect monotony, while capturing the shock profile and its evolution in time. The use of high-order kernels allows us to achieve higher accuracy and confirm visual convergence to shock capture, while preserving the spectral convergence in the smooth parts of the solution away from the discontinuity. They, however, result in the reappearance of the Gibbs oscillations in a reduced neighbourhood of the discontinuity. This is in line with what has been observed for the SVV method and for filtering approaches that use high-order kernels~\cite{hesthaven2017numerical}. We note here that the parameters $(\alpha,\gamma)$, that minimise artefacts and best capture shock profiles, vary based on the problem at hand, the initial data, as well as the spatial resolution used in a given computation.\par

We also note that the extension of SR and SP schemes to Chebyshev pseudospectral methods is easy; we use Chebyshev SR schemes to model challenging test problems for the 1D compressible Euler equations. We do not need to perform additional boundary treatment for SR schemes; the CCM method for pure pseudospectral methods works accurately. Therefore, we conclude that spectral relaxation and purging schemes are very effective at capturing the shock profiles in discontinuous solutions that are known otherwise to challenge traditional spectral methods.\par

Finally, we applied the SR schemes for the approximation of singular flows in the 1D HL model. Here, we find that SR approximations using the high-order de La Vall\'ee Poussin kernel are capable of following the reference solution closely, while remaining free from Gibbs oscillations and tygers. The estimates obtained for $\delta(t)$ using the SR approximations agreed with the reference to an appreciable degree. This is a key success for SR schemes; they can be used for tracking the finite-time singularity for much longer times compared to a traditional pseudospectral scheme that has the same spatial resolution. We highlight, however, that the choice of kernel and parameters is crucial to the extension and has to be checked by using resolution studies to achieve a good scaling regime. We are yet to test the method with systems where there is a paradigm change in the decay of $\delta(t)$. This will be the subject of future work, in addition to the application of this method for three-dimensional PDEs. 

The parameters $(\alpha,\gamma)$ used for the best approximations were chosen manually for all test problems, except for the case of the 1D inviscid Burgers equation. This was done so as to ensure maximum overlap between the reference solution and the approximate one, while minimising any spurious oscillations. Further tuning of the parameters $(\alpha,\gamma)$ may produce better approximations for all the cases that we have discussed. Optimisation strategies such as the gradient descent method can be used to automate the search for the optimal parameter set, once we determine the optimal range. For this step, machine-learning strategies may also be employed.\par

\section*{Acknowledgement}
The authors thank SERB, CSIR, NSM, UGC (India), Indo--French Center for Applied Mathematics (IFCAM) and Universit\'e C\^ote d'Azur (UniCA) for their support. This work has been supported by the French government through the UniCA\textsuperscript{JEDI} investments in the Future project managed by the National Research Agency (ANR) with reference number ANR-15-IDEX-01. 

\printbibliography

@book{canuto2007spectral,
  title={Spectral methods: evolution to complex geometries and applications to fluid dynamics},
  author={Canuto, Claudio and Hussaini, M Yousuff and Quarteroni, Alfio and Zang, Thomas A},
  year={2007},
  publisher={Springer Science \& Business Media}
}

@book{canuto2006spectral,
  title={Spectral methods: fundamentals in single domains},
  author={Canuto, Claudio and Hussaini, M Yousuff and Quarteroni, Alfio and Zang, Thomas A},
  year={2007},
  publisher={Springer Science \& Business Media}
}

@article{gottlieb1981spectral,
  title={Spectral calculations of one-dimensional inviscid compressible flows},
  author={Gottlieb, David and Lustman, Liviu and Orszag, Steven A},
  journal={SIAM Journal on Scientific and Statistical Computing},
  volume={2},
  number={3},
  pages={296--310},
  year={1981},
  publisher={SIAM}
}

@unpublished{bessemain,
  title={Analysis of the spectral relaxation method and the spectral purging method for nonlinear conservation laws},
  author={Besse, Nicolas}, 
note   = {Submitted.}}

@book{leveque2002finite,
  title={Finite volume methods for hyperbolic problems},
  author={LeVeque, Randall J},
  volume={31},
  year={2002},
  publisher={Cambridge university press}
}

@article{thompson1990time,
  title={Time-dependent boundary conditions for hyperbolic systems, II},
  author={Thompson, Kevin W},
  journal={Journal of computational physics},
  volume={89},
  number={2},
  pages={439--461},
  year={1990},
  publisher={Elsevier}
}

@article{kosloff1993modified,
  title={A modified Chebyshev pseudospectral method with an O (N-1) time step restriction},
  author={Kosloff, Dan and Tal-Ezer, Hillel},
  journal={Journal of Computational Physics},
  volume={104},
  number={2},
  pages={457--469},
  year={1993},
  publisher={Elsevier}
}

@article{sod1978survey,
  title={A survey of several finite difference methods for systems of nonlinear hyperbolic conservation laws},
  author={Sod, Gary A},
  journal={Journal of computational physics},
  volume={27},
  number={1},
  pages={1--31},
  year={1978},
  publisher={Elsevier}
}

@article{lax1954weak,
  title={Weak solutions of nonlinear hyperbolic equations and their numerical computation},
  author={Lax, Peter D},
  journal={Communications on pure and applied mathematics},
  volume={7},
  number={1},
  pages={159--193},
  year={1954},
  publisher={Wiley Online Library}
}

@article{tadmor2005adaptive,
  title={Adaptive filters for piecewise smooth spectral data},
  author={Tadmor, Eitan and Tanner, Jared},
  journal={IMA journal of numerical analysis},
  volume={25},
  number={4},
  pages={635--647},
  year={2005},
  publisher={Oxford University Press}
}

@article{tadmor1990shock,
  title={Shock capturing by the spectral viscosity method},
  author={Tadmor, Eitan},
  journal={Computer Methods in Applied Mechanics and Engineering},
  volume={80},
  number={1-3},
  pages={197--208},
  year={1990},
  publisher={Elsevier}
}

@article{tadmor1993super,
  title={Super viscosity and spectral approximations of nonlinear conservation laws},
  author={Tadmor, Eitan},
  journal={Numerical Methods for Fluid Dynamics},
  volume={4},
  pages={69--81},
  year={1993},
  publisher={Clarendon Press Oxford}
}

@article{maday1993legendre,
  title={Legendre pseudospectral viscosity method for nonlinear conservation laws},
  author={Maday, Yvon and Kaber, Sidi M Ould and Tadmor, Eitan},
  journal={SIAM Journal on Numerical Analysis},
  volume={30},
  number={2},
  pages={321--342},
  year={1993},
  publisher={SIAM}
}

@article{guo2001spectral,
  title={Spectral vanishing viscosity method for nonlinear conservation laws},
  author={Guo, Ben-yu and Ma, He-ping and Tadmor, Eitan},
  journal={SIAM journal on numerical analysis},
  volume={39},
  number={4},
  pages={1254--1268},
  year={2001},
  publisher={SIAM}
}

@article{majda1978fourier,
  title={The Fourier method for nonsmooth initial data},
  author={Majda, Andrew and McDonough, James and Osher, Stanley},
  journal={Mathematics of Computation},
  volume={32},
  number={144},
  pages={1041--1081},
  year={1978}
}

@article{karamanos2000spectral,
  title={A spectral vanishing viscosity method for large-eddy simulations},
  author={Karamanos, GE and Karniadakis,GS},
  journal={Journal of Computational Physics},
  volume={163},
  number={1},
  pages={22--50},
  year={2000},
  publisher={Elsevier}
}

@article{sun2002practical,
  title={A practical spectral method for hyperbolic conservation laws},
  author={Sun, Yu-Hui and Zhou, YC and Wei, GW},
  journal={arXiv preprint math/0212346},
  year={2002}
}

@article{gu2003conjugate,
  title={Conjugate filter approach for shock capturing},
  author={Gu, Yun and Wei, GW},
  journal={Communications in Numerical Methods in Engineering},
  volume={19},
  number={2},
  pages={99--110},
  year={2003},
  publisher={Wiley Online Library}
}

@article{liwen2004regularised,
  title={The regularized Whittaker-Kotel'nikov-Shannon sampling theorem and its application to the numerical solutions of partial differential equations},
  author={Liwen, Qian},
  year={2004},
note={PhD thesis}
}

@article{butzer1971fourier,
  title={Fourier analysis and approximation, Vol. 1},
  author={Butzer, Paul L and Nessel, Rolf J},
  journal={Reviews in Group Representation Theory, Part A (Pure and Applied Mathematics Series, Vol. 7)},
  year={1971},
  publisher={Marcel Dekker Inc. New York, NY, USA}
}

@article{woodward1984numerical,
  title={The numerical simulation of two-dimensional fluid flow with strong shocks},
  author={Woodward, Paul and Colella, Phillip},
  journal={Journal of computational physics},
  volume={54},
  number={1},
  pages={115--173},
  year={1984},
  publisher={Elsevier}
}

@article{brehm2015comparison,
  title={A comparison of higher-order finite-difference shock capturing schemes},
  author={Brehm, Christoph and Barad, Michael F and Housman, Jeffrey A and Kiris, Cetin C},
  journal={Computers \& Fluids},
  volume={122},
  pages={184--208},
  year={2015},
  publisher={Elsevier}
}

@article{zhao2019general,
  title={A general framework for the evaluation of shock-capturing schemes},
  author={Zhao, Guoyan and Sun, Mingbo and Memmolo, Antonio and Pirozzoli, Sergio},
  journal={Journal of Computational Physics},
  volume={376},
  pages={924--936},
  year={2019},
  publisher={Elsevier}
}

@article{moura2016eigensolution,
  title={Eigensolution analysis of spectral/hp continuous Galerkin approximations to advection--diffusion problems: Insights into spectral vanishing viscosity},
  author={Moura, Rodrigo Costa and Sherwin, Spencer J and Peir{\'o}, Joaquim},
  journal={Journal of Computational Physics},
  volume={307},
  pages={401--422},
  year={2016},
  publisher={Elsevier}
}

@article{klockner2011viscous,
  title={Viscous shock capturing in a time-explicit discontinuous Galerkin method},
  author={Kl{\"o}ckner, Andreas and Warburton, Tim and Hesthaven, Jan S},
  journal={Mathematical Modelling of Natural Phenomena},
  volume={6},
  number={3},
  pages={57--83},
  year={2011},
  publisher={EDP Sciences}
}

@article{hartmann2006adaptive,
  title={Adaptive discontinuous Galerkin methods with shock-capturing for the compressible Navier--Stokes equations},
  author={Hartmann, Ralf},
  journal={International Journal for Numerical Methods in Fluids},
  volume={51},
  number={9-10},
  pages={1131--1156},
  year={2006},
  publisher={Wiley Online Library}
}

@article{kumar2019efficient,
  title={Efficient seventh order WENO schemes of adaptive order for hyperbolic conservation laws},
  author={Kumar, Rakesh and Chandrashekar, Praveen},
  journal={Computers \& Fluids},
  volume={190},
  pages={49--76},
  year={2019},
  publisher={Elsevier}
}

@article{gottlieb1997gibbs,
  title={On the Gibbs phenomenon and its resolution},
  author={Gottlieb, David and Shu, Chi-Wang},
  journal={SIAM review},
  volume={39},
  number={4},
  pages={644--668},
  year={1997},
  publisher={SIAM}
}

@article{kreiss1979stability,
  title={Stability of the Fourier method},
  author={Kreiss, Heinz-Otto and Oliger, Joseph},
  journal={SIAM Journal on Numerical Analysis},
  volume={16},
  number={3},
  pages={421--433},
  year={1979},
  publisher={SIAM}
}

@inproceedings{hussaini1983spectral,
  title={Spectral methods for the Euler equations},
  author={Hussaini, M and Kopriva, D and Salas, M and Zang, T},
  booktitle={6th Computational Fluid Dynamics Conference Danvers},
  pages={1942},
  year={1983}
}

@article{gottlieb2001spectral,
  title={Spectral methods for hyperbolic problems},
  author={Gottlieb, David and Hesthaven, Jan S},
  journal={Journal of Computational and Applied Mathematics},
  volume={128},
  number={1-2},
  pages={83--131},
  year={2001},
  publisher={Elsevier}
}

@article{brachet2013ideal,
  title = {Ideal evolution of magnetohydrodynamic turbulence when imposing Taylor-Green symmetries},
  author = {Brachet, M. E. and Bustamante, M. D. and Krstulovic, G. and Mininni, P. D. and Pouquet, A. and Rosenberg, D.},
  journal = {Phys. Rev. E},
  volume = {87},
  issue = {1},
  pages = {013110},
  numpages = {14},
  year = {2013},
  publisher = {American Physical Society},
  doi = {10.1103/PhysRevE.87.013110},
  url = {https://link.aps.org/doi/10.1103/PhysRevE.87.013110}
}

@article{pereira2023adaptive,
  title={Are Adaptive Galerkin Schemes Dissipative?},
  author={Pereira, Rodrigo M and Nguyen van yen, Natacha and Schneider, Kai and Farge, Marie},
  journal={SIAM Review},
  volume={65},
  number={4},
  pages={1109--1134},
  year={2023},
  publisher={SIAM}
}

@article{kolluru2022insights,
  title={Insights from a pseudospectral study of a potentially singular solution of the three-dimensional axisymmetric incompressible Euler equation},
  author={Kolluru, Sai Swetha Venkata and Sharma, Puneet and Pandit, Rahul},
  journal={Physical Review E},
  volume={105},
  number={6},
  pages={065107},
  year={2022},
  publisher={APS}
}

@misc{clawpack,
  title={Clawpack version 4.0 user’s guide},
  author={LeVeque, Randall J},
url={\url{http://depts.washington.edu/clawpack/clawpack-4.3/}{CLAWPACK 4.3}},
year={2003},
  publisher={University of Washington},

}

@misc{hypar,
  title = {HyPar - Finite-Difference Hyperbolic-Parabolic PDE Solver on Cartesian Grids },
  author={Ghosh, Debojyoti},
  howpublished = {\url{https://bitbucket.org/deboghosh/hypar}},
}

@article{cai1989essentially,
  title={Essentially nonoscillatory spectral Fourier methods for shock wave calculations},
  author={Cai, Wei and Gottlieb, David and Shu, Chi-Wang},
  journal={Mathematics of Computation},
  volume={52},
  number={186},
  pages={389--410},
  year={1989}
}

@article{cai1993uniform,
  title={Uniform high-order spectral methods for one-and two-dimensional Euler equations},
  author={Cai, Wei and Shu, Chi-Wang},
  journal={Journal of Computational Physics},
  volume={104},
  number={2},
  pages={427--443},
  year={1993},
  publisher={Elsevier}
}

@inproceedings{gottlieb1985recovering,
  title={Recovering pointwise values of discontinuous data within spectral accuracy},
  author={Gottlieb, David and Tadmor, Eitan},
  booktitle={Progress and supercomputing in computational fluid dynamics},
  pages={357--375},
  year={1985},
  organization={Springer}
}

@article{canuto1982approximation,
  title={Approximation results for orthogonal polynomials in Sobolev spaces},
  author={Canuto, Claudio and Quarteroni, Alfio},
  journal={Mathematics of Computation},
  volume={38},
  number={157},
  pages={67--86},
  year={1982}
}

@book{kreiss1973methods,
  title={Methods for the approximate solution of time dependent problems},
  author={Kreiss, Heinz and Oliger, Joseph},
  number={10},
  year={1973},
  publisher={International Council of Scientific Unions, World Meteorological Organization}
}

@article{tadmor1986exponential,
  title={The exponential accuracy of Fourier and Chebyshev differencing methods},
  author={Tadmor, Eitan},
  journal={SIAM Journal on Numerical Analysis},
  volume={23},
  number={1},
  pages={1--10},
  year={1986},
  publisher={SIAM}
}

@article{maday1989analysis,
  title={Analysis of the spectral vanishing viscosity method for periodic conservation laws},
  author={Maday, Yvon and Tadmor, Eitan},
  journal={SIAM Journal on Numerical Analysis},
  volume={26},
  number={4},
  pages={854--870},
  year={1989},
  publisher={SIAM}
}

@book{hesthaven2017numerical,
  title={Numerical methods for conservation laws: From analysis to algorithms},
  author={Hesthaven, Jan S},
  year={2017},
  publisher={SIAM}
}

@article{sun2006windowed,
  title={A windowed Fourier pseudospectral method for hyperbolic conservation laws},
  author={Sun, Yuhui and Zhou, YC and Li, Shu-Guang and Wei, Gang W},
  journal={Journal of Computational Physics},
  volume={214},
  number={2},
  pages={466--490},
  year={2006},
  publisher={Elsevier}
}

@article{hussaini1985spectral,
  title={Spectral methods for the Euler equations. i-Fourier methods and shock capturing},
  author={Hussaini, MY and Kopriva, DA and Salas, MD and Zang, TA},
  journal={AIAA journal},
  volume={23},
  number={1},
  pages={64--70},
  year={1985}
}

@article{tadmor1989convergence,
  title={Convergence of spectral methods for nonlinear conservation laws},
  author={Tadmor, Eitan},
  journal={SIAM Journal on Numerical Analysis},
  volume={26},
  number={1},
  pages={30--44},
  year={1989},
  publisher={SIAM}
}

@article{vandeven1991family,
  title={Family of spectral filters for discontinuous problems},
  author={Vandeven, Herv{\'e}},
  journal={Journal of Scientific Computing},
  volume={6},
  pages={159--192},
  year={1991},
  publisher={Springer}
}

@article{tadmor2007filters,
  title={Filters, mollifiers and the computation of the Gibbs phenomenon},
  author={Tadmor, Eitan},
  journal={Acta Numerica},
  volume={16},
  pages={305--378},
  year={2007},
  publisher={Cambridge University Press}
}

@article{luo2014potentially,
  title={Potentially singular solutions of the {3D} axisymmetric {Euler} equations},
  author={Luo, Guo and Hou, Thomas Y},
  journal={Proceedings of the National Academy of Sciences},
  volume={111},
  number={36},
  pages={12968--12973},
  year={2014},
  publisher={National Acad Sciences}
}

@article{frigo1999fftw,
  title={{FFTW} user’s manual},
  author={Frigo, Matteo and Johnson, Steven G},
  journal={Massachusetts Institute of Technology},
  year={1999}
}

@article{frisch2003singularities,
  title={Singularities of {Euler} flow? {Not} out of the blue!},
  author={Frisch, U and Matsumoto, T and Bec, J},
  journal={Journal of statistical physics},
  volume={113},
  number={5-6},
  pages={761--781},
  year={2003},
  publisher={Springer}
}

@article{sulem1983tracing,
  title={Tracing complex singularities with spectral methods},
  author={Sulem, Catherine and Sulem, Pierre-Louis and Frisch, H{\'e}l{\`e}ne},
  journal={Journal of Computational Physics},
  volume={50},
  number={1},
  pages={138--161},
  year={1983},
  publisher={Elsevier}
}

@article{cichowlas2005evolution,
  title={Evolution of complex singularities in {Kida}--{Pelz} and {Taylor}--{Green} inviscid flows},
  author={Cichowlas, C and Brachet, ME},
  journal={Fluid Dynamics Research},
  volume={36},
  number={4-6},
  pages={239},
  year={2005},
  publisher={IOP Publishing}
}

@article{gibbon2008euler,
  title={The three-dimensional {Euler} equations: how much do we know?},
  author={Gibbon, JD},
  journal={Physica D: Nonlinear Phenomena},
  volume={237},
  number={14-17},
  pages={1894--1904},
  year={2008},
  publisher={Elsevier}
}

@article{1d_1,
title={On the finite-time blowup of a {1D} model for the {3D} incompressible {Euler} equations},
author={Hou, Thomas Y and Luo, Guo},
journal={arXiv preprint arXiv:1311.2613},
year={2013}
}

@article{caflisch2015complex,
title = "Complex singularities and {PDE}s",
author = "Caflisch, {Russel E.} and Francesco Gargano and Marco Sammartino and Vincenzo Sciacca",
year = "2015",
volume = "6",
pages = "69--133",
journal = "Rivista di Matematica della Universita di Parma",
issn = "0035-6298",
publisher = "Universita degli Studi di Parma",
number = "1",

}

@article{krstulovic2009cascades,
  title={Cascades, thermalization, and eddy viscosity in helical Galerkin truncated Euler flows},
  author={Krstulovic, G and Mininni, Pablo Daniel and Brachet, ME and Pouquet, A},
  journal={Physical Review E},
  volume={79},
  number={5},
  pages={056304},
  year={2009},
  publisher={APS}
}

@article{feng2017thermalized,
  title={Thermalized solution of the Galerkin-truncated Burgers equation: From the birth of local structures to thermalization},
  author={Feng, Peihua and Zhang, Jiazhong and Cao, Shengli and Prants, Sergey V and Liu, Yan},
  journal={Communications in Nonlinear Science and Numerical Simulation},
  volume={45},
  pages={104--116},
  year={2017},
  publisher={Elsevier}
}

@article{kida1986study,
  title={Study of complex singularities by filtered spectral method},
  author={Kida, Shigeo},
  journal={Journal of the Physical Society of Japan},
  volume={55},
  number={5},
  pages={1542--1555},
  year={1986},
  publisher={The Physical Society of Japan}
}

@article{pereira2013wavelet,
  title={Wavelet methods to eliminate resonances in the Galerkin-truncated Burgers and Euler equations},
  author={Pereira, Rodrigo M and Farge, Marie and Schneider, K and others},
  journal={Physical Review E},
  volume={87},
  number={3},
  pages={033017},
  year={2013},
  publisher={APS}
}

@article{farge2017wavelet,
  title={Wavelet-based regularization of the Galerkin truncated three-dimensional incompressible Euler flows},
  author={Farge, Marie and Okamoto, Naoya and Schneider, Kai and Yoshimatsu, Katsunori},
  journal={Physical Review E},
  volume={96},
  number={6},
  pages={063119},
  year={2017},
  publisher={APS}
}

@article{ray2015thermalized,
  title={Thermalized solutions, statistical mechanics and turbulence: An overview of some recent results},
  author={Ray, Samriddhi Sankar},
  journal={Pramana},
  volume={84},
  number={3},
  pages={395--407},
  year={2015},
  publisher={Springer}
}

@article{ray2011resonance,
  title = {Resonance phenomenon for the {Galerkin}-truncated {Burgers} and {Euler} equations},
  author = {Ray, Samriddhi Sankar and Frisch, Uriel and Nazarenko, Sergei and Matsumoto, Takeshi},
  journal = {Phys. Rev. E},
  volume = {84},
  issue = {1},
  pages = {016301},
  numpages = {21},
  year = {2011},
  publisher = {American Physical Society},
  doi = {10.1103/PhysRevE.84.016301},
  url = {https://link.aps.org/doi/10.1103/PhysRevE.84.016301}
}

@article{venkataraman2017onset,
  title={The onset of thermalization in finite-dimensional equations of hydrodynamics: insights from the {Burgers} equation},
  author={Venkataraman, Divya and Sankar Ray, Samriddhi},
  journal={Proceedings of the Royal Society A: Mathematical, Physical and Engineering Sciences},
  volume={473},
  number={2197},
  pages={20160585},
  year={2017},
  publisher={The Royal Society Publishing}
}

@article{di2018dynamics,
  title={Dynamics of partially thermalized solutions of the {Burgers} equation},
  author={Di Leoni, Patricio Clark and Mininni, Pablo D and Brachet, Marc E},
  journal={Physical Review Fluids},
  volume={3},
  number={1},
  pages={014603},
  year={2018},
  publisher={APS}
}

@article{frisch2008hyperviscosity,
  title={Hyperviscosity, Galerkin truncation, and bottlenecks in turbulence},
  author={Frisch, Uriel and Kurien, Susan and Pandit, Rahul and Pauls, Walter and Ray, Samriddhi Sankar and Wirth, Achim and Zhu, Jian-Zhou},
  journal={Physical review letters},
  volume={101},
  number={14},
  pages={144501},
  year={2008},
  publisher={APS}
}

@article{choi2017finite,
  title={On the Finite-Time Blowup of a One-Dimensional Model for the Three-Dimensional Axisymmetric Euler Equations},
  author={Choi, Kyudong and Hou, Thomas Y and Kiselev, Alexander and Luo, Guo and Sverak, Vladimir and Yao, Yao},
  journal={Communications on Pure and Applied Mathematics},
  volume={70},
  number={11},
  pages={2218--2243},
  year={2017},
  publisher={Wiley Online Library}
}

@article{murugan2020suppressing,
  title={Suppressing thermalization and constructing weak solutions in truncated inviscid equations of hydrodynamics: Lessons from the Burgers equation},
  author={Murugan, Sugan Durai and Frisch, Uriel and Nazarenko, Sergey and Besse, Nicolas and Ray, Samriddhi Sankar},
  journal={Physical Review Research},
  volume={2},
  number={3},
  pages={033202},
  year={2020},
  publisher={APS}
}

@book{trefethen2019approximation,
  title={Approximation Theory and Approximation Practice, Extended Edition},
  author={Trefethen, Lloyd N},
  year={2019},
  publisher={SIAM}
}

@book{whitham2011linear,
  title={Linear and nonlinear waves},
  author={Whitham, Gerald Beresford},
  year={2011},
  publisher={John Wiley \& Sons}
}

@article{shu1988efficient,
  title={Efficient implementation of essentially non-oscillatory shock-capturing schemes},
  author={Shu, Chi-Wang and Osher, Stanley},
  journal={Journal of computational physics},
  volume={77},
  number={2},
  pages={439--471},
  year={1988},
  publisher={Elsevier}
}

@article{cartes2021galerkin,
  title={The {Galerkin}-truncated {Burgers} equation: {Crossover} from inviscid-thermalised to {Kardar}-{Parisi}-{Zhang} scaling},
  author={Cartes, C and Tirapegui, E and Pandit, R and Brachet, M},
  journal={arXiv preprint arXiv:2105.06170},
  year={2021}
}

@book{carrier2005functions,
  title={Functions of a complex variable: theory and technique},
  author={Carrier, George F and Krook, Max and Pearson, Carl E},
  year={2005},
  publisher={SIAM}
}

@article{barkley2020fluid,
  title={A fluid mechanic’s analysis of the teacup singularity},
  author={Barkley, Dwight},
  journal={Proceedings of the Royal Society A},
  volume={476},
  number={2240},
  pages={20200348},
  year={2020},
  publisher={The Royal Society Publishing}
}

\appendix

\section{Filters}
\label{app:filters}
Here, we give functional forms of the kernels that we use in the manuscript. We first summarise the positive kernels and then kernels that do not satisfy this property. Then, we briefly detail the implementation of the Spectrally Vanishing Viscosity (SVV) method for the 1D inviscid Burgers equation.

In real space, a kernel is denoted by $K_m(x)$. We then define the convolution operator $\mathcal{M}_m(\cdot) = K_m * (\cdot)$ used in Eqs.~\eqref{eq:burg_k_sr} and~\eqref{eq:burg_k_sp}. Note that the convolution operation, in real space, simplifies to a product in Fourier spectral space. For a $N$-point Fourier pseudospectral discretisation, $K_m(x) = \sum_{|k|\le N} \widehat{K}_m(k) e^{i\frac{2\pi}{L}kx}$, where $\{\widehat{K}_m(k)\}_{|k|\leq N}$ are the Fourier series coefficients of the kernel computed according to Eq.~\eqref{eq:fourcoeffs}. 

\subsection{Positive kernels}
For the following kernels, we use the parameter set $(\alpha,\gamma)$ given in the definitions of $m(N)$ and $\tau(N)$ in Eqs.~\eqref{eq:sr_params}:
\begin{align}
        m  = N^\gamma ,  \quad 0 < \gamma <1 ,\qquad \tau  = N^{-\alpha},  \quad {\alpha >0}.
\end{align}
These kernels $K_m(x)$ are positive in real space. Their analytical expressions and further discussion can be found in~\cite{butzer1971fourier}. Here, we only give the Fourier coefficients  $\{\widehat{K}_m(k)\}_{|k|\leq N}$.

\subsubsection{Fej\'er--Korovkin kernel}
\begin{subequations}
\begin{equation}
\widehat{K}_m(k)= \begin{cases}\left(1-\frac{|k|}{m+2}\right) \cos \left(\frac{|k| \pi}{m+2}\right)+\frac{1}{m+2} \cot \left(\frac{\pi}{m+2}\right) \sin \left(\frac{|k| \pi}{m+2}\right) & \text { if }|k| \leq m, \\ 0 & \text { if }|k|>m .\end{cases}
\end{equation}
\end{subequations}

\subsubsection{Jackson kernel}
\begin{subequations}


\begin{equation}
\widehat{K}_{2 m-2}(k)= \begin{cases}\frac{1}{2 m\left(2 m^2+1\right)}\left(3|k|^3-6 m|k|^2-3|k|+4 m^3+2 m\right) & \text { if }|k| \leq m, \\ \frac{1}{2 m\left(2 m^2+1\right)}\left(-|k|^3+6 m|k|^2-\left(12 m^2-1\right)|k|+8 m^3-2 m\right) & \text { if } m<|k| \leq 2 m-2, \\ 0 & \text { if }|k|>2 m-2 .\end{cases}
\end{equation}
\end{subequations}

\subsubsection{Jackson--de La Vall\'ee Poussin kernel}
\begin{subequations}
    

\begin{equation}    
\widehat{K}_{2 m-1}(k)= \begin{cases}1-\frac{3}{2}\left|\frac{k}{m}\right|^2+\frac{3}{4}\left|\frac{k}{m}\right|^3 & \text { if }|k| \leq m, \\ \frac{1}{4}\left(2-\left|\frac{k}{m}\right|\right)^3 & \text { if } m<|k| \leq 2 m-1, \\ 0 & \text { if }|k|>2 m-1 .\end{cases}
\end{equation}
\end{subequations}

\begin{figure}[h!]
 \centering
 \begin{tikzpicture}
\draw (0,0) node[inner sep=0]{\includegraphics[width=0.45\linewidth]{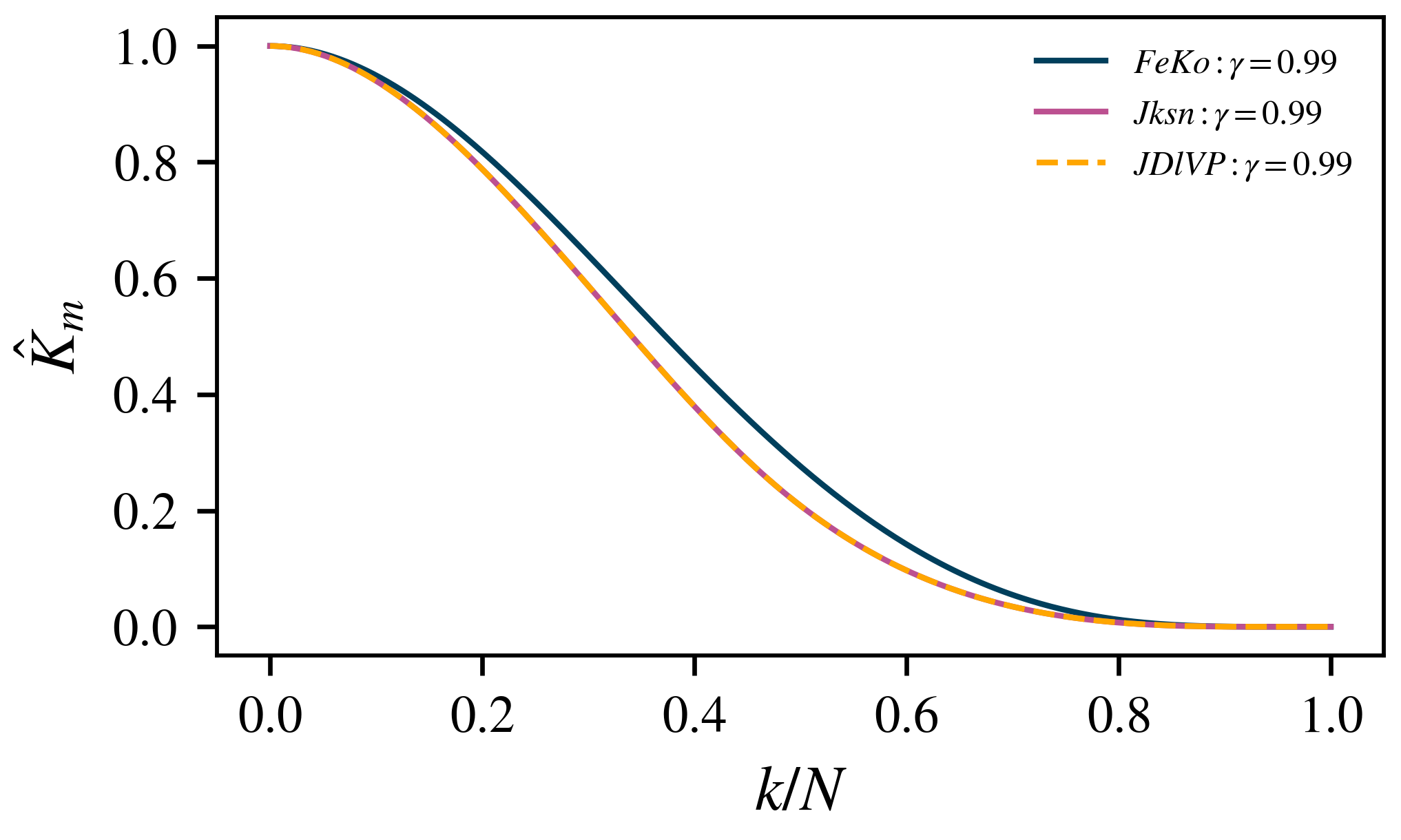}};
    \draw (2.8,-1.2) node {$(a)$};
\end{tikzpicture}
\begin{tikzpicture}
\draw (0,0) node[inner sep=0]{\includegraphics[width=0.45\linewidth]{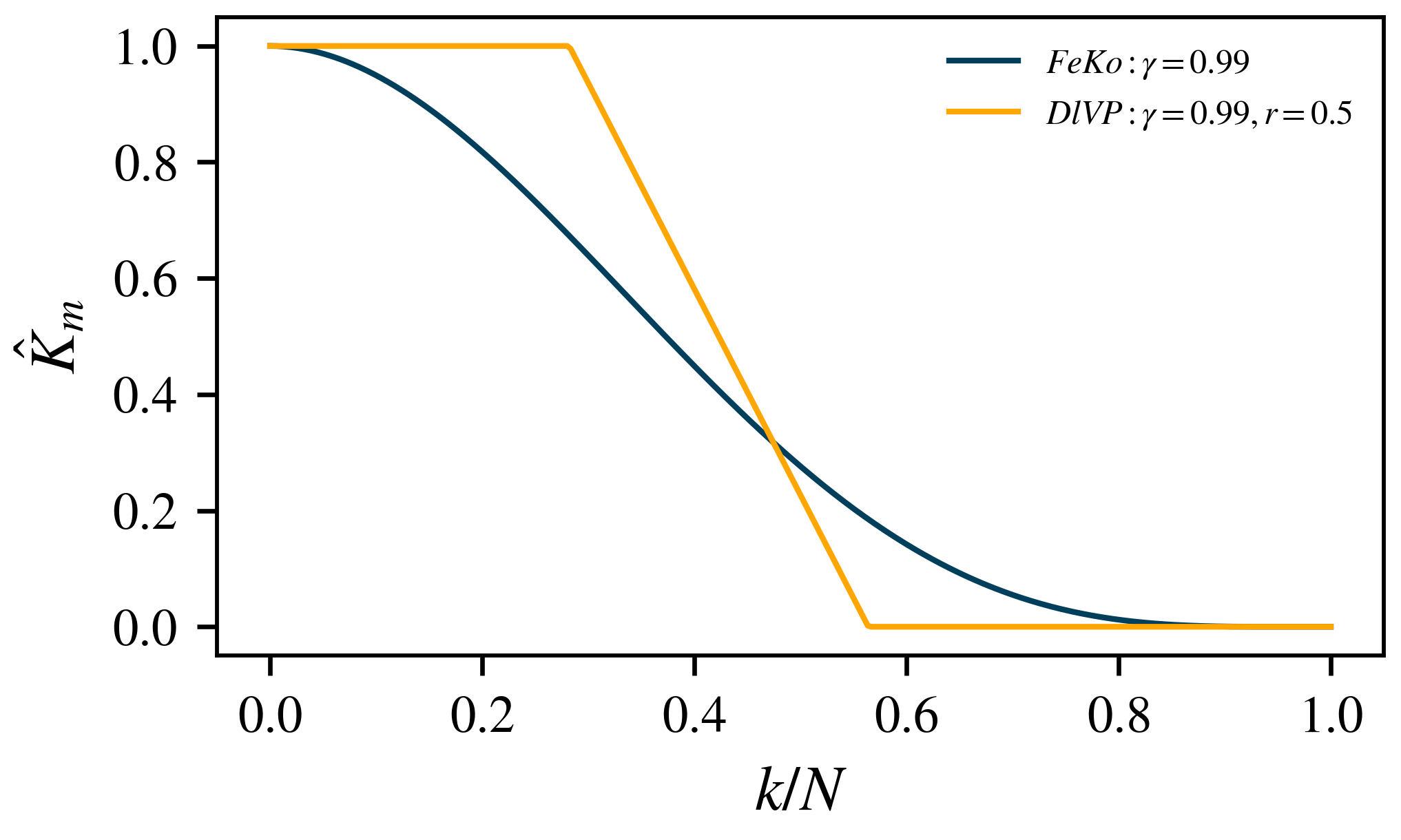}};
    \draw (2.8,-1.2) node {$(b)$};
\end{tikzpicture}
 \caption{Plots vs. $|k|/N$ of Fourier coefficients $\widehat{K}_m(k)$ of the kernels $(a)$ Fej\'er--Korovkin (FeKo), Jackson (Jksn) and Jackson--de La Vall\'ee Poussin (JDlVP) for $\gamma=0.99$. $(b)$ Analogous plots for the de La Vall\'ee Poussin (DlVP) kernel. We also show the profile of Fej\'er--Korovkin for comparison.}
 \label{fig:app_filters1}
\end{figure}

\subsection{Other kernels}
Here, we would like to summarise the other kernels that we use in the convolution operator in Fig.~\ref{fig:SR-diffkernels}. These kernels are not positive. We use the parameter $\alpha$ as in the definition for the relaxation time scale $\tau = N^{-\alpha}$. Our use of $\gamma$ and other parameters change according to the attributes of the kernel and are discussed below.

\begin{subequations}    
\subsubsection{de La Vall\'ee Poussin kernel (DlVP)}
The de La Vall\'ee Poussin kernel function with parameters $(n,p)$ is given as
\begin{align}
    \widehat{K}_m(k) = 
    \begin{cases}
          1         & \text{if } |k| \le n,  \\
  \frac{n+p-|k|}{p} & \text{if } n \le |k| \le n+p, \\
          0         & \text{if } |k| \ge n+p.  \\
    \end{cases}
\end{align}
\end{subequations}

Here, we have introduced a parameter $r$ such that $n=r {N}^\gamma$ and $p=(1-r)N^\gamma$. For the studies presented in Sections~\ref{subsec:Results_burg}-\ref{subsec:Results_Euler}, we use $r=0.5$. For the 1D HL model in Section~\ref{subsec:Results_lhmodel}, we use $r=0.9$. This kernel is a high-order approximation kernel compared to the ones defined previously.

\subsubsection{Adaptive spectral filter from Tadmor et al. (2005) (TT05)}
The adaptive spectral filter of~\cite{tadmor2005adaptive,tadmor2007filters} is given by 
\begin{subequations}
\begin{align}
    \widehat{K}_p(k) = \begin{cases}
          \exp \left( \frac{{(|k|/N)}^p}{{(|k|/N)}^2 - 1} \right)  & \text{if } |k| \le N,  \\
          0  & \text{if } |k| \ge N. 
    \end{cases} 
\end{align}
\end{subequations}

The order $p = N^{\gamma}$ and $\tau=N^{-\alpha}$. This kernel does not recover root-exponential accuracy since we do not implement an adaptive choice of the order $p_N(x)$. The non-adaptive choice of $p \sim N^\gamma$, which we use, is not dependent on the evaluation of $d_x$ and was shown to lead to a loss of convergence near the discontinuity. The optimal parameter set that was used in Fig.~\ref{fig:SR-diffkernels} is $(\alpha=-0.95,\gamma=0.1)$.

\subsubsection{Spectral filter from Majda et al. (1978) (MMO78)}
The smoothing kernel used by~\cite{majda1978fourier} are partitions of unity of the following form,
\begin{align}
    \widehat{K}_m(k) = \begin{cases}
          1  & \text{if } |k| \le m(N),  \\
  e^{-\xi(k-m)^{2p}} & \text{if } m(N) \le |k| \le N.
    \end{cases}
\end{align}

Here, $m(N) = N^{\gamma}$ and $\xi=10^{-\beta}$. The optimal parameter set that was used in Fig.~\ref{fig:SR-diffkernels} is $(\alpha=0.87,\gamma=0.87,\beta=2.5,p=1)$.


\subsubsection{Gaussian Regularised Shannon Kernel of Sun et al. (2002)}
The kernel used in the work of~\cite{gu2003conjugate,liwen2004regularised} is the Gaussian regularised Shannon Kernel. Its Fourier coefficients are given by
\begin{align}
    \widehat{\delta}_{\sigma,\Delta}(k) = \frac{1}{2} 
   \left[\text{erf} \left( \frac{\sigma}{\sqrt{2}} \left(\frac{\pi}{\Delta}
   - k\right) \right)+\text{erf} \left( \frac{\sigma}{\sqrt{2}} \left(\frac{\pi}{\Delta}
   + k\right) \right)\right] \qquad \text{for $|k|\le N$.}
\end{align}
where $\Delta = N^{-\gamma}$ is comparable to the grid spacing and $\text{sinc} (t) = \sin(\pi t)/(\pi t)$ is the normalised sinc function, which is defined at the origin as, $\text{sinc}(0)=1$. Furthermore, $\sigma=r \Delta$, where we can choose the order of the filter by the choice of $r = 1 ,1.5 ,2, \dots$.

The optimal parameter set that was used in Fig.~\ref{fig:SR-diffkernels} is $(\alpha=0.85,\gamma=0.6,r=1.5)$.


\begin{figure}[h!]
 \centering
 \begin{tikzpicture}
\draw (0,0) node[inner sep=0]{\includegraphics[width=0.45\linewidth]{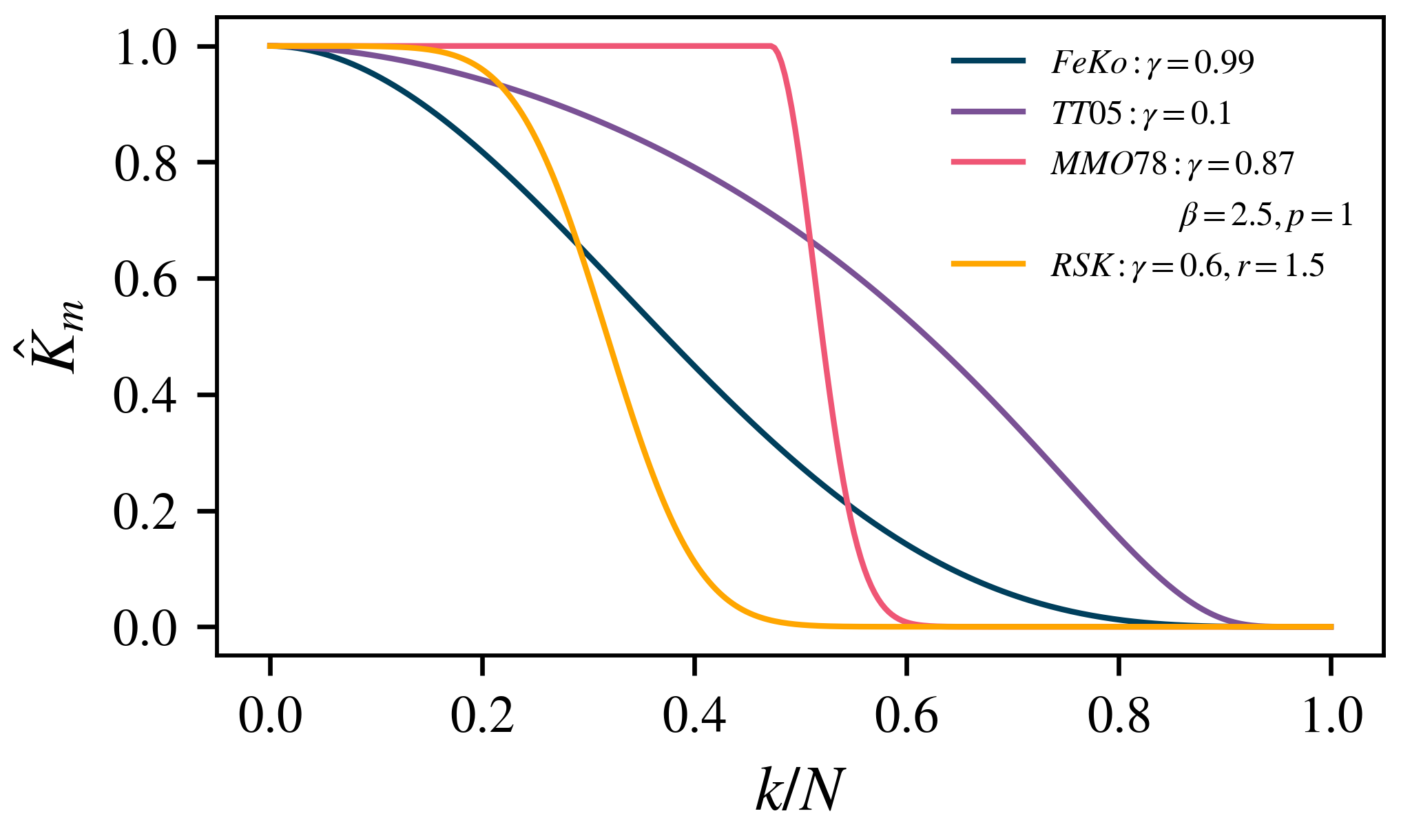}};
    \draw (2.8,-0.9) node {$(a)$};
\end{tikzpicture}
\begin{tikzpicture}
\draw (0,0) node[inner sep=0]{\includegraphics[width=0.45\linewidth]{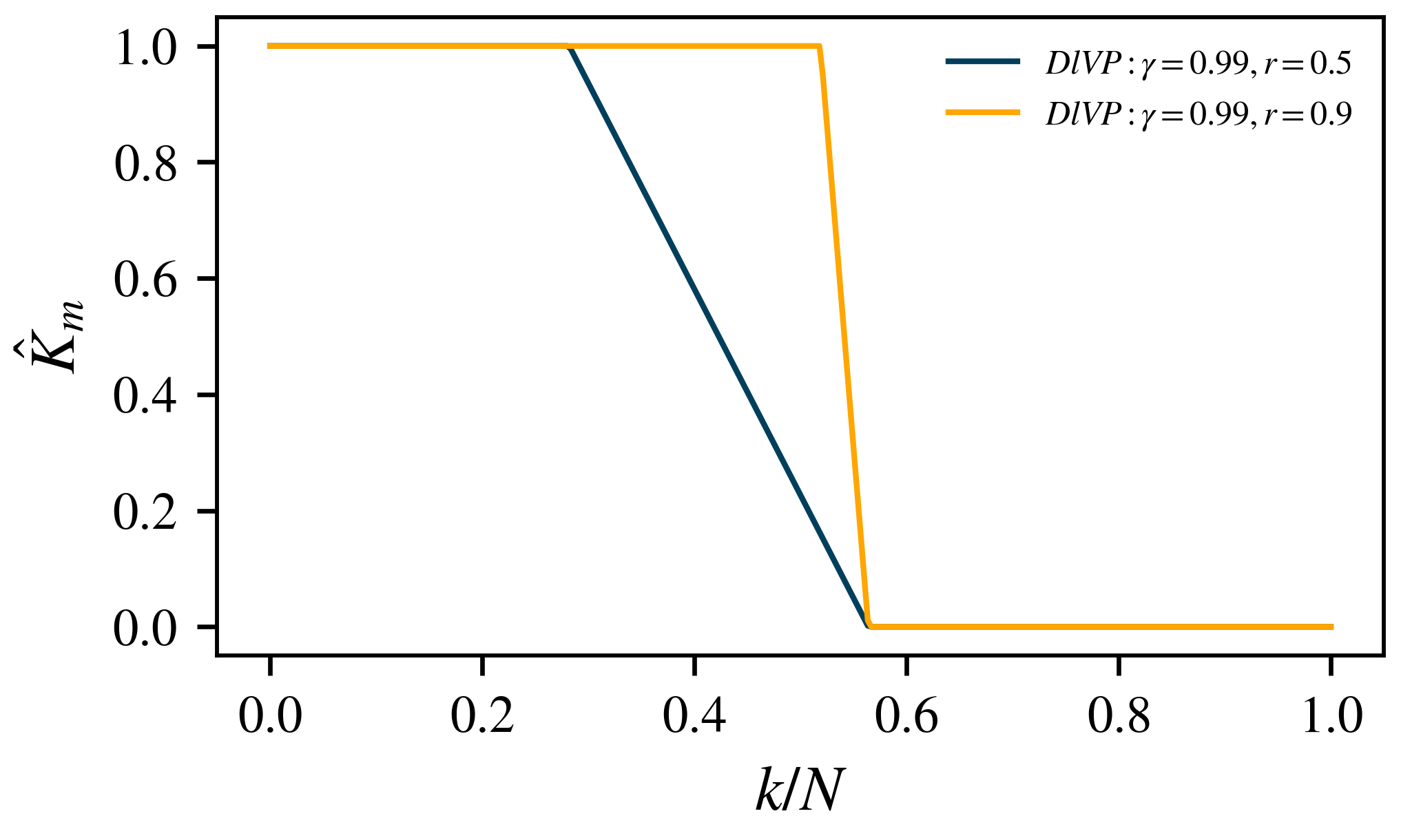}};
    \draw (2.8,-0.9) node {$(b)$};
\end{tikzpicture}
 \caption{ Plots vs. $|k|/N$ of the Fourier coefficients $\widehat{K}_m(k)$ of the kernels $(a)$ Tadmor et al.'s adaptive spectral filter (TT05), Spectral filter of Majda et al. (MMO58) and the Gaussian Regularised Shannon kernel of Sun et al. (RSK) for the parameters listed in the legend. These kernels are used for the results shown in Fig.~\ref{fig:SR-diffkernels}. $(b)$ Analogous plots for the de La Vall\'ee Poussin kernel used for the 1D HL model (orange) in Fig.~\ref{fig:axeul} and the rest of the models (blue) of Sections~\ref{subsec:Results_burg}-~\ref{subsec:Results_Euler}. We observe that the higher value of $r=0.9$, allows us to retain unmodified modes for a larger range of wavenumbers. The dissipation acts in the sloped area which has an extent of $p$ wavenumbers.}
 \label{fig:app_filters2}
\end{figure}

\subsection{The Spectrally Vanishing Viscosity (SVV) method}
The 1D inviscid Burgers equation~\eqref{eq:burg} is modified by including viscous regularisation~\cite{tadmor1989convergence} on the RHS, giving:
\begin{subequations}
    \begin{align}
    \partial_t u_N + \partial_x \left( \mathcal{P}_N \tfrac{1}{2} u_N^2 \right) = \epsilon \partial_x \left( Q_N * \partial_x u_N \right).
\end{align}
Here, $u_N(x,t)$ is the Fourier pseudospectral approximation of $u(x,t)$ on a uniform grid of $N$ points. Periodic boundary conditions are used on the domain $x \in [0,L]$, where $L=1$ for the IBVP discussed in Sections~\ref{subsec:Results_burg} and~\ref{subsec:Results_burg_ics}.

We can then write the above equation in Fourier spectral space, where the convolution operation in real space becomes a product in spectral space, as follows:
\begin{align}
\partial_t \widehat{u}(k,t) + i \left( \frac{2\pi}{L}k \right) \reallywidehat{P_N (\tfrac{1}{2}u_N^2)}(k,t) = - \epsilon  \left( \frac{2\pi}{L}k \right)^2 \widehat{Q}(k) \widehat{u}(k,t) ,
\end{align}
where $\tfrac{2\pi}{L}k$ is the Fourier wavenumber. 
Reference~\cite{maday1989analysis} uses the following $C^{\infty}$ kernel as it improves the SVV method:
\begin{align}
    \widehat{Q}(k)  = \exp \left(- \frac{(k-N)^2}{(k-M)^2} \right), \qquad k>M,
\end{align}
where $\epsilon = N^{-1}$ and $M = 2 \sqrt{N}$ for the results presented in Fig.~\ref{fig:SRSPSVV}. Different parameters were used in Fig.~\ref{fig:svv-deal-conv}; they are listed in the legend. 
\end{subequations}

\section{Implementation of characteristic boundary conditions for 1D compressible Euler equations in Chebyshev spectral methods}
\label{app:cbc}
We use the formalism established in~\cite{thompson1990time} and briefly summarise characteristic boundary treatment for the 1D compressible Euler equation. We then discuss two sets of boundary conditions: $(a)$ reflecting solid walls at both boundaries and $(b)$ supersonic inflow at $x_L$ and outflow at $x_R$, which we use in Section~\ref{subsec:Results_Euler}. Note that~\cite{thompson1990time} considers the case of the 3D compressible Euler equation; we can use the same treatment for the 1D case by simply neglecting the terms corresponding to the $x_2$ and $x_3$ coordinates. Furthermore, we do not consider the effect of gravitation on the fluid and thus we can set $g=0$ to derive the boundary conditions for our case. 

The 1D compressible Euler equations are given in conserved form in Eqs.~\eqref{eq:1deulereq}. We define $\mathbf{\Tilde{U}}(\mathbf{x},t) = (\rho,E, \rho u)$ as the vector of conserved quantities. We reformulate the system in terms of primitive variables $\mathbf{U}=(\rho,p,u)$ as this simplifies the boundary treatment. Now, Eqs.~\eqref{eq:1deulereq} can be written as 
    \begin{align}
\mathbf{U}_t + \mathbf{A} \mathbf{U}_x = 0\,,
\label{eq:1deulereq_ql}
    \end{align}
    
where the Jacobian matrix $\mathbf{A}$ is a $3\times3$ matrix is given by
    \begin{align}
\mathbf{A}= \left( \begin{array}{ccc}
    u & 0 & \rho \\
    0 & u & \gamma p \\
    0 & 1/\rho & u 
\end{array} \right).
    \end{align}
Since the above system of equations is hyperbolic, $\mathbf{A}$ is diagonalisable and has real eigenvalues. 
\subsection{The Characteristic Analysis}
The boundary treatment of hyperbolic PDEs is greatly simplified by using the characteristic variable formulation. In this formulation, the PDEs represent the propagation of waves of characteristic variables $z_i$ at their respective velocities $\lambda_i$; here, $z_i$ and corresponding $\lambda_i$ are functions of $\mathbf{U}(x,t)$. Depending on the sign of the $\lambda_i$ at the boundaries, the associated characteristic lines point into or out of the domain. The characteristics that point out of the domain at a given boundary are defined entirely by the solution at or within the boundary, and we are not required to specify boundary conditions for these waves. However, the inward pointing characteristics represent waves that are incoming and require boundary conditions to specify their behaviour completely. We thus require a boundary treatment for these characteristic variables. We begin by solving the eigenvalue problem for the Jacobian matrix $\mathbf{A}$ which gives us the eigenvalues/characteristic velocities $\lambda_i$ and the eigenvectors which are used to transform primitive variables to characteristic variables.   

The eigenvalues of the matrix $\mathbf{A}$ are
\begin{align}
    \lambda_1 = u-c, \quad \lambda_2 = u,\quad \lambda_3 = u+c, \label{eq:eigvals}
\end{align}
where $c=\sqrt{\frac{\gamma p}{\rho}}$ is the speed of sound. 
The corresponding left eigenvectors $\mathbf{l}_i$ of  $\mathbf{A}$ are given by
\begin{subequations}    
\begin{align}
    \mathbf{l}^T_1 &= (0,1,-\rho c), \\
    \mathbf{l}^T_2 &= (c^2,1,0), \\
    \mathbf{l}^T_3 &= (0,1,\rho c). 
\end{align}
\end{subequations}

Diagonalisation for $\mathbf{A}$ is obtained by using the matrices $\mathbf{S}$ where the columns are the right eigenvectors $\{ \mathbf{r}_j\}_{j=1,2,3}$ and $\mathbf{S}^{-1}$ where the rows are the left eigenvectors $\{\mathbf{l}_i\}_{i=1,2,3}$. The diagonalisation is then given by
\begin{align}
    \mathbf{S}^{-1} \mathbf{A} \mathbf{S} = \mathbf{\Lambda},
\end{align}
where $\mathbf{\Lambda}$ is the diagonal matrix of eigenvalues with $\Lambda_{ij}=\lambda_j \delta_{ij}$. Applying this transformation to Eq.~\eqref{eq:1deulereq_ql}, we obtain
\begin{align}
    \mathbf{S}^{-1} \mathbf{U}_t + \mathbf{\Lambda}\mathbf{S}^{-1} \mathbf{U}_x =0
    \label{eq:1deul_ql1}
\end{align}
whose components are
\begin{align}
    \mathbf{l}^T_i \mathbf{U}_t + \lambda_i  \mathbf{l}^T_i \mathbf{U}_x =0 , \qquad i=1,2,3.
\end{align}

We define a vector $\mathcal{L}$ for the boundary analysis given by
\begin{align}
    \mathcal{L}_i \equiv \lambda_i \mathbf{l}^T_i \mathbf{U}_x.
\end{align}

For the 1D compressible Euler equations, the components $\mathcal{L}_i$ are given by
\begin{subequations}
\begin{align}
    \mathcal{L}_1 &= \lambda_1 (\partial_x p - \rho c \ \partial_x u ), \\ 
        \mathcal{L}_2 &=  \lambda_2 ( c^2 \partial_x \rho  - \partial_x p ), \\
        \mathcal{L}_3 &= \lambda_3 (\partial_x p +    \rho c \ \partial_x u). 
\end{align}    
\label{eq:Lcomp}
\end{subequations}
Here, we remark that we use $\mathcal{L}_3$ where~\cite{thompson1990time} uses $\mathcal{L}_5$ and so on. The problem of implementing boundary conditions is reduced to the problem of computing the appropriate values for $\mathcal{L}_i$. We later discuss the specifics of the computation of $\mathcal{L}_i$ in the domain and at boundaries, and expand on the role played by boundary conditions in these specifications; but we first lay down the boundary treatment, once the values of $\mathcal{L}_i$ are computed and known. 
The primitive formulation in Eq.~\eqref{eq:1deul_ql1} can be written in terms of $\mathcal{L}_i$ as follows,
\begin{align}
   \mathbf{U}_t + \mathbf{S} \mathcal{L} = 0.
   \label{eq:1deuler_ql2}
\end{align}
Here, we define $\mathbf{d}\equiv \mathbf{S}\mathcal{L}$ whose components are given by
 \begin{align}
    \mathbf{d} &= \begin{pmatrix}
           d_{1} \\
           d_{2} \\
           d_{3}
         \end{pmatrix} 
         = \begin{pmatrix}
           (1/c^2)(\mathcal{L}_2 -\frac{1}{2}(\mathcal{L}_3 + \mathcal{L}_1)) \\
           (1/2)(\mathcal{L}_3 + \mathcal{L}_1) \\
           (1/2\rho c)(\mathcal{L}_3 - \mathcal{L}_1)
         \end{pmatrix}.
         \label{eq:dcomp}
  \end{align}

Since we work with the conservative form of the 1D compressible Euler equations in Section~\ref{subsec:Results_Euler}, we use the transformation $\mathbf{P}_{ij} = \partial{\Tilde{\mathbf{U}}_i}/\partial{\mathbf{U}_j} $. Then, the conservative formulation for Eqs.~\eqref{eq:1deuler_ql2} in terms of $\mathbf{d}$ is given by
\begin{align}
   \Tilde{\mathbf{U}}_t + \mathbf{Pd} = 0.
\end{align}
Equations ~\eqref{eq:1deulereq} are now rewritten in terms of the vector $\mathcal{L}_i$ for simple boundary treatment as follows,
\begin{subequations}    
\begin{align}
    \partial_t \rho + d_1 &= 0, \\
        \partial_t E + \frac{1}{2} u^2 d_1 + \frac{d_2}{\gamma -1} + \rho u d_3 &=0, \\
    \partial_t (\rho u) + u d_1 + \rho d_3   &=0. \label{eq:1deul_u_bt}
\end{align}
\label{eq:1deul_bt}
\end{subequations}

\subsection{Boundary treatment}
We now describe the procedure for calculation of $\mathcal{L}_i$ and their dependence on boundary conditions. 
\begin{itemize}
    \item For all points in the interior of the domain $x_L < x < x_R$, we compute $\mathcal{L}_i$ using Eqs.~\eqref{eq:Lcomp} and then compute components of $\mathbf{d}$ using Eq.~\eqref{eq:dcomp}. We then use these values to solve for $\Tilde{\mathbf{U}}$ from Eqs.~\eqref{eq:1deul_bt}.
    
    \item  For the left (right) boundary at $x=x_L(x_R)$, we first compute the characteristic velocities $\{\lambda_i,i=1,2,3\}$ using Eq.~\eqref{eq:eigvals}. When the velocity points out of the domain, i.e., for the left boundary, $\lambda_i(x_L) \le 0$ and for the right boundary $\lambda_i(x_R) \ge 0$, we compute the corresponding $\mathcal{L}_i$ using Eq.~\eqref{eq:Lcomp} since the associated outgoing characteristic wave is completely determined by the solution in the solution domain. 
    When the velocity points into the domain, i.e., $\lambda_i(x_L)>0$ or $\lambda_i(x_R)<0$, we compute the corresponding $\mathcal{L}_i$ using the boundary conditions since the associated characteristic wave is incoming and its specification depends on the solution near the boundary as well as the boundary conditions. We discuss the specifics below for 2 sets of boundary conditions that we use in the manuscript:
\begin{itemize}
    \item[\textbf{(a)}] \textbf{Reflecting solid walls at $x_L$ and $x_R$:} At a solid wall, the velocity $u=0$ at all times. The initial data reflects this boundary condition, as seen in the case of the shock tubes and the blast-waves problems in Section~\ref{subsec:Results_Euler}. 
    
    At the left boundary $x_L$, $\lambda_1=-c<0$ is associated with an outgoing wave and hence $\mathcal{L}_1$ must be computed from the definition in Eq.~\eqref{eq:Lcomp}. At this boundary, $\lambda_2=0$ and thus, $\mathcal{L}_2=0$. For the incoming characteristic wave associated to $\lambda_3=c>0$, we specify $\mathcal{L}_3$ by requiring that $u|_{x_L}=0$ for all time. We see that $u(x_L,t)$ will remain $0$ in Eq.~\eqref{eq:1deul_u_bt}, if we set $\mathcal{L}_3=\mathcal{L}_1$. 
    
    We have the inverse scenario at the right boundary $x_R$. Eigenvalues have the same values, but now $\lambda_3$ is associated with an outgoing wave. We thus compute $\mathcal{L}_3$ from the values of the solution at the boundary using Eq.~\eqref{eq:Lcomp}. Again, we have $\lambda_2=0$ and thus, $\mathcal{L}_2=0$. The incoming characteristic wave is now associated with $\lambda_1>0$ and has to be set by the relation $\mathcal{L}_1=\mathcal{L}_3$.
    
    \item[\textbf{(b)}]\textbf{ Supersonic outflow at $x_L$ and non-reflecting subsonic inflow at $x_R$:} Supersonic inflow at $x_L$ is characterised by $u>c$. This implies that for all $i$, $\lambda_i>0$. Thus, we are required to specify the values of all components of $\mathcal{L}$ using the boundary conditions. Since supersonic flow in the case of the Shu--Osher problem is implemented as a steady-state boundary condition, we set $\mathcal{L}=0$.

    At the right boundary $x_R$, we have subsonic outflow, and, in particular, we use a non-reflecting boundary so that the waves produced in the domain exit the domain without causing any reflections or artificial resonance. We impose this condition by requiring that $\mathcal{L}_1=0$, as described in Section 3.2.6.1 of~\cite{thompson1990time}. 
\end{itemize}

\end{itemize}


\end{document}